%% file: hopf.tex
\documentclass[%
a4paper,%
arxiv,%
]%
{myclass}

\usepackage{amsmath}

\ifLMS
\else
  \usepackage{amsthm}
\fi

\ifpxfonts
  \usepackage{pxfonts}
\else
  \usepackage{amssymb}
  \newcommand{\coloneqq}{\colon=}
\fi

\usepackage{hyphenat}
\usepackage[matrix,arrow]{xy}

\usepackage[nowarnings]{hopf}
\usepackage{hopf.abbrev}

\ifLMS
  \shiftpage{30pt}
\fi

\newtheorem{thm}{Theorem}
\ifLMS
\else
  \newtheoremstyle{thm}{3pt}{3pt}{\itshape}{}{\bfseries}{}{.5em}{}
  \newtheoremstyle{thmsub}{3pt}{3pt}{\upshape}{}{\bfseries}{}{.5em}{}
  \theoremstyle{thm}
\fi
\newtheorem{theorem}{Theorem}[section]
\newtheorem{lemma}[theorem]{Lemma}
\newtheorem{proposition}[theorem]{Proposition}
\ifLMS
  \newnumbered{defn}[theorem]{Definition}
\else
  \newtheorem{defn}[theorem]{Definition}
\fi
\newtheorem{corollary}[theorem]{Corollary}
\ifLMS
\else
  \theoremstyle{thmsub}
\fi
\ifLMS
  \newnumbered{examples}[theorem]{Examples}
\else
  \newtheorem{examples}[theorem]{Examples}
\fi

\numberwithin{equation}{section}

\begin{document}

\ifLMS
  \bibliographystyle{lms}
\else
  \bibliographystyle{alpha}
\fi

\title{%
The Hunting of the Hopf Ring
}  
\date{\today}
\author{%
Andrew Stacey
\and
Sarah Whitehouse%
}

\hopfMSC{
Primary: 55S25;
Secondary: 55N20,
16W99
}

\ifAMS
\fi

\hopfthanks{The authors acknowledge the support of the EPSRC,
  grant no.: GR/S76823/01.}

\ifAMS
\else
  \maketitle
\fi

\begin{abstract}
We provide a new algebraic description of the structure on the set of
all unstable cohomology operations for a suitable generalised
cohomology theory, \(\Efunc(-)\). Our description is as a graded
and completed version of a \emph{Tall\hyp{}Wraith monoid}. The
\(\Efunc\)\hyp{}cohomology of a space \(X\) is a module for this
Tall\hyp{}Wraith monoid.  We also show that the corresponding Hopf
ring of unstable \co{operations} is a module for the Tall\hyp{}Wraith
monoid of unstable operations. Further examples are provided by
considering operations from one theory to another.
\end{abstract}

\ifAMS
  \maketitle
\fi

\include{hopf.intro}
\include{hopf.genalg}

\include{hopf.filter}
\include{hopf.tall}

\include{hopf.algtop}

\bibliography{hopf}

\end{document}

%% file: hopf.intro.tex
\section{Introduction}
\label{sec:intro}

In this paper we provide a new algebraic description of 
the structure on the set of all unstable
cohomology operations for a suitable generalised cohomology theory,
say \(\Efunc(-)\).  The bigraded set of unstable
operations for \(\Efunc(-)\) is identified by the usual Yoneda
lemma argument with the bigraded set of cohomology groups of the
representing spaces, \(\{\Efunc[k](\Erep[l])\}\).  This has
considerable structure and it is natural to ask how to best describe
it.

To date, the most comprehensive work in this area is due to Boardman,
Johnson, and Wilson~\cite{jbdjww}.  They provide the following descriptions of
this structure.
\begin{enumerate}
\item The unstable operations of a suitable cohomology theory define a
monad on the category of complete, Hausdorff, filtered, graded,
commutative, unital \(\Ecoef\)\hyp{}algebras.
\item The unstable operations of a suitable cohomology theory are dual
to the \emph{enriched Hopf ring} of the corresponding \co{operations}.
\end{enumerate}

The \emph{Hopf ring} part --
  i.e., ignoring the \emph{enrichment} -- of the second answer
  is, of course, a well\hyp{}established and
  important notion in algebraic topology.  Since the work of Ravenel
  and Wilson~\cite{dcrwsw} the language of Hopf rings has 
  been widely used and
  the structures of the Hopf rings associated to many important
  cohomology theories have been calculated; see, for
  example~\cite{%
    dscmns,%
    dcrwsw,%
    dcrwsw96,%
    wsw84%
}.  A useful
  introduction to Hopf rings, with further references, can be found
  in~\cite{wsw00}.

In considering a Hopf ring, we have not yet taken into account one of
the most obvious pieces of structure: that operations may be composed.
A Hopf ring does not include any structure which dualises to
composition of operations.  The \emph{enrichment} in the second answer
encodes the dual of composition.  The enriched Hopf ring structures of
several important cohomology theories are described in~\cite{jbdjww}.

The first answer, describing operations as a monad on a suitable
category, certainly includes the composition as a fundamental part of
the structure.  However, this answer does not explicitly describe the
internal structure on
the set of operations; instead it specifies the
action of operations on some category.  In particular, this approach
does not lend itself to specifying the structure of operations via
generators and relations.

Of these descriptions, that of an unenriched Hopf ring has
proved to be more amenable to further study than either that of an
enriched Hopf ring or that of a monad on a suitable category of
algebras.
  
Our description of operations may be thought of as a monoidal
reinterpretation of the first answer.  It is algebraic in nature and
we employ the language of general or universal algebra to express it.
One advantage of this approach is that it describes the structure of
unstable operations, including composition, directly.  Another is that
it allows for descriptions of unstable operations via generators and
relations.  Such descriptions are expected to shed light on the
relationships between stable, additive, and unstable operations. For
the Morava K\hyp{}theories and related cohomology theories, these
relationships were studied in \cite{math.AT/0605471}; the results of
this paper provide the foundations for a very explicit description,
via generators and relations, of the splitting of stable operations
from unstable operations given in \cite{math.AT/0605471}.  We expect
to develop this point of view for familiar examples in future work.

Roughly speaking, our answer is that the unstable operations of a
suitable theory have the structure of a graded, completed
\emph{Tall\hyp{}Wraith monoid}; a term that we shall now explain.

Let \(\vcat\) be a variety of algebras, in the sense of general or
universal algebra.  A \emph{Tall\hyp{}Wraith \(\vcat\)\hyp{}monoid} is
a set with precisely the algebraic structure required for it to act on
\(\vcat\)\hyp{}algebras.  To make this precise, one considers \vCvcat.
This is equivalent to the category of representable functors from
\(\vcat\) to \(\vcat\) and so it has a monoidal structure
corresponding to composition of functors.  A Tall\hyp{}Wraith
\(\vcat\)\hyp{}monoid is then defined to be a monoid in this category.
One example is very familiar: a ring is a Tall\hyp{}Wraith
\(\vcat\)\hyp{}monoid for \(\vcat\) the category of abelian groups.
The case originally considered by Tall and Wraith in~\cite{dtgw},
under the name \emph{biring triple}, was for \(\vcat\) the category of
commutative unital rings.  Recently, Borger and Weiland~\cite{jbbw}
rediscovered this and extended it to the case where \(\vcat\) is the
category of commutative unital \(k\)\hyp{}algebras, for a commutative
unital ring \(k\).  They adopted the term \emph{plethory} in that
situation; thus a plethory is
that\hyp{}which\hyp{}acts\hyp{}on\hyp{}algebras.  This is clearly
relevant to our purposes as unstable cohomology operations of
multiplicative cohomology theories act on the cohomology algebras.

It is also clear that there remains work to be done, because the
cohomology theories that we are considering are graded and
topologised.  The grading introduces no technical difficulties:
varieties of graded algebras have been studied for almost as long as
varieties of ordinary algebras; they are also known as
\emph{many\hyp{}sorted algebras} in the literature.  We arrive at
Tall\hyp{}Wraith \(\Gvcat\)\hyp{}monoids for \(\Gvcat\) a variety of
graded algebras.  These are naturally bi\hyp{}graded, as are unstable
cohomology operations.

The main work of this paper comes in the extension of Tall\hyp{}Wraith
monoids to a suitably topologised context.  For \(\Efunc(-)\) a
multiplicative cohomology theory and \(X\) a \(\m{C W}\)\hyp{}complex,
\(\Efunc(X)\) is given the pro\hyp{}finite topology; that is, the
filtration topology for the filtration by the ideals \(\{\ker
  i^*\colon\Efunc(X)\to \Efunc(X_f)\}\) for all inclusions of finite
subcomplexes \(i\colon X_f\to X\).  We therefore develop a theory of
filtered Tall\hyp{}Wraith monoids so that our description of unstable
operations takes into account the pro\hyp{}finite filtrations.  While
there are various notions of filtered objects in a (suitable) category
in the literature, we have not been able to find a set\hyp{}up suited
to our needs.  Therefore we introduce a suitable definition of
filtered objects in a category to model the topology; in this setting
we introduce \emph{iso\hyp{}filtrations} as the generalisation of
complete, Hausdorff spaces.  This allows us to formulate the notion of
a Tall\hyp{}Wraith \(\KGvcat\)\hyp{}monoid, where \(\KGvcat\) denotes
the category of iso\hyp{}filtered objects in the variety of graded
algebras, \(\Gvcat\).

Once we have established the general theory of graded filtered
Tall\hyp{}Wraith monoids, the applications to suitable generalised
cohomology theories are straightforward.  We can now state precisely
our monoidal reformulation of the first description of unstable
operations.  We adopt standard notation related to a cohomology theory
\(\Efunc(-)\), so the representing spaces are denoted by
\(\Erep[k]\) (for \(k\in\Z\)), the corresponding homology theory
is denoted by \(\Ehom(-)\), and the coefficient ring by
\(\Ecoef\).

\begin{thm}
\label{th:twcoh}
Let \(\Efunc(-)\) be a graded, commutative, multiplicative,
cohomology theory.  Let \(\Gvcat\) be the
variety of graded, commutative, unital \(\Ecoef\)\hyp{}algebras.
Suppose that \(\Ehom(\Erep[k])\) is free as an
\(\Ecoef\)\hyp{}module for each \(k\).  Then
\(\Efunc(\Erep)\) is a Tall\hyp{}Wraith
  \(\KGvcat\)\hyp{}monoid.
\end{thm}

As noted in~\cite{jbdjww}, the freeness hypothesis of
theorem~\ref{th:twcoh} is satisfied for \(H\F_p\), \(B P\), \(M U\),
\(K(n)\), and \(K U\).

\smallskip

There is a natural notion of a module for a Tall\hyp{}Wraith monoid.  The
cohomology of spaces provides examples of modules for the Tall\hyp{}Wraith
monoid of unstable operations.

\begin{thm}
\label{th:twcohmod}
Let \(\Efunc(-)\) and \(\Gvcat\) be as in
theorem~\ref{th:twcoh}.  Let \(X\) be a topological space.  The
natural morphisms of sets
\[
  \Efunc[k](X) \to
  \Hom{\Gvcat}{\Efunc(\Erep[k])}{\Efunc(X)}
\]
make the completed cohomology of \(X\), \(\mBo{\Efunc}(X)\),
   into a module for the Tall\hyp{}Wraith
  \(\KGvcat\)\hyp{}monoid of unstable operations
  \(\Efunc(\Erep)\).
\end{thm}

It turns out that the Hopf ring of co\hyp{}operations is also a module
for the Tall\hyp{}Wraith monoid of operations, and it is this extra
structure that is encoded in the term \emph{enriched}.  Further
examples come from considering operations and \co{operations} from one
theory to another.

\begin{thm}
\label{th:twhopf}
Let \(\Efunc(-)\) and \(\Gvcat\) be as in
theorem~\ref{th:twcoh}.
Let \(\Ffunc(-)\) be another graded, commutative, multiplicative 
cohomology theory.  Suppose that each
\(\Fhom(\Erep[k])\) is free as an
\(\Fcoef\)\hyp{}module.  Then the following statements are
true.
\begin{enumerate}
\item The bigraded set \(\Ffunc(\Erep)\) 
of unstable operations \(\Efunc(-) \to \Ffunc(-)\)
has the structure of a module for 
the Tall\hyp{}Wraith \(\KGvcat\)\hyp{}monoid 
\(\Efunc(\Erep)\).

\item The Hopf ring \(\Fhom(\Erep)\) is a module for
\(\Efunc(\Erep)\).

\item Let \(X\) be a space such that \(\Fhom(X)\) is free as an
\(\Fcoef\)\hyp{}module.  Let \(\Fcoalgcat\) be \Fcoalgcat.  Then the
natural morphisms of sets
\[
  \Efunc[k](X) \to
  \Hom{\Fcoalgcat}{\Fhom(X)}{\Fhom(\Erep[k])}
\]
extend to a morphism of modules for \(\Efunc(\Erep)\).
\end{enumerate}
\end{thm}

\medskip

This paper is organised as follows.  Section~\ref{sec:genalg} covers
background material from general algebra, in both the ungraded and
graded contexts.  The main aim of this section is to establish the
necessary conditions to consider Tall\hyp{}Wraith monoids and certain
important related concepts.  Section~\ref{sec:filter} is concerned
with setting out the necessary details of the theory of filtered
objects in a category.  In this section we consider several types of
filtrations and the relationships between them.  As we shall want to
apply the constructions of general algebra in such categories we are
also concerned with establishing the categorical properties of these
categories of filtered objects.  The main technical work of this paper
is in this section and concerns functors between categories of
filtered objects.  Section~\ref{sec:filtall} brings together the work
of the preceding sections by considering Tall\hyp{}Wraith monoids in
the filtered context.  Section~\ref{sec:algtop} applies the results to
the examples arising in algebraic topology from suitable generalised
cohomology theories, thus proving the theorems stated above.

The reader more interested in the results than the method by which
they are demonstrated may prefer to read this paper in reverse.

Finally, we wish to acknowledge the work of Boardmann, Johnson, and
Wilson in understanding unstable operations.  Although this paper does
not depend on \cite{jbdjww} mathematically, it was an invaluable
resource as a guide to determine the form of our final answer.

%% file: hopf.genalg.tex
\section{General Algebra}
\label{sec:genalg}

In this section we shall expand a little on the basic constructions of
general algebra.  The results quoted in this section are all standard
results from that field.  For ungraded algebra objects in \scat, these
results can be found in any good introduction to the subject, for
example \cite{gb}.  The more general cases can be found in the wider
literature, for example in \cite{wkjmmpis}.  We record these results
here without proofs to establish notation and as a quick reference for
the rest of the paper.  For those initiates of the deeper secrets of
general algebra we mention that we are only considering algebras of
finite arity and so we can assume that our identities are specified by
finite sets.

We start by summarising the results that we need in the arena of
ungraded algebras, also known as \emph{single\hyp{}sorted} or
\emph{homogeneous} algebras.
We shall then explain how this generalises to graded algebras, also known
as \emph{many\hyp{}sorted} or \emph{heterogeneous} algebras.

\subsection{Ungraded Algebra}

\begin{defn}
Let \(\dcat\) be a category with finite products.

\begin{enumerate}
\item A \emph{type} \(\otype\) is a pair \((\abs{\otype}, n)\) where
\(\abs{\otype}\) is a set and \(n \colon \abs{\otype} \to \N_0\) is a
morphism of sets called the \emph{arity} morphism.

\item An \emph{\doobj[\doobj]} consists of an \dobj, \(\abs{\doobj}\),
together with, for each \(\oop \in \sabs{\otype}\), a
\(\dcat\)\hyp{}morphism \(\oop_{\doobj} \colon \abs{\doobj}^{n(\oop)}
\to \abs{\doobj}\); these morphisms are called the \emph{operations}
of the \doobj.  A morphism of \doobjs is a morphism of the underlying
\dobjs which intertwines the operations.

\item An \emph{\oalg{}} is an \soobj.

\item We denote \docat by \(\docat\) and \socat by \(\ocat\).  We
refer to the functor \(\docat \to \dcat\) which assigns to an \doobj
the underlying \dobj as the \emph{forgetful functor}.  We write the
underlying \dobj of an \doobj[\doobj] as \(\abs{\doobj}\).
\end{enumerate}
\end{defn}

We trust to context to distinguish between \(\otype\) the type and
\(\ocat\) the resulting category of algebras.

One of the key initial results in general algebra is the following
characterisation of \oalgobj structures.

\begin{proposition}
\label{prop:conlift}
To give an \dobj[\abs{\doobj}] the structure of an
\oalgobj is equivalent to giving a lift of the
contravariant hom\hyp{}functor \(\Hom{\dcat}{-}{\abs{\doobj}} \colon
\dcat \to \scat\) to a functor \(\dcat \to \ocat\).  \noproof
\end{proposition}

If \(\doobj\) is an \doobj and \(\dobj\) is an \dobj we shall
write
\[
  \Hom{\dcat}{\dobj}{\doobj}
\]
for the corresponding \oalg with underlying set
\(\Hom{\dcat}{\dobj}{\abs{\doobj}}\).  The operations on
\(\Hom{\dcat}{\dobj}{\doobj}\) are straightforward; let \(\oop \in
\abs{\otype}\) be an element of arity \(n\) and let \(\oop_{\doobj}
\colon \abs{\doobj}^n \to \abs{\doobj}\) be the operation on \(V\).
Then the corresponding operation on \(\Hom{\dcat}{\dobj}{\doobj}\) is
the \(\dcat\)\hyp{}morphism
\[
\dobj \xrightarrow{\Delta^n} \dobj^n \xrightarrow{f_1 \times \dotsb
  \times f_n} \sabs{\doobj}^n \xrightarrow{\oop_{\doobj}} \sabs{\doobj}
\]
where \(\Delta^n \colon \dobj \to \dobj^n\) is the \(n\)\hyp{}fold
diagonal morphism.

\oalgs and \oalgobjs simply have operations, they are not
constrained to satisfy any particular identities.  To consider
identities, we need to know about free \oalgs.

\begin{proposition}
\label{prop:free}
Let \(\dcat\) be a category with finite products which is closed under
filtered \co{}limits and such that finite products commute with
filtered \co{}limits.  Then the forgetful functor \(\docat \to \dcat\)
has a left adjoint, \(\free{o} \colon \dcat \to \docat\), which is
called the \free{o}. \noproof
\end{proposition}

For identities, we only need to know about free \oalgs in \scat.

\begin{defn}
Let \(\otype\) be a type.  An \emph{identity} for \oalgs is a triple
\((\sobj,u,v)\) where \(\sobj\) is a (finite) set and \(u, v \in
\abs{\free{o}(\sobj)}\).
\end{defn}

Let \(\doobj\) be an \doobj.  An identity for \oalgs, \((\sobj,u,v)\),
induces natural \dmors
\[
  u_{\doobj}, v_{\doobj} \colon \abs{\doobj}^{\sobj} \to \abs{\doobj}.
\]
These are defined as follows.

The canonical projections \(\abs{\doobj}^{\sobj} \to \abs{\doobj}\) in
\(\dcat\) define a set morphism
\[
  \sobj \to \Hom{\dcat}{\abs{\doobj}^{\sobj}}{\abs{\doobj}}.
\]
As \(\doobj\) is an \doobj, the right hand side of this is
the underlying set of an \oalg.  Using the
adjunction we therefore have a morphism of \oalgs
\[
  \free{o}(\sobj) \to \Hom{\dcat}{\abs{\doobj}^{\sobj}}{\doobj}
\]
and thus the elements \(u, v \in \abs{\free{o}(\sobj)}\) define
elements in the underlying set of the \oalg
\(\Hom{\dcat}{\abs{\doobj}^{\sobj}}{\doobj}\) which is
\(\Hom{\dcat}{\abs{\doobj}^{\sobj}}{\abs{\doobj}}\).  Thus we have the
required \(\dcat\)\hyp{}morphisms \(u_{\doobj}, v_{\doobj} \colon
\abs{\doobj}^{\sobj} \to \abs{\doobj}\).

\begin{defn}
An \doobj[\doobj] is said to \emph{satisfy} the identity
\((\sobj,u,v)\) if the two induced \(\dcat\)\hyp{}morphisms
\(u_{\doobj}, v_{\doobj} \colon \abs{\doobj}^{\sobj} \to
\abs{\doobj}\) are the same.
\end{defn}

\begin{defn}
A \emph{variety of algebras}, \(\vcat\), is specified by a type,
\(\otype\), and a set of identities for \oalgs, \(J\).  It is the full
subcategory of \(\ocat\) consisting of all \oalgs which satisfy all of
the identities in \(J\).

The pair \((\otype, J)\) is a \emph{presentation} of \(\vcat\).

Let \(\dcat\) be a category with finite products, \(\vcat\) a variety
of algebras with presentation \((\otype, J)\).  \dvcatu, \(\dvcat\),
is the full subcategory of \(\docat\) consisting of all \doobjs which
satisfy all of the identities in \(J\).
\end{defn}

As is well\hyp{}known, presentations are not unique.

Proposition~\ref{prop:conlift} holds in the presence of identities.

\begin{proposition}
\label{prop:conliftid}
To give an \dobj[\abs{\doobj}] the structure of a
\valgobj is equivalent to giving a lift of the
contravariant hom\hyp{}functor \(\Hom{\dcat}{-}{\abs{\doobj}} \colon
\dcat \to \scat\) to a functor \(\dcat \to \vcat\).  \noproof
\end{proposition}

To get free \valgobjs we need to know that we can ``impose''
identities on an \oalgobj.

\begin{theorem}
\label{th:impid}
Let \(\dcat\) be a complete category with finite products.  Suppose
that \(\dcat\) is an (\(\m{E}\), \(\m{M}\)) category for some classes
\(\m{E}\) of epimorphisms and \(\m{M}\) of monomorphisms, that
\(\m{E}\) is closed under taking finite products, and that \(\dcat\)
is \(\m{E}\)\hyp{}\co{well}\hyp{}powered.  Then the inclusion functor
\(\dvcat \to \docat\) has a left adjoint, \(\docat \to \dvcat\),
called \emph{imposition of identities}. \noproof
\end{theorem}

\begin{corollary}
\label{cor:free}
Under the conditions of theorem~\ref{th:impid} and
proposition~\ref{prop:free}, the forgetful functor \(\dvcat \to
\dcat\) has a left adjoint, \(\free{v} \colon \dcat \to \dvcat\),
which is called the \free{v}. \noproof
\end{corollary}

Dual to \valgobjs are \Cvalgobjs.

\begin{defn}
Let \(\vcat\) be a variety of algebras, \(\dcat\) a category with
finite \co{}products.  A \emph{\dCvobj{}} is a \Odvobj.
A \emph{morphism of \dCvobjs} is a morphism in \(\dcat\) which
intertwines the \Cvalg structures.  We denote \dCvcat by \(\dCvcat\).
\end{defn}

The morphisms are chosen such that there is an isomorphism of
categories
\[
  \dCvcat \cong \mOo{(\Odvcat)}
\]
and there is a natural forgetful functor \(\dCvcat \to \dcat\).
The analogue of proposition~\ref{prop:conlift} is the following.

\begin{proposition}
\label{prop:covlift}
To give an \dobj[\abs{\dCvobj}] the structure of a
\Cvalgobj is equivalent to giving a lift of the
covariant hom\hyp{}functor \(\Hom{\dcat}{\abs{\dCvobj}}{-} \colon
\dcat \to \scat\) to a functor \(\dcat \to \vcat\).  \noproof
\end{proposition}

We shall use similar notation: \(\Hom{\dcat}{\dCvobj}{\dobj}\) will
denote the \valg with underlying set
\(\Hom{\dcat}{\abs{\dCvobj}}{\dobj}\).

The main tool of our analysis is the link between representable
functors and functors with adjoints.  This result is a standard one
from general algebra, although one of its corollaries is perhaps the
best known result.

\begin{theorem}
\label{th:adj}
Let \(\dcat\) be a category that has finite products, is
\co{complete}, is \co{well\hyp{}powered}, is an (\(\m{E}\), \(\m{M}\))
category where \(\m{E}\) is closed under finite products, and is such that 
its finite products commute with filtered \co{limits}.  Let \(\vcat\) be a
variety of algebras.  Let \(\fcat\) be a category with
\co{equalisers}.  Let \(\func{G} \colon \fcat \to \dvcat\) be a
covariant functor.  Then the following statements are equivalent.
\begin{enumerate}
\item \(\func{G}\) has a left adjoint, \(\ladj{\func{G}}\).
\item The composition \(\abs{\func{G}} \colon \fcat \to \dcat\) of
\(\func{G}\) with the forgetful functor \(\dvcat \to \dcat\) has a
left adjoint, \(\ladj{\sabs{\func{G}}}\). \noproof
\end{enumerate} 
\end{theorem}

The relationship between the two adjoints is that there is a
\co{equaliser} sequence in \(\fcat\), natural in \(\dvobj\),
\[
  \xymatrix{
    \ladj{\sabs{\func{G}}}(\sabs{\free{v}(\sabs{\dvobj})})
    \ar@<2pt>[r]^(.6){r_{\dvobj}} \ar@<-2pt>[r]_(.6){s_{\dvobj}} &
    \ladj{\sabs{\func{G}}}(\sabs{\dvobj}) \ar[r]^{p_{\dvobj}} &
    \ladj{\func{G}}(\dvobj).
  }
\]

Working in the opposite category we obtain the corresponding result on
\co{algebra} objects. 

\begin{corollary}
\label{cor:adj}
Let \(\dcat\) be a category that has finite products, is
\co{complete}, is \co{well\hyp{}powered}, is an (\(\m{E}\), \(\m{M}\))
category where \(\m{E}\) is closed under finite products, and its
finite products commute with filtered \co{limits}.  Let \(\vcat\) be a
variety of algebras.  Let \(\fcat\) be a category with equalisers.
Let \(\func{G} \colon \fcat \to \dvcat\) be a contravariant functor.
Then the following statements are equivalent.
\begin{enumerate}
\item \(\func{G}\) is one of a mutually right adjoint pair.
\item The composition \(\abs{\func{G}} \colon \fcat \to \dcat\) of
\(\func{G}\) with the forgetful functor \(\dvcat \to \dcat\) is one of
a mutually right adjoint pair. \noproof
\end{enumerate}
\end{corollary}

Since a functor from a \co{complete} category into \(\scat\) is
representable if and only if it has a left adjoint, the following
standard result of general algebra readily follows from
proposition~\ref{prop:covlift}.

\begin{corollary}
\label{cor:adjrep}
\begin{enumerate}
\item Let \(\fcat\) be a \co{complete} category, \(\vcat\) a variety
of algebras.  For a covariant functor \(\func{G} \colon \fcat \to
\vcat\), the following statements are equivalent.
\begin{enumerate}
\item \(\func{G}\) has a left adjoint.
\item \(\func{G}\) is representable by a \fCvobj.
\item \(\abs{\func{G}}\) is representable by an \fobj.
\end{enumerate}

\item Let \(\fcat\) be a complete category, \(\vcat\) a variety
of algebras.  For a contravariant functor \(\func{G} \colon \fcat \to
\vcat\), the following statements are equivalent.
\begin{enumerate}
\item \(\func{G}\) is one of a mutually right adjoint pair.
\item \(\func{G}\) is representable by a \fvobj.
\item \(\abs{\func{G}}\) is representable by an \fobj. \noproof
\end{enumerate}
\end{enumerate}
\end{corollary}

We shall need to know various categorical properties of \(\vcat\).

\begin{theorem}
As a category, \(\vcat\) is complete, \co{complete},
well\hyp{}powered, extremally \co{well}\hyp{}powered, and is an
(extremal epi, mono) category.

A morphism is an extremal epimorphism if and only if the underlying
morphism of sets is an epimorphism.  Moreover, all extremal
epimorphisms are regular epimorphisms.  \noproof
\end{theorem}

\subsection{Graded Algebras}

We turn now to graded algebras.  A graded algebra has components
indexed by some (fixed) set and its operations go from components to
components rather than being globally defined.  The theory of graded
algebras is very similar to that of ungraded algebras.

We fix some (non\hyp{}empty) grading set \(\iset\).  We shall regard
this both as a set and as a (small) discrete category.  We write
\(\Zdcat\) for \Zdcat.  As \(\iset\) is a discrete category, there is
no distinction between covariant and contravariant functors from
\(\iset\).

An \Zdobj[\Zdobj] represents both a covariant and a contravariant
functor \(\dcat \to \Zscat\) via
\begin{align*}
  \cov{\Zdobj}(\dobj') &= \big(\isetelt \mapsto
  \Hom{\dcat}{\Zdobj(\isetelt)}{\dobj'}\big), \\
  \con{\Zdobj}(\dobj') &= \big(\isetelt \mapsto
  \Hom{\dcat}{\dobj'}{\Zdobj(\isetelt)}\big).
\end{align*}

To define a graded algebra, we first need to define the graded version
of a type.

\begin{defn}
A \emph{\(\iset\)\hyp{}graded type}, \(\Gotype\), is a triple
\((\abs{\Gotype}, i, o)\) where \(\abs{\Gotype}\) is a set, and
\[
i \colon \abs{\Gotype} \to \coprod_{m \in \N_0} \iset^m \qquad \text{
  and } \qquad
o \colon \abs{\Gotype} \to \iset
\]
are morphisms of sets.  For an operation \(\Goop\), we call
\(i(\Goop)\) the \emph{input} and \(o(\Goop)\) the \emph{output} of
\(\Goop\).  We define the \emph{arity} function, \(n \colon
\abs{\Gotype} \to \N_0\) by the composition
\[
  \abs{\Gotype} \xrightarrow{i} \coprod_{m \in \N_0} \iset^m \to
  \coprod_{m \in \N_0} \{*\} \cong \N_0.
\]
\end{defn}

We think of \(\coprod_{m \in \N_0} \iset^m\) as the set of finite
ordered sets of elements of \(\iset\) and so we interpret the element
in \(\iset^0\) as representing the empty subset of \(\iset\).  Under
the assumption that \(\dcat\) has finite products, for an
element \((\isetelt_1, \dotsc, \isetelt_m)\) and an \Zdobj[\Zdobj], we
write
\[
  \Zdobj(\isetelt_1, \dotsc, \isetelt_m) = \prod_{j=1}^m
  \Zdobj(\isetelt_j)
\]
with \(\Zdobj(\emptyset) = \term{d}\), the terminal object in
\(\dcat\).

\begin{defn}
Let \(\dcat\) be a category with finite products.

\begin{enumerate}
\item An \emph{\dGoobj[\dGoobj]{}} consists of an \Zdobj[\abs{\dGoobj}]
together with, for each \(\Goop \in \abs{\Gotype}\), a
\(\dcat\)\hyp{}morphism
\[
  \Goop_{\dGoobj} \colon \sabs{\dGoobj}({i(\Goop)}) \to
  \sabs{\dGoobj}(o(\Goop)).
\]
A morphism of \dGoobjs is a morphism of the underlying \Zdobjs which
intertwines the operations.

\item An \emph{\Goalg{}} is an \sGoobj.

\item We denote \dGocat by \(\dGocat\) and \sGocat by \(\Gocat\).  We
refer to the functor \(\dGocat \to \Zdcat\) which assigns to an \dGoobj
the underlying \Zdobj as the \emph{forgetful functor}.  We write the
underlying \Zdobj of an \dGoobj[\dGoobj] as \(\abs{\dGoobj}\).
\end{enumerate}
\end{defn}

All of the results for ungraded algebras have their counterparts in
graded algebras.

\begin{proposition}
\label{prop:conliftgrade}
To give an \Zdobj[\abs{\dGoobj}] the structure of an
\Goalgobj is equivalent to giving a lift of the
contravariant hom\hyp{}functor \(\Hom{\dcat}{-}{\abs{\dGoobj}} \colon
\dcat \to \Zscat\) to a functor \(\dcat \to \Gocat\).  \noproof
\end{proposition}

As before, if \(\dGoobj\) is an \dGoobj and \(\dobj\) is an \dobj we
shall write
\[
  \Hom{\dcat}{\dobj}{\dGoobj}
\]
for the corresponding \Goalg with underlying \Zsobj,
\[
  \isetelt \mapsto \Hom{\dcat}{\dobj}{\abs{\dGoobj}(\isetelt)}.
\]

Free \Goalgs exist under suitable circumstances.

\begin{proposition}
\label{prop:freegrade}
Let \(\dcat\) be a category with finite products which is closed under
filtered \co{}limits and such that finite products commute with
filtered \co{}limits.  Then the forgetful functor \(\dGocat \to
\Zdcat\) has a left adjoint, \(\free{Go} \colon \Zdcat \to \dGocat\),
which is called the \free{Go}. \noproof
\end{proposition}

Identities are defined by modifying the ungraded definition in the
obvious way.

\begin{defn}
Let \(\Gotype\) be a graded type.  An \emph{identity} for \Goalgs is a
triple \((\Zsobj,u,v)\) where \(\Zsobj\) is an \Zsobj such that
\(\coprod_{\isetelt \in \iset} \Zsobj(\isetelt)\) is a finite set and
\(u, v \in \abs{\free{Go}(\Zsobj)}(\isetelt)\) for some \(\isetelt \in
\iset\).
\end{defn}

Let \((\Zsobj, u, v)\) be an identity for \Goalgs with \(u, v \in
\abs{\free{Go}(\Zsobj)}(\isetelt_1)\).  Let \(\dGoobj\) be an \dGoobj.
Consider the \Zsobj
\[
  \isetelt_0 \mapsto \Hom{\dcat\big}{\prod_{\isetelt \in \iset}
    \abs{\dGoobj}(\isetelt)^{\Zsobj(\isetelt)}}{\abs{\dGoobj}(\isetelt_0)\big}.
\]
As \(\coprod_{\isetelt \in \iset} \Zsobj(\isetelt)\) is finite, the
product on the left is finite.  For each \(\isetelt_0 \in \iset\),
there is an obvious projection morphism in \(\dcat\)
\[
  \prod_{\isetelt \in \iset}
  \abs{\dGoobj}(\isetelt)^{\Zsobj(\isetelt)} \to
  \abs{\dGoobj}(\isetelt_0)^{\Zsobj(\isetelt_0)}
\]
and thus for \(\Zselt \in \Zsobj(\isetelt_0)\) we have a canonical
projection morphism in \(\dcat\)
\[
  \prod_{\isetelt \in \iset}
  \abs{\dGoobj}(\isetelt)^{\Zsobj(\isetelt)} \to
  \abs{\dGoobj}(\isetelt_0).
\]
This yields a \(\scat\) morphism
\[
  \Zsobj(\isetelt_0) \to \Hom{\dcat\big}{\prod_{\isetelt \in \iset}
    \abs{\dGoobj}(\isetelt)^{\Zsobj(\isetelt)}}{\abs{\dGoobj}(\isetelt_0)\big}
\]
and thus a natural transformation of functors, equivalently a
\(\Zscat\)\hyp{}morphism,
\[
  \Zsobj \to \Hom{\dcat\big}{\prod_{\isetelt \in \iset}
    \abs{\dGoobj}(\isetelt)^{\Zsobj(\isetelt)}}{\abs{\dGoobj}\big}.
\]
The same adjunction argument as in the ungraded case now produces
\(\dcat\)\hyp{}morphisms
\[
  u_{\dGoobj}, v_{\dGoobj} \colon \prod_{\isetelt \in \iset}
  \abs{\dGoobj}(\isetelt)^{\Zsobj(\isetelt)} \to
  \abs{\dGoobj}(\isetelt_1).
\]

\begin{defn}
An \dGoobj[\dGoobj] is said to \emph{satisfy} the identity
\((\Zsobj,u,v)\) if the two induced morphisms \(u_{\dGoobj},
v_{\dGoobj}\) are the same.
\end{defn}

\begin{defn}
A \emph{variety of graded algebras}, \(\Gvcat\), is specified by a
graded type, \(\Gotype\), and a set of identities for \Goalgs, \(J\).
It is the full subcategory of \(\Gocat\) consisting of all \Goalgs
which satisfy all of the identities in \(J\).

The pair \((\Gotype, J)\) is a \emph{presentation} of \(\Gvcat\).

Let \(\dcat\) be a category with finite products, \(\Gvcat\) a variety
of graded algebras with presentation \((\Gotype, J)\).  \dGvcatu,
\(\dGvcat\), is the full subcategory of \(\dGocat\) consisting of all
\dGoobjs which satisfy all of the identities in \(J\).
\end{defn}

Proposition~\ref{prop:conliftgrade} holds in the presence of identities.

\begin{proposition}
\label{prop:conliftidgrade}
To give an \Zdobj[\abs{\dGoobj}] the structure of a
\Gvalgobj is equivalent to giving a lift of the
contravariant hom\hyp{}functor \(\Hom{\dcat}{-}{\abs{\dGoobj}} \colon
\dcat \to \Zscat\) to a functor \(\dcat \to \Gvcat\).  \noproof
\end{proposition}

The same conditions as in the ungraded case allow us to impose identities
and so get free \Gvalgobjs.

\begin{theorem}
\label{th:impidgrade}
Let \(\dcat\) be a complete category with finite products.  Suppose
that \(\dcat\) is an (\(\m{E}\), \(\m{M}\)) category for some classes
\(\m{E}\) of epimorphisms and \(\m{M}\) of monomorphisms, that
\(\m{E}\) is closed under taking finite products, and that \(\dcat\)
is \(\m{E}\)\hyp{}\co{well}\hyp{}powered.  Then the inclusion functor
\(\dGvcat \to \dGocat\) has a left adjoint, \(\dGocat \to \dGvcat\),
called \emph{imposition of identities}. \noproof
\end{theorem}

\begin{corollary}
Under the conditions of theorem~\ref{th:impidgrade} and
proposition~\ref{prop:freegrade}, the forgetful functor \(\dGvcat \to
\Zdcat\) has a left adjoint, \(\free{Gv} \colon \Zdcat \to \dGvcat\),
which is called the \free{Gv}. \noproof
\end{corollary}

Dual to \Gvalgobjs are \CGvalgobjs.

\begin{defn}
Let \(\Gvcat\) be a variety of graded algebras, \(\dcat\) a category
with finite \co{}products.  A \emph{\dCGvobj{}} is a \OdGvobj.  A
\emph{morphism of \dCGvobjs{}} is a morphism in \(\Zdcat\) which
intertwines the \dCGvobj structures.  We denote \dCGvcat by
\(\dCGvcat\).
\end{defn}

The analogue of proposition~\ref{prop:conliftgrade} is the following.

\begin{proposition}
\label{prop:covliftgrade}
To give an \Zdobj[\abs{\dCGvobj}] the structure of a
\CGvalgobj is equivalent to giving a lift of the
covariant hom\hyp{}functor \(\Hom{\dcat}{\abs{\dCGvobj}}{-} \colon
\dcat \to \Zscat\) to a functor \(\dcat \to \Gvcat\).  \noproof
\end{proposition}

We shall use similar notation: \(\Hom{\dcat}{\dCGvobj}{\dobj}\) will
denote the \Gvalg with underlying object
\(\Hom{\dcat}{\abs{\dCGvobj}}{\dobj}\) in \(\Zscat\).

Theorem~\ref{th:adj} easily generalises to the graded case.

\begin{theorem}
\label{th:adjgrade}
Let \(\dcat\) be a category that has finite products, is
\co{complete}, is \co{well\hyp{}powered}, is an (\(\m{E}\), \(\m{M}\))
category where \(\m{E}\) is closed under finite products, and its
finite products commute with filtered \co{limits}.  Let \(\Gvcat\) be a
variety of graded algebras.  Let \(\fcat\) be a category with
\co{equalisers}.  Let \(\func{G} \colon \fcat \to \dGvcat\) be a
covariant functor.  Then the following statements are equivalent.
\begin{enumerate}
\item \(\func{G}\) has a left adjoint, \(\ladj{\func{G}}\).
\item The composition \(\abs{\func{G}} \colon \fcat \to \Zdcat\) of
\(\func{G}\) with the forgetful functor \(\dGvcat \to \Zdcat\) has a
left adjoint, \(\ladj{\sabs{\func{G}}}\). \noproof
\end{enumerate}
\end{theorem}

We have the same relationship between the two adjoints as in the
ungraded case: there is a \co{equaliser} sequence in \(\fcat\),
natural in \(\dGvobj\),
\[
  \xymatrix{
    \ladj{\sabs{\func{G}}}(\sabs{\free{Gv}(\sabs{\dGvobj})})
    \ar@<2pt>[r]^(.6){r_{\dGvobj}} \ar@<-2pt>[r]_(.6){s_{\dGvobj}} &
    \ladj{\sabs{\func{G}}}(\sabs{\dGvobj}) \ar[r]^{p_{\dGvobj}} &
    \ladj{\func{G}}(\dGvobj).
  }
\]

The graded version of corollary~\ref{cor:adj} follows immediately.  To
get the graded version of corollary~\ref{cor:adjrep} we need to
understand the relationship between adjunctions and representability
in the graded case.

\begin{lemma}
\label{lem:gradjrep}
Let \(\dcat\) be a \co{complete} category.  A covariant functor
\(\func{G} \colon \dcat \to \Zscat\) has a left adjoint if and only if
it is representable by an \Zdobj.
\end{lemma}

\begin{proof}
Suppose that \(\func{G}\) has a left adjoint, say \(\func{H}
\colon \Zscat \to \dcat\).  We extend this to a functor
\(\func[mZo]{H} \colon \ZZscat \to \Zdcat\) in the obvious way.  Let
\(I\) be the \ZZsobj defined by
\[
  \isetelt \mapsto \Bigg( \isetelt' \mapsto \begin{cases} \{*\} &
  \text{if } \isetelt = \isetelt' \\
  \emptyset &\text{otherwise}
  \end{cases} \Bigg).
\]
Then for an \Zsobj[\Zsobj] we have isomorphisms in \(\Zscat\)
\begin{align*}
  \Hom{\Zscat}{I}{\Zsobj} &= \big( \isetelt \mapsto
  \Hom{\Zscat}{I(\isetelt)}{\Zsobj} \big) \\
&\cong \big( \isetelt \mapsto \prod_{\isetelt' \in \iset}
\Hom{\scat}{I(\isetelt)(\isetelt')}{\Zsobj(\isetelt')} \big) \\
&\cong \big(\isetelt \mapsto \Hom{\scat}{\{*\}}{\Zsobj(\isetelt)}
\times \prod_{\isetelt' \ne \isetelt}
\Hom{\scat}{\emptyset}{\Zsobj(\isetelt')} \big) \\
&\cong \big(\isetelt \mapsto \Zsobj(\isetelt) \times \prod_{\isetelt'
  \ne \isetelt} \{*\} \big) \\
&\cong \big(\isetelt \mapsto \Zsobj(\isetelt) \big) \\
&\cong \Zsobj,
\end{align*}
all natural in \(\Zsobj\).  Hence for \(\dobj\) an \dobj there are
natural isomorphisms
\[
  \func{G}(\dobj) \cong \Hom{\Zscat}{I}{\func{G}(\dobj)} \cong
  \Hom{\dcat}{\func[mZo]{H}(I)}{\dobj}
\]
and so \(\func{G}\) is represented by the \Zdobj[\mZo{\func{H}}(I)].

Conversely, suppose that \(\func{G}\) is represented by the
\Zdobj[G].  Let \(\Zsobj\) be an \Zsobj.  We have the following
natural isomorphisms of sets
\begin{align*}
\Hom{\Zscat}{\Zsobj}{\func{G}(\dobj)}
&\cong \Hom{\Zscat}{\Zsobj}{\Hom{\dcat}{G}{\dobj}} \\
&\cong \prod_{\isetelt \in \iset}
\Hom{\scat}{\Zsobj(\isetelt)}{\Hom{\dcat}{G(\isetelt)}{\dobj}} \\
&\cong \prod_{\isetelt \in \iset}
\Hom{\dcat}{G(\isetelt)}{\dobj}^{\Zsobj(\isetelt)} \\
&\cong \prod_{\isetelt \in \iset}
\Hom{\dcat\big}{\coprod_{\Zsobj(\isetelt)} G(\isetelt)}{\dobj\big} \\
&\cong \Hom{\dcat\big}{\coprod_{\isetelt \in
    \iset}\coprod_{\Zsobj(\isetelt)} G(\isetelt)}{\dobj\big}.
\end{align*}

Therefore we define the functor \(\func{H} \colon \Zscat \to
\dcat\) on objects by
\[
  \func{H}(\Zsobj) = \coprod_{\isetelt \in \iset}
  \coprod_{\ZZsobj(\isetelt)} G(\isetelt)
\]
and in the obvious way on morphisms.  This is the required left
adjoint.
\end{proof}

As a corollary we deduce the graded version of
corollary~\ref{cor:adjrep}.

\begin{corollary}
\label{cor:adjrepgrade}
\begin{enumerate}
\item Let \(\dcat\) be a \co{complete} category, \(\Gvcat\) a variety
of graded algebras.  For a covariant functor \(\func{G} \colon \dcat
\to \Gvcat\), the following statements are equivalent.
\begin{enumerate}
\item \(\func{G}\) has a left adjoint.
\item \(\func{G}\) is representable by a \dCGvobj.
\item \(\abs{\func{G}}\) is representable by an \Zdobj.
\end{enumerate}

\item Let \(\dcat\) be a complete category, \(\Gvcat\) a variety of
graded algebras.  For a contravariant functor \(\func{G} \colon \dcat
\to \Gvcat\), the following statements are equivalent.
\begin{enumerate}
\item \(\func{G}\) is  part of a mutually right
adjoint pair.
\item \(\func{G}\) is representable by an \dGvobj.
\item \(\abs{\func{G}}\) is representable by an \Zdobj. \noproof
\end{enumerate}
\end{enumerate}
\end{corollary}

The categorical properties of \(\Gvcat\) are the same as those of
\(\vcat\).

\begin{theorem}
\label{th:gvcatprop}
As a category, \(\Gvcat\) is complete, \co{complete},
well\hyp{}powered, extremally \co{well}\hyp{}powered, and is an
(extremal epi, mono) category.

A morphism is an extremal epimorphism if and only if the underlying
morphism of \Zsobjs is an epimorphism.  Moreover, all extremal
epimorphisms are regular epimorphisms.  \noproof
\end{theorem}

The work of the following sections can be viewed simply as
applications of corollaries~\ref{cor:adjrep}
and~\ref{cor:adjrepgrade}.

\subsection{The Tall\hyp{}Wraith Monoidal Structure}
\label{sec:twmon}

The categories \(\vCvcat\) and \(\GvCGvcat\) have a monoidal structure
corresponding to composition of (representable) functors.  The first
trace of this that we have discovered in the literature is \cite{pf2}
where it is referred to as the \emph{tensor product of algebras}.  The
first study of the corresponding monoids that we have found is
\cite{dtgw} which deals with the category of commutative, unital
rings.  As we are similarly interested in the monoids, we have elected
to call them \emph{Tall\hyp{}Wraith monoids}.  For consistency, and
because the terminology of tensor products is already somewhat
overused, we name the monoidal structure on \(\vCvcat\) and
\(\GvCGvcat\) the \emph{Tall\hyp{}Wraith monoidal structure}.

Although, as we have just said, the ideas in this section go back at
least to \cite{pf2}, we have not been able to find a reference which
covers all that we wish to do; in particular,
theorem~\ref{th:varswitch} and the extensions to graded algebras.  On
the other hand, these results are not central to this paper but rather
are a guide to what to expect in the filtered context and so we have
not included their proofs here.  The missing proofs can be found in
\cite{assw3}.

In addition to \cite{pf2} and \cite{dtgw} mentioned above, similar
ideas occur in \cite{jbbw} and \cite{gbah}.

\begin{theorem}
\label{th:twstr}
Let \(\vcat\) be a variety of algebras.  The category \(\vCvcat\) has
a monoidal structure which we write as \((\vCvcat, \twprod, I)\).  The
functor \(\vCvcat \to \covfun(\vcat,\vcat)\) given by sending a
\vCvobj to the covariant functor that it represents, is strong
monoidal.  The underlying \vobj of the unit, \(I\), is isomorphic to
\(\free{v}(\{*\})\). \noproof
\end{theorem}

We shall not give a full proof of this result here; the idea can be
found in \cite{pf2} and a full proof is in \cite{assw3}.  We shall
give a description of the product pairing as this will be important
later.

For a \vCvobj[\vCvobj] let us write \(\cov{\vCvobj} \colon \vcat \to
\vcat\) for the covariant functor that it represents.  By
corollary~\ref{cor:adjrep}, this functor has a left adjoint which we
denote by \(\ladj{\vCvobj}\).  Now \vCvobjs are \vobjs with extra
structure; this extra structure involves morphisms from the underlying
\vobj to iterated \co{products} of it.  As \(\ladj{\vCvobj}\) is a
left adjoint, it preserves \co{products} and thus lifts to a functor
\(\Cladj{\vCvobj} \colon \vCvcat \to \vCvcat\).  The assignment
\(\vCvobj \mapsto \Cladj{\vCvobj}\) is functorial in \(\vCvobj\).  The
pairing on \(\vCvcat\) is, up to natural isomorphism, given on objects
by
\[
  (\vCvobj_1, \vCvobj_2) \mapsto
  \Cladj{\vCvobj_2} (\vCvobj_1).
\]
It has the property that we have natural isomorphisms
\[
    \vcat(\vCvobj_1,\vcat(\vCvobj_2, \vCvobj_3))
    \cong\vcat(\vCvobj_1\twprod \vCvobj_2, \vCvobj_3).
\]

Theorem~\ref{th:twstr} readily adapts to the following situations.

\begin{proposition}
Let \(\dcat\) be a \co{complete} category, \(\vcat\) and \(\wcat\)
varieties of algebras.  There are products
\begin{align*}
\vCvcat \times \dCvcat &\to \dCvcat, \\
\mOo{(\vCvcat)} \times \dvcat &\to \dvcat, \\
\vCwcat \times \vCvcat &\to \vCwcat,
\end{align*}
all compatible with the monoidal structure of \(\vCvcat\) and with
composition of representable functors.  \noproof
\end{proposition}

We write all of the pairings using the notation \(-\twprod-\).

There are two things to note from this generalisation.  Firstly that
there are two pairings which involve \(\vCvcat\) and \(\vcat\).  The
first views \(\vcat\) as \(\svcat\) and so comes from the middle
pairing above with \(\dcat=\scat\); in terms of functors we have
\[
  \con{(\vCvobj \twprod \vobj)}(\sobj) = \Hom{\scat}{\sobj}{\vCvobj
    \twprod \vobj} \cong
  \Hom{\vcat}{\vCvobj}{\Hom{\scat}{\sobj}{\vobj}}.
\]
The second views \(\vcat\) as \(\vCscat\) and so comes from the third
pairing with \(\wcat=\scat\); in terms of functors we have
\[
  \cov{(\vobj \twprod \vCvobj)}(\vobj') = \Hom{\vcat}{\vobj \twprod
    \vCvobj}{\vobj'} \cong
  \Hom{\vcat}{\vobj}{\Hom{\vcat}{\vCvobj}{\vobj'}}.
\]
This latter pairing was the one considered in \cite{dtgw} with
\(\vcat\) the category of commutative, unital rings.

The second thing to note from this generalisation is the annoyance of
having a partially contravariant pairing.  Providing \(\dcat\) is
sufficiently structured this can be countered;  again, the
details can be found in \cite{assw3}.

\begin{theorem}
\label{th:varswitch}
Let \(\dcat\) be a category satisfying the conditions of
theorem~\ref{th:adj}.  Then there is a pairing
\[
  \vCvcat \times \dvcat \to \dvcat, \qquad (\vCvobj, \dvobj) \mapsto
  \vCvobj \circledast \dvobj,
\]
which is covariant in both arguments and satisfies
\[
  \Hom{\dvcat}{\vCvobj \circledast \dvobj}{\dvobj'} \cong
  \Hom{\dvcat}{\dvobj}{\vCvobj \twprod \dvobj'}
\] 
naturally in all arguments. \noproof
\end{theorem}

In a monoidal category it is natural to consider monoids.

\begin{defn}
Let \(\vcat\) be a variety of algebras.  A \emph{Tall\hyp{}Wraith
  \(\vcat\)\hyp{}monoid} is a monoid in \(\vCvcat\).  We write the
category of such monoids as \(\vCvTtcat\).
\end{defn}

These were discussed briefly in \cite[chs 63, 64]{gbah}, though
without explicit reference to the underlying monoidal structure on
\(\vCvcat\).

With a monoid one can consider modules for that monoid.  Since the
monoidal category \(\vCvcat\) acts on other categories we can consider
modules that are not \vCvobjs.  That is, if \(\vCvTtobj\) is a
Tall\hyp{}Wraith \(\vcat\)\hyp{}monoid and \(\dcat\) is a
\co{complete} category then we can consider \dCvobjs[\dCvobj] for
which there is a \(\dCvcat\)\hyp{}morphism
\[
  \vCvTtobj \twprod \dCvobj \to \dCvobj
\]
satisfying the obvious coherences.

In \cite{gbah} the authors show that the category of \vobjs with an
action of a Tall\hyp{}Wraith \(\vcat\)\hyp{}monoid is again a variety
of algebras.  Extending this, we easily see that a \dvobj or \dCvobj
is a module for a Tall\hyp{}Wraith \(\vcat\)\hyp{}monoid if and only
if the corresponding functor \(\dcat \to \vcat\) factors through the
category of \vobjs with an action of the Tall\hyp{}Wraith
\(\vcat\)\hyp{}monoid.

Two remarks are worth making at this juncture.  Firstly, if \(\wcat\)
is another variety of algebras then the structure of a
\(\vCvTtobj\)\hyp{}module on a \vCwobj does not have such an
interpretation since a \vCwobj represents a functor \emph{out} of
\vcat.  Secondly, due to the variance shift, a
\(\vCvTtobj\)\hyp{}module in \(\dvcat\) is better thought of as a
\(\vCvTtobj\)\hyp{}\co{module} as the required morphism is
\[
  \dvobj \to \vCvTtobj \twprod \dvobj.
\]
We can surmount this using the product \(\circledast\) since the
adjunction turns a \co{action} as above into a more
normal\hyp{}looking action.  That is to say, if a \dvobj[\dvobj] is a
\(\vCvTtobj\)\hyp{}\co{module} for \(\twprod\) with \co{action}
morphism
\[
  \dvobj \to \vCvTtobj \twprod \dvobj
\]
then it is a \(\vCvTtobj\)\hyp{}module for \(\circledast\) with action
morphism
\[
  \vCvTtobj \circledast \dvobj \to \dvobj.
\]

\medskip

The adaptation of all this to the graded situation is straightforward.

\begin{theorem}
Let \(\Gvcat\) be a variety of graded algebras.  The category
\(\GvCGvcat\) has a monoidal structure which we write as \((\GvCGvcat,
\twprod, I)\).  The functor \(\GvCGvcat \to \covfun(\Gvcat,\Gvcat)\),
given by sending a \GvCGvobj to the covariant functor that it
represents, is strong monoidal.  \noproof
\end{theorem}

As before we shall give a description of the product.  A
\GvCGvobj[\GvCGvobj] represents a functor \(\cov{\GvCGvobj} \colon
\Gvcat \to \Gvcat\).  By corollary~\ref{cor:adjrepgrade} this functor
has a left adjoint \(\ladj{\GvCGvobj} \colon \Gvcat \to \Gvcat\).  We
extend this to a functor \(\Zladj{\GvCGvobj} \colon \ZGvcat \to
\ZGvcat\) in the obvious way.  \Co{products} in a graded category are
formed by taking component\hyp{}by\hyp{}component \co{products},
whence \(\Zladj{\GvCGvobj}\) preserves \co{products} because
\(\ladj{\GvCGvobj}\) does.  Hence it lifts to a functor
\(\Cladj{\GvCGvobj} \colon \GvCGvcat \to \GvCGvcat\).  This has the
property that the adjunction isomorphism lifts to an isomorphism of
\Gvalgs
\[
  \Hom{\Gvcat}{\GvCGvobj_1}{\Hom{\Gvcat}{\GvCGvobj_2}{\Gvobj}} \cong
  \Hom{\Gvcat}{\Cladj{\GvCGvobj_2}(\GvCGvobj_1)}{\Gvobj}.
\]
Thus there is a natural isomorphism of \GvCGvobjs
\[
  \GvCGvobj_1 \twprod \GvCGvobj_2 \cong
  \Cladj{\GvCGvobj_2}(\GvCGvobj_1).
\]

The other part of the structure that needs a word of explanation is
the representing object for the unit of the monoidal structure.  We
saw in the proof of lemma~\ref{lem:gradjrep} that the identity functor
\(\Zscat \to \Zscat\) is representable by an \ZZsobj, labelled \(I\)
in that proof. The free \Gvalg on the components of this
\ZZsobj represents the identity on \(\Gvcat\).

\begin{proposition}
Let \(\dcat\) be a \co{complete} category, \(\Gvcat\) and \(\Gwcat\)
varieties of graded algebras.  There are products
\begin{align*}
\GvCGvcat \times \dCGvcat &\to \dCGvcat, \\
\mOo{(\GvCGvcat)} \times \dGvcat &\to \dGvcat, \\
\GvCGwcat \times \GvCGvcat &\to \GvCGwcat,
\end{align*}
all compatible with the monoidal structure of \(\GvCGvcat\) and with
composition of representable functors.  \noproof
\end{proposition}

We write all of the pairings using the notation \(-\twprod-\).

We remark that the varieties of graded algebras, \(\Gvcat\) and
\(\Gwcat\), could be graded by different indexing sets.  This allows
us to take, for example, \(\Gwcat = \scat\) and so get the obvious
pairing
\[
  \Gvcat \times \GvCGvcat \to \Gvcat.
\]

We can remove the variance switch in the middle product by means of
the graded analogue of theorem~\ref{th:varswitch}.

\begin{theorem}
Let \(\dcat\) be a category satisfying the conditions of
theorem~\ref{th:adj}.  Then there is a pairing
\[
  \GvCGvcat \times \dGvcat \to \dGvcat
\]
which is covariant in both arguments and satisfies
\[
  \Hom{\dGvcat}{\GvCGvobj \circledast \dGvobj}{\dGvobj'} \cong
  \Hom{\dGvcat}{\dGvobj}{\GvCGvobj \twprod \dGvobj'}
\] 
naturally in all arguments. \noproof
\end{theorem}

\begin{defn}
Let \(\Gvcat\) be a variety of graded algebras.  A \emph{Tall\hyp{}Wraith
  \(\Gvcat\)\hyp{}monoid} is a monoid in \(\GvCGvcat\).  We write the
category of such monoids as \(\GvCGvTtcat\).
\end{defn}

The remarks regarding modules (and \co{modules}) for a
Tall\hyp{}Wraith \(\vcat\)\hyp{}monoid carry over to the graded case.

%% file: hopf.filter.tex
\section{Filtered Categories}
\label{sec:filter}

The purpose of this section is to introduce a categorical version of
a very specific type of topology.  What we wish to generalise is the
following situation: we have a topological space whose topology is the
projective topology for a family of maps into discrete spaces.  This
particular case is easy to put into a general categorical situation
and we do not need any of the usual machinery used to meld topology
and category theory.

In the first part we introduce the basic idea: filtered objects.  To
give a set, \(\sobj\), a topology in this fashion it is sufficient to
give a family of maps with source \(\sobj\).  Putting this into a
categorical context leads to \emph{projectively filtered objects} in
an arbitrary category.  We shall also define \emph{inductively
  filtered objects} since we shall need to consider how
contravarient functors transform filtered objects in one category into
filtered objects in another category.

In the example of topological spaces, we only need to
consider surjective morphisms and we can reduce an arbitrary
filtration to one in which all the morphisms are surjective.  We
cannot mirror this reduction in all categories and, moreover,
the condition that a functor preserve epimorphisms is more restrictive
than we wish to impose.  However, in certain categories there is a
reduction functor and we shall examine the extra features of the
theory that this introduces.

When we can consider these reduced filtrations it makes sense to
consider variations on the themes of being complete and being
Hausdorff.  Completion is not a purely topological concept, rather it
is a notion from the theory of \emph{uniform spaces}.  The correct
generalisation of these two notions to reduced filtered objects
involves examining the projective limit of the filtration.  There is a
morphism from the underlying object to this limit and we can ask
whether this morphism is a monomorphism, epimorphism, or isomorphism.
Being a monomorphism corresponds to the topology being Hausdorff
whilst being an epimorphism generalises the notion of completeness.

We start by introducing the most general form of filtrations before
moving on to the reduced version.  Once we have that then we can
consider the projective limit.

\subsection{Projective Filtrations}
\label{sec:projind}

We start with the general case of a filtration on an object in a
category.

\begin{defn}
Let \(\dcat\) be a category, \(\dobj\) an \dobj.  We define
\(\mDo{\dobj}\) to be the quasi\hyp{}ordered class whose elements are
\dmors with source \(\dobj\) and whose ordering is given by \(\delt_1
\ge \delt_2\) if there is a \dmor \(h\) such that \(h \delt_1 =
\delt_2\).

A \emph{projective filtration}, \(\Pdobj\), on \(\dobj\) is a
non\hyp{}empty, saturated, directed subclass of \(\mDo{\dobj}\).

We say that \(\Pdobj_1\) is \emph{stronger} than \(\Pdobj_2\) if
\(\Pdobj_2 \subseteq \Pdobj_1\).

If we are given a projective filtration \(\Pdobj\) on an \dobj without
having specified the underlying \dobj we shall write it as
\(\sabs{\Pdobj}\).  An element of \(\Pdobj\) is a \dmor which we shall
write as \(\Pdelt \colon \sabs{\Pdobj} \to \Pdobj_{\Pdelt}\).
\end{defn}

In this context, saturated means that if \(\delt_1 \ge \delt_2\) and
\(\delt_1 \in \Pdobj\) then \(\delt_2 \in \Pdobj\).

Let \(\Pdobj\) be a projective filtration on \(\dobj\).  Let \(f
\colon \dobj' \to \dobj\) be a \dmor.  The family of all elements of
\(\mDo{\dobj'}\) of the form \(\delt f\) for \(\delt\) an element of
\(\Pdobj\) is a non\hyp{}empty, saturated, directed subclass of
\(\mDo{\dobj'}\) and hence a projective filtration on \(\dobj'\).

\begin{defn}
We refer to this filtration as the \emph{pull back of \(\Pdobj\) by
  \(f\)} and write it as \(\con{f}(\Pdobj)\).
\end{defn}

This construction is strictly associative.

\begin{lemma}
Given \dmors \(\dobj_1 \xrightarrow{f} \dobj_2 \xrightarrow{g}
\dobj_3\) and a projective filtration \(\Pdobj\) on \(\dobj_3\), the
projective filtrations \(\con{f}(\con{g}(\Pdobj))\) and \(\con{(g
  f)}(\Pdobj)\) on \(\dobj_1\) are the same. \noproof
\end{lemma}

With these notions we can define a category of projective filtrations
on \dobjs.

\begin{defn}
We define \emph{\Pdcat, \(\Pdcat\)}.  Its objects are projective
filtrations on \dobjs.  A \Pdmor \(f \colon \Pdobj_1 \to \Pdobj_2\) is
a \dmor \(\sabs{f} \colon \sabs{\Pdobj_1} \to \sabs{\Pdobj_2}\) with
the property that \(\Pdobj_1\) is stronger than \(\con{f}(\Pdobj_2)\).
\end{defn}

By construction, the obvious functor \(\Pdcat \to \dcat\) is faithful.

Any projective filtration is completely determined by an initial
subclass, which per force is directed, and any non\hyp{}empty directed
subclass of \(\mDo{\dobj}\) determines a projective filtration by
\emph{saturation}; that is, if \(\delt_1, \delt_2\) are in
\(\mDo{\dobj}\) with \(\delt_1 \ge \delt_2\) and \(\delt_1\) is in the
specified class then we include \(\delt_2\).  It is clear that the
original directed class is initial for the resulting projective
filtration.  If \(\dcat\) has finite products then any subclass of
\(\mDo{\dobj}\), directed or not, determines a projective filtration:
first we include all finite products and then we saturate it.
Therefore we could choose to work with directed subclasses of
\(\mDo{\dobj}\), or even arbitrary subclasses, but the above
formulation of saturated subclasses is most directly analogous to a
topology on a set.  The correspondences are: projective filtration and
topology, directed subclass and a basis of the topology, subclass and
a subbasis of the topology.

We shall find it useful to have a characterisation of when a \dmor
lifts to a \Pdmor in terms of choices of initial subclasses of the
projective filtrations involved.  Let \(\Pdobj_1\) and \(\Pdobj_2\) be
projective filtrations in \(\dcat\) and let \(f \colon \sabs{\Pdobj_1}
\to \sabs{\Pdobj_2}\) be a \dmor on the underlying \dobjs.  Suppose
that we have initial subclasses of \(\Pdobj_1\) and \(\Pdobj_2\)
indexed by \(\Lambda_1\) and \(\Lambda_2\) respectively.  Then \(f\)
lifts to a \Pdmor if and only if for each \(\lambda_2\) in
\(\Lambda_2\) there is a \(\lambda_1\) in \(\Lambda_1\) and a \dmor
\(f_{\lambda_1, \lambda_2} \colon \Pdobj_{1, \lambda_1} \to \Pdobj_{2,
  \lambda_2}\) such that the following diagram of \dmors commutes
\[
  \xymatrix{
    \sabs{\Pdobj_1} \ar[r]^f \ar[d]_{\Pdelt_{1,\lambda_1}} &
    \sabs{\Pdobj_2} \ar[d]^{\Pdelt_{2,\lambda_2}} \\
    \Pdobj_{1,\lambda_1} \ar[r]^{f_{\lambda_1, \lambda_2}} &
    \Pdobj_{2,\lambda_2}.
  }
\]

\medskip

\begin{lemma}
The assignment \(\dcat \to \Pdcat\) underlies a \(2\)\hyp{}functor
of 2\hyp{}categories \(\cat \to \cat\).
\end{lemma}

\begin{proof}
For a covariant functor \(\func{G} \colon \dcat \to \ecat\) we define
a covariant functor \(\pfunc{G} \colon \Pdcat \to \Pecat\) in the
obvious way: for \(\Pdobj\) an \Pdobj, \(\pfunc{G}(\Pdobj)\) is the
saturation of the non\hyp{}empty, directed subclass of
\(\mDo{\func{G}(\sabs{\Pdobj})}\) consisting of \(\func{G}(\Pdelt)\)
for \(\Pdelt \in \Pdobj\).  For a \Pdmor \(f \colon \Pdobj_1 \to
\Pdobj_2\), \(\pfunc{G}(f)\) has underlying \emor
\(\func{G}(\sabs{f})\).  That this is an \Pemor is obvious.  This
construction is compatible with composition in that \(\pfunc{G H} =
\pfunc{G} \pfunc{H}\).  It is obvious that identity functors map to
identity functors.

Similarly, if \(\nu \colon \func{G} \to \func{H}\) is a natural
transformation of covariant functors \(\dcat \to \ecat\), we define a
natural transformation \(\mPo{\,\nu\,} \colon \pfunc{G} \to
\pfunc{H}\).  For an \Pdobj[\Pdobj], the \Pemor
\(\mPo{\,\nu\,}_{\Pdobj}\) 
 has underlying \emor
\(\nu_{\sabs{\Pdobj}}\).  Again, this construction is obviously
compatible with composition and identity natural transformations.
\end{proof}

An important consequence of this is the following result.

\begin{proposition}
\label{prop:padjoint}
Let \(\func{G} \colon \dcat \to \ecat\) be a covariant functor and
suppose that it has a left adjoint, say \(\func{H} \colon \ecat \to
\dcat\), then \(\pfunc{H}\) is left adjoint to
\(\pfunc{G}\). \noproof
\end{proposition}

If \(\dcat\) has a terminal object then the class of all projective
filtrations of a fixed \dobj is a (large!) complete lattice.
Its top and bottom elements provide adjoints to the forgetful functor
\(\Pdcat \to \dcat\).

\begin{proposition}
The forgetful functor \(\Pdcat \to \dcat\) has a left adjoint \(\disfunc
\colon \dcat \to \Pdcat\).  If \(\dcat\) has a terminal object then
the forgetful functor \(\Pdcat \to \dcat\) has a right adjoint
\(\indfunc \colon \dcat \to \Pdcat\).

For an \dobj[\dobj], \(\disfunc(\dobj)\) is \(\mDo{\dobj}\) whilst
\(\indfunc(\dobj)\) is the subclass of \(\mDo{\dobj}\)
consisting of all \dmors from \(\dobj\) to terminal objects in
\(\dcat\).
\end{proposition}

\begin{proof}
Clearly the descriptions given of \(\disfunc(\dobj)\) and
\(\indfunc(\dobj)\) do produce projective filtrations on \(\dobj\) and
if \(\Pdobj\) is another projective filtration on \(\dobj\) then we
have \(\indfunc(\dobj) \subseteq \Pdobj \subseteq \disfunc(\dobj)\);
the second inclusion by definition and the first as \(\Pdobj\) is
non\hyp{}empty.

From this, it is clear that if \(f \colon \dobj_1 \to \dobj_2\) is a
\dmor then it underlies \Pdmors \(\disfunc(\dobj_1) \to
\disfunc(\dobj_2)\) and \(\indfunc(\dobj_1) \to \indfunc(\dobj_2)\).
As the forgetful functor \(\Pdcat \to \dcat\) is faithful, this is
sufficient to define \(\disfunc\) and \(\indfunc\) on morphisms.

Clearly, if we apply either \(\disfunc \colon \dcat \to \Pdcat\) or
\(\indfunc \colon \dcat \to \Pdcat\) and the forgetful functor
\(\Pdcat \to \dcat\) then the resulting composition is the identity
functor on \(\dcat\).

Finally, from above  we see that the identity on \(\dobj\) lifts to
morphisms \(\disfunc(\dobj) \to \Pdobj \to \indfunc(\dobj)\).  These
provide the required natural transformations for the adjunctions.
\end{proof}

\begin{defn}
For an \dobj[\dobj] we refer to \(\disfunc(\dobj)\) as the
\emph{discrete (projective) filtration} on \(\dobj\) and
\(\indfunc(\dobj)\) as the \emph{indiscrete (projective) filtration} on
\(\dobj\) (assuming that \(\dcat\) has a terminal object).
\end{defn}

These two functors are very simple examples of a more general type of
functor.

\begin{defn}
To \emph{filter} a category is to give a functor \(\dcat \to \Pdcat\)
which is right inverse to the forgetful functor.  We call such a functor
a \emph{projective filtration functor}.
\end{defn}

\begin{examples}
\label{ex:filtrations}
\begin{enumerate}
\item
\label{ex:profin}
\let\tempobj\textfobj \renewcommand{\textfobj}{\hopfLower{F}inite
  \hopfLower{O}bject\hopfSingle{}{s}}
The first example is of a \emph{pro\hyp{}finite} filtration.  Let
\(\fcat\) be a non\hyp{}empty full subcategory of \(\dcat\) which is
closed under finite products.  We refer to \fobjalts as
\emph{\fobjs{}}.

Let \(\dobj\) be an \dobj.  We define a projective filtration on
\(\dobj\) as follows.  We start with the subclass of
\(\mDo{\dobj}\) consisting of all \dmors with target a \fobj.  Our
assumption on \(\fcat\) ensures that this is directed.  We saturate it
to produce a projective filtration.

It is straightforward to show that the assignment to an \dobj of its
pro\hyp{}finite filtration is functorial.  We therefore have the
pro\hyp{}finite filtration functor on \(\dcat\).
\let\textfobj\tempobj

\item
\label{ex:filfunc}
The second example of a category that can be filtered is that of
a filtered category.  We shall define a filtration functor
\[
  \Pdcat \to \PPdcat.
\]

Let \(\Pdobj\) be an \Pdobj.  We start by observing that for \(\Pdelt
\in \Pdobj\), the \dmor \(\Pdelt \colon \sabs{\Pdobj} \to
\Pdobj_{\Pdelt}\) is the underlying \dmor of a \Pdmor \(\Pdobj \to
\disfunc(\Pdobj_{\Pdelt})\).  Let us write \(\mPo{\,\Pdelt\,} \colon
\Pdobj \to \disfunc(\Pdobj_{\Pdelt})\) for the resulting \Pdmor.

The subclass of \(\mDo{\Pdobj}\) consisting of the elements
\(\mPo{\,\Pdelt\,}\) is directed, as the original projective filtration
was directed, and hence saturates to a projective filtration.

It is clear from its construction that this is functorial in
\(\Pdobj\).
\end{enumerate}
\end{examples}

\subsection{Reduced Filtrations}

Having looked at general filtrations, we now turn to a particular type
of filtration.  Let us consider the example of a topology on a set
defined by a projective filtration.  The structure of the category of
sets means that we can ensure that each of the maps in the filtration
is a surjection.  This has certain advantages which we wish to mirror
in our more general filtered categories.  Although the definition
below does not depend on any additional properties of the underlying
category, in order to do anything useful with it we need to assume
that this category is an extremally \co{well}\hyp{}powered (extremal
epi, mono) category.  We also want to know that the forgetful functor
\(\Pdcat \to \dcat\) has a right adjoint; for this we need to know
that \(\dcat\) has a terminal object.

\begin{defn}
Let \(\dcat\) be a category.  A projective filtration \(\Pdobj\) on an
\dobj is said to be \emph{reduced} if \(\Pdobj\) has an initial
subclass consisting of extremal epimorphisms. 

We write \(\Qdcat\) for the full subcategory of \(\Pdcat\) consisting
of reduced projective filtrations.  Let \(\qipfunc \colon
\Qdcat \to \Pdcat\) be the inclusion functor.
\end{defn}

We could broaden our definition by first choosing a reasonable class
of epimorphisms and considering those filtrations with morphisms in
that class, but we shall only be interested in extremal epimorphisms
and so we confine our attention to those.

Under the right conditions, the inclusion functor \(\qipfunc
\colon\Qdcat \to \Pdcat\) has a right adjoint.

\begin{proposition}
Let \(\dcat\) be an (extremal epi, mono) category.  Then there is a
\emph{reduction functor} \(\redfunc \colon \Pdcat \to \Qdcat\) which
is faithful.  The composition \(\redfunc\qipfunc \colon \Qdcat \to
\Qdcat\) is the identity functor.  The functor \(\redfunc\) is right
adjoint to \(\qipfunc\).
\end{proposition}

\begin{proof}
Let \(\Pdobj\) be a projective filtration in \(\dcat\).  We define
another projective filtration with the same underlying \dobj as
follows.  Each element \(\Pdelt \in \Pdobj\) is a \dmor \(\Pdelt
\colon \sabs{\Pdobj} \to \Pdobj_{\Pdelt}\).  By assumption on
\(\dcat\) this has a factorisation as \(m_{\Pdelt} \mQo{\,\Pdelt\,}\)
where \(\mQo{\,\Pdelt\,}\) is an extremal epimorphism and
\(m_{\Pdelt}\) a monomorphism.  Let us write \(\mQo{\Pdobj}_{\Pdelt}\)
for the target of \(\mQo{\,\Pdelt\,}\).  We claim that the class of
morphisms consisting of these \(\mQo{\,\Pdelt\,}\) is directed.  This
follows from the diagonal property of an (extremal epi, mono)
category.  Suppose that \(\Pdelt_1 \ge \Pdelt_2\) in \(\Pdobj\).  Then
there is some \dmor \(f \colon \Pdobj_{\Pdelt_1} \to
\Pdobj_{\Pdelt_2}\) such that \(\Pdelt_2 = f \Pdelt_1\).  Thus the
following is a commutative diagram in \(\dcat\).
\[
  \xymatrix{
    \sabs{\Pdobj} \ar[r]^{\mQo{\Pdelt_1}} \ar[d]_{\mQo{\Pdelt_2}} &
    \mQo{\Pdobj}_{\Pdelt_1} \ar[r]^{m_{\Pdelt_1}} &
    \Pdobj_{\Pdelt_1} \ar[ddll]^{f} \\
    \mQo{\Pdobj}_{\Pdelt_2} \ar[d]_{m_{\Pdelt_2}} \\
    \Pdobj_{\Pdelt_2} 
  }
\]
The diagonal property of the (extremal epi, mono)\hyp{}factorisations
now implies the existence of a \dmor \(\mQo{\Pdobj}_{\Pdelt_1} \to
\mQo{\Pdobj}_{\Pdelt_2}\) which fits into the above diagram.  Hence in
\(\mDo{\sabs{\Pdobj}}\) we have \(\mQo{\Pdelt_1} \ge
\mQo{\Pdelt_2}\).  Thus as \(\Pdobj\) is directed, the class
\(\{\mQo{\,\Pdelt\,}\}\) is also directed.  Its saturation is thus a
projective filtration which, by construction, is reduced.  Let us
write the result as \(\mQo{\Pdobj}\).  We define \(\redfunc\) on
objects by \(\redfunc(\Pdobj) = \mQo{\Pdobj}\).

To define \(\redfunc\) on morphisms we need to examine its interaction
with pull backs.  Let \(\Pdobj\) be a projective filtration on an
\dobj, let \(\dobj\) be an \dobj, and let \(f \colon \dobj \to
\sabs{\Pdobj}\) be a \dmor.  We wish to compare
\(\con{f}(\mQo{\Pdobj})\) with \(\mQo{\con{f}(\Pdobj)}\).  We obtain
an initial subclass for \(\con{f}(\mQo{\Pdobj})\) by taking the
extremal epimorphisms appearing in the
(extremal epi, mono)\hyp{}factorisations of elements of \(\Pdobj\) and
composing with \(f\).  That is, it consists of the \dmors
\(\mQo{\,\Pdelt\,} f \colon \dobj \to \mQo{\Pdobj_{\Pdelt}}\) where
\(\Pdelt \colon \sabs{\Pdobj} \to \Pdobj_{\Pdelt}\) is an element of
\(\Pdobj\) with (extremal epi, mono)\hyp{}factorisation \(m_{\Pdelt}
\mQo{\,\Pdelt\,}\) and intervening \dobj[\mQo{\Pdobj_{\Pdelt}}].
On the other hand, an initial subclass of \(\mQo{\con{f}(\Pdobj)}\)
consists of the extremal epimorphisms appearing in the (extremal epi,
mono)\hyp{}factorisations of the \dmors \(\Pdelt f\) for \(\Pdelt\) in
\(\Pdobj\).  For \(\Pdelt\) in \(\Pdobj\) we therefore have the
diagram
\[
  \xymatrix{
    \dobj \ar[r]^f \ar[d]^e &
    \sabs{\Pdobj} \ar[r]^{\mQo{\,\Pdelt\,}} 
    & \mQo{\Pdobj_{\Pdelt}} \ar[d]^{m_{\Pdelt}} \\
    \dobj_{\Pdelt} \ar[rr]^m &&
    \Pdobj_{\Pdelt}
  }
\]
where \(m e\) is the (extremal epi, mono)\hyp{}factorisation of
\(\Pdelt f\).  As \(\dcat\) is an (extremal epi, mono) category there
is a \dmor \(\dobj_{\Pdelt} \to \mQo{\Pdobj_{\Pdelt}}\) which fits
into this diagram.  Hence \(\con{f}(\mQo{\Pdobj}) \subseteq
\mQo{\con{f}(\Pdobj)}\).

The proof that the defining subclass of \(\mQo{\Pdobj}\) is directed
easily extends to show that if \(\Pdobj_1\) and \(\Pdobj_2\) are
projective filtrations on the same underlying \dobj with \(\Pdobj_1
\subseteq \Pdobj_2\) then \(\mQo{\Pdobj_1} \subseteq \mQo{\Pdobj_2}\).
Putting this together with the above, we see that if \(f \colon
\Pdobj_1 \to \Pdobj_2\) is a \Pdmor then
\(\con{\sabs{f}}(\mQo{\Pdobj_2}) \subseteq
\mQo{\con{\sabs{f}}(\Pdobj_2)} \subseteq \mQo{\Pdobj_1}\) and so
\(\sabs{f}\) also underlies a \Qdmor \(\mQo{\Pdobj_1} \to
\mQo{\Pdobj_2}\) which we denote by \(\redfunc(f)\).  As \(f\) and
\(\redfunc(f)\) have the same underlying \dmor, this assignment
clearly respects composition and identities whence we have a functor
\(\redfunc \colon \Pdcat \to \Qdcat\).

As (extremal epi, mono)\hyp{}factorisations in \(\dcat\) are unique up
to canonical isomorphism, the composition \(\Qdcat \to \Pdcat \to
\Qdcat\) is clearly the identity (saturation ensures that it is
actually the identity, rather than just isomorphic to the identity).
The obvious inclusion \(\Pdobj \subseteq \mQo{\Pdobj}\) provides the other
natural transformation in the adjunction.

Since both forgetful functors \(\Pdcat \to \dcat\) and \(\Qdcat \to
\dcat\) are faithful, and since the reduction functor
\(\Pdcat \to \Qdcat\) covers the identity on \(\dcat\), the reduction
functor must be faithful.
\end{proof}

The discrete filtration functor, \(\disfunc \colon \dcat \to \Pdcat\),
factors through \(\Qdcat\) with no modification since
\(\disfunc(\dobj)\) contains the initial subclass \(\{\dobj
  \xrightarrow{1} \dobj\}\) and every
isomorphism is an extremal epimorphism.  This provides us with a left
adjoint to the forgetful functor \(\Qdcat \to \dcat\).  The indiscrete
filtration functor, \(\indfunc \colon \dcat \to \Pdcat\), (assuming
that \(\dcat\) has a terminal object) is not so fortunate.  It is not
even true that the terminal morphism \(\termm{\dobj} \colon \dobj \to
\term{d}\) is always an epimorphism.  Thus we need to define
\(\mQo{\indfunc} \colon \dcat \to \Qdcat\) as the composition of
\(\indfunc\) with the reduction functor \(\Pdcat \to \Qdcat\).  This
is right adjoint to the forgetful functor \(\Qdcat \to \dcat\).

\subsection{Iso\hyp{}Filtrations}

In this section we consider those filtered objects that are analogous
to Hausdorf spaces and to complete uniform spaces.  In order to work
with these objects we need to make additional assumptions on our
underlying category, namely that it be complete and extremally
\co{well}\hyp{}powered.  As these filtered objects are a subclass of
the reduced filtered objects we also retain the assumptions of the previous
section.  Thus we assume that \(\dcat\) is a complete, extremally
\co{well}\hyp{}powered, (extremal epi, mono) category.  Note that
completeness implies the existence of a terminal object.

We start by defining another functor \(\Qdcat \to \dcat\) by taking
the limit of a reduced filtration.  We need to be in \(\Qdcat\) rather
than \(\Pdcat\) to ensure that the definition makes sense.

We start with the most general definition.  Let \(\dcat\) be a
category and \(\Qdobj\) an \Qdobj.  We define a category \((\Qdobj,
\dcat)\) by
\begin{description}
\item[Objects] elements of \(\Qdobj\),
\item[Morphisms] a morphism \(\Qdelt_1 \to \Qdelt_2\) is a \dmor[\dmor]
from the target of \(\Qdelt_1\) to the target of \(\Qdelt_2\) such
that \(\dmor \Qdelt_1 = \Qdelt_2\).
\end{description}
An alternative description of this, which explains the notation, is as
the full subcategory of the \emph{comma category} \((\sabs{\Qdobj},
\dcat)\) with class of objects the same as \(\Qdobj\).

Note that the quasi\hyp{}ordered class \(\Qdobj\) when viewed as a
category is a quotient of \((\Qdobj, \dcat)\) under the relation \(f
\sim g\) if \(f\) and \(g\) have the same source and target.

There is an obvious functor \((\Qdobj, \dcat) \to \dcat\) which sends
an object of \((\Qdobj, \dcat)\) to its target, and which regards a
morphism in \((\Qdobj, \dcat)\) as a \dmor.  Although the category
\((\Qdobj, \dcat)\) is usually large, the functor \((\Qdobj, \dcat)
\to \dcat\) still might have a limit.  By abuse of notation, we shall
refer to the limit of the functor \((\Qdobj, \dcat) \to \dcat\) as the
limit of \(\Qdobj\) and write it as \(\varprojlim \Qdobj\).

When this limit exists, it is obvious
that there is a natural \dmor
\[
  \sabs{\Qdobj} \to \varprojlim \Qdobj.
\]
The standard properties of limits show that if \(f \colon \Qdobj_1 \to
\Qdobj_2\) is a \Qdmor and \(\Qdobj_1\) and \(\Qdobj_2\) are such that
both of the appropriate limits exist, then there is a corresponding
\dmor \(f' \colon \varprojlim \Qdobj_1 \to \varprojlim \Qdobj_2\)
compatible with the above natural morphisms.

\begin{defn}
\label{def:isofilt}
If \(\Qdobj\) is an \Qdobj such that \(\varprojlim \Qdobj\) exists and
the natural morphism \(\sabs{\Qdobj} \to \varprojlim \Qdobj\) is an
isomorphism then we say that \(\Qdobj\) is an
\emph{iso\hyp{}filtration} on \(\sabs{\Qdobj}\).  We write \(\Kdcat\)
for the full subcategory of \(\Qdcat\) consisting of all such \Qdobjs.
We write \(\kiqfunc \colon \Kdcat \to \Qdcat\) for the inclusion
functor.
\end{defn}

Under reasonable conditions on \(\dcat\) we can show that, in fact,
every \Qdobj has a limit.

\begin{proposition}
If \(\dcat\) is an extremally \co{well}\hyp{}powered complete
category, then every \Qdobj has a limit.
\end{proposition}

\begin{proof}
To prove the required result, we observe that if \(\Qdobj\) is an
\Qdobj, it has an initial subclass consisting of extremal epimorphisms
with source \(\sabs{\Qdobj}\).  As \(\dcat\) is extremally
\co{well}\hyp{}powered, we can take this subclass to be small.  It is,
per force, directed.  We claim that this is an initial full
subcategory of \((\Qdobj, \dcat)\).  Let \(\Qdelt_1\) and \(\Qdelt_2\)
in \(\Qdobj\) be extremal epimorphisms.  Then as \(\Qdelt_1\) is an
epimorphism, there can be at most one \dmor with the property that
\(\dmor \Qdelt_1 = \Qdelt_2\).  Comparison of the ordering on elements
of \(\Qdobj\) with the definition of morphisms in \((\Qdobj, \dcat)\)
now shows that \(\Hom{(\Qdobj, \dcat)}{\Qdelt_1}{\Qdelt_2}\) has at
most one element and it has precisely one element if and only if
\(\Qdelt_1 \ge \Qdelt_2\).  Hence our initial subclass of \(\Qdobj\)
is a full subcategory of \((\Qdobj, \dcat)\).  It is clearly initial
in \((\Qdobj, \dcat)\).

As \(\dcat\) is complete we can therefore find a limit of the functor
\((\Qdobj, \dcat) \to \dcat\) by taking a limit of its restriction to
our small initial subclass of \(\Qdobj\).  This is unique up to
canonical isomorphism.  It also depends functorially on \(\Qdobj\)
since if we have a \Qdmor \(\Qdobj_1 \to \Qdobj_2\) then the pull back
of \(\Qdobj_2\) is contained in \(\Qdobj_1\), whence we get a
canonical \dmor on the limits, as they exist.
\end{proof}

By the remarks preceeding definition~\ref{def:isofilt}, these limits
fit together to form a functor \(\Qdcat \to \dcat\) and there is a
natural transformation of functors \(\Qdcat \to \dcat\) from the
forgetful functor to this limit functor.

\begin{defn}
We shall refer to the functor constructed above as the
\emph{projective limit functor} and write it as \(\limfunc \colon
\Qdcat \to \dcat\).
\end{defn}

Using the fact that a category of filtered objects is itself naturally
filtered, we obtain the following important construction.

\begin{proposition}
Let \(\dcat\) be a complete, extremally \co{well}\hyp{}powered
(extremal epi, mono) category.  Then there is a natural lift of the
projective limit functor \(\limfunc \colon \Qdcat \to \dcat\) to a
functor \(\comfunc \colon \Qdcat \to \Kdcat\) which is left adjoint to
the inclusion functor \(\kiqfunc \colon \Kdcat \to \Qdcat\).
Moreover, the composition \(\comfunc\kiqfunc \colon \Kdcat \to
\Kdcat\) is naturally isomorphic to the identity functor on
\(\Kdcat\).
\end{proposition}

Since the forgetful functor is naturally isomorphic to the projective
limit functor when restricted to \(\Kdcat\) we do not need
to specify which functor \(\Kdcat \to \dcat\) we are lifting along.

\begin{proof}
Let \(\Qdobj\) be an \Qdobj.  Let \(\iota \colon \sabs{\Qdobj} \to
\limfunc(\Qdobj)\) be the canonical \dmor.  Let \(\Qdelt\) be an
element of \(\Qdobj\).  By the definition of a limit, there is a
canonical \dmor
\[
  \mKo{\Qdelt} \colon \limfunc(\Qdobj) \to \Qdobj_{\Qdelt}
\]
such that \(\Qdelt = \mKo{\Qdelt}\iota\).

We claim that the subclass of \(\mDo{\limfunc(\Qdobj)}\) consisting of
the morphisms \(\mKo{\Qdelt}\) which arise in this fashion is a
projective filtration.  It is obviously non\hyp{}empty and directed.
To show that it is saturated, let \(\dmor \colon \limfunc(\Qdobj) \to
\dobj\) be a \dmor and suppose that it factors through
\(\mKo{\Qdelt}\) for some \(\Qdelt\) in \(\Qdobj\), say \(\dmor = f
\mKo{\Qdelt}\).  Then \(\dmor \iota = f \Qdelt\) whence \(\dmor
\iota\) is in \(\Qdobj\).  The \dmor \(f\) defines a morphism in
\((\Qdobj, \dcat)\) from \(\Qdelt\) to \(\dmor \iota\).  Hence, by the
definition of a limit, the canonical morphisms \(\mKo{(\dmor \iota)}
\colon \limfunc(\Qdobj) \to \dobj\) and \(\mKo{\Qdelt} \colon
\limfunc(\Qdobj) \to \Qdobj_{\Qdelt}\) satisfy \(\mKo{(\dmor \iota)} =
f \mKo{\Qdelt}\), whence \(\mKo{(\dmor \iota)} = \dmor\) and so
\(\dmor\) is in our chosen subclass.  Hence this subclass is
saturated and we have a projective filtration.

This construction clearly depends functorially on \(\Qdobj\).

Let us show that this is, in fact, a reduced projective filtration.
We need to show that it has an initial subclass of extremal
epimorphisms.  It certainly has an initial subclass consisting of
elements \(\mKo{\Qdelt}\) where \(\Qdelt \colon \sabs{\Qdobj} \to
\Qdobj_{\Qdelt}\) is an extremal epimorphism.  Since \(\Qdelt\) is an
epimorphism, so is \(\mKo{\Qdelt}\).  Suppose that we have a
factorisation \(\mKo{\Qdelt} = m f\) where \(m\) is a monomorphism;
let \(\dobj\) be the intervening \dobj.  As \(\dcat\) is an (extremal
epi, mono) category, the morphism \(f \iota \colon \sabs{\Qdobj} \to
\dobj\) has an (extremal epi, mono)\hyp{}factorisation, say \(f \iota
= m'e'\).  Then \(m m' e' = m f \iota = \mKo{\Qdelt} \iota = \Qdelt\).
As both \(m\) and \(m'\) are monomorphisms, their composition is still
a monomorphism whence as \(\Qdelt\) is an extremal epimorphism, \(m
m'\) is an isomorphism.  The monomorphism \(m\) is therefore a
retraction, with right inverse \(m' (m m')^{-1}\), and so is an
isomorphism.  Hence \(\mKo{\Qdelt}\) is an extremal epimorphism and so
the projective filtration that we have defined on \(\limfunc(\Qdobj)\)
is reduced.

Let \(\Qdobj'\) denote this reduced projective filtration on
\(\limfunc(\Qdobj)\).  We now claim that \(\Qdobj'\) is an
iso\hyp{}filtration.  Consider the categories \((\Qdobj, \dcat)\) and
\((\Qdobj', \dcat)\).  As \(\sabs{\Qdobj'} = \limfunc(\Qdobj)\) is the
limit of \((\Qdobj, \dcat) \to \dcat\), the assignment \(\Qdelt
\mapsto \mKo{\Qdelt}\) satisfies \(\mKo{(f \Qdelt)} = f \mKo{\Qdelt}\)
and thus we have a covariant functor \((\Qdobj, \dcat) \to (\Qdobj',
\dcat)\) which on objects is \(\Qdelt \mapsto \mKo{\Qdelt}\) and which
leaves morphisms alone (when viewed as \dmors).  The \dmor \(\iota
\colon \sabs{\Qdobj} \to \limfunc(\Qdobj)\) defines a covariant
functor \((\Qdobj', \dcat) \to (\Qdobj, \dcat)\).  Since
\(\mKo{\Qdelt} \iota = \Qdelt\), the composition is the identity
functor on \((\Qdobj, \dcat)\).  The argument above which showed that
the family \(\{\mKo{\Qdelt}\}\) is saturated shows that the
composition in the other direction is the identity functor on
\((\Qdobj', \dcat)\).  Since these functors do not change the targets
of the objects when viewed as \dmors and do not change the morphism
sets it is clear that this isomorphism intertwines the two functors
\((\Qdobj, \dcat) \to \dcat\) and \((\Qdobj', \dcat) \to \dcat\).
Hence as \(\limfunc(\Qdobj)\) is the limit of \((\Qdobj, \dcat) \to
\dcat\) it is also the limit of \((\Qdobj', \dcat) \to \dcat\) and
thus \(\Qdobj'\) is an iso\hyp{}filtration.

We therefore have a functor \(\comfunc \colon \Qdcat \to \Kdcat\) as
required.

It is clear from this construction that if we start with an \Kdobj
then all we do is replace the underlying \dobj by an isomorphic one
(with a specified isomorphism), whence the composition \(\comfunc
\kiqfunc\) is naturally isomorphic to the identity functor.

It is also clear from the construction that the natural transformation
from the forgetful functor to \(\limfunc\) underlies a natural
transformation from the identity on \(\Qdcat\) to the composition
\(\kiqfunc \comfunc\).

It is obvious that these natural transformations produce the
adjunction as stated.
\end{proof}

\begin{defn}
We shall refer to the functor \(\comfunc \colon \Qdcat \to \Kdcat\) as
the \emph{filtered projective limit functor}.
\end{defn}

We have already noted that the discrete filtration functor, \(\disfunc
\colon \dcat \to \Pdcat\), factors through \(\Qdcat\); it is equally
easy to see that it factors through \(\Kdcat\).  For comparison,
\(\mQo{\indfunc} \colon \dcat \to \Qdcat\) does not factor through
\(\Kdcat\), and there is little point in considering the composition
\(\comfunc \mQo{\indfunc}\).

\subsection{Categorical Properties}
\label{sec:catprop}

We wish to determine the categorical properties of the various
categories of filtered objects.  Certain results on \(\Pdcat\),
\(\Qdcat\), and \(\Kdcat\) have depended on the categorical properties
of \(\dcat\).  We shall want to work with categories such as
\(\PQdcat\) and so we need to know whether the various properties of
\(\dcat\) lift to, say, \(\Qdcat\).  We shall not consider categories
such as \(\QPdcat\) where the second type of filtration is more
restrictive than the first.  We also want to be able to apply the
results of section~\ref{sec:genalg} and therefore we need to know
other categorical properties to ensure that these apply.

In summary, we want to know the following.
\begin{enumerate}
\item Conditions on \(\dcat\) to ensure that \(\Qdcat\) is an
extremally \co{well}\hyp{}powered (extremal epi, mono) category with a
terminal object.

\item Conditions on \(\dcat\) to ensure that \(\Kdcat\) is a complete
extremally \co{well}\hyp{}powered (extremal epi, mono) category.

\item Conditions on \(\dcat\) to ensure that \(\Kdcat\) is
\co{complete}.

\item How to form products and (finite) \co{products} in each of
\(\Qdcat\) and \(\Kdcat\).
\end{enumerate}

Let us state all the conditions on \(\dcat\) that we need so that they
are collected in one place.  We assume that \(\dcat\) is
\begin{enumerate}
\item complete,
\item \co{complete},
\item an (extremal epi, mono) category, and
\item extremally \co{well}\hyp{}powered.
\end{enumerate}

Let us illustrate the various functors that we have.  We denote the
forgetful functor \(\Pdcat \to \dcat\) by \(\forfunc\), though we
shall still use the notation \(\sabs{\Pdobj}\) for
\(\forfunc(\Pdobj)\).  We shall also find it useful to have a notation
for the forgetful functor \(\Kdcat \to \dcat\) so we denote this by
\(\kfofunc \colon \Kdcat \to \dcat\).  Let us write the discrete
filtration functor \(\disfunc\) as a functor into \(\Kdcat\) rather
than \(\Pdcat\).  We have the following (non-commuting!)  diagram.
\[
  \xymatrix@=40pt{
    \Kdcat \ar@<1pt>[r]^{\kiqfunc} \ar@<2pt>[rrd]|\hole^(.7){\kfofunc}
    &
    \Qdcat \ar@<3pt>[l]^(.3){\comfunc} \ar@<1pt>[r]^{\qipfunc}
    \ar[d]^(.3){\limfunc} &
    \Pdcat \ar@<3pt>[l]^{\redfunc} \ar@<2pt>[d]^{\forfunc} \\
    & \dcat &
    \dcat \ar@<2pt>[ull]|\hole^(.2){\disfunc} \ar@<2pt>[u]^{\indfunc}
  }
\]
We have the following identities and adjunctions
\begin{align*}
\redfunc\qipfunc &= 1, &
\comfunc\kiqfunc &\cong 1, &
\kfofunc\disfunc &= 1, &
\forfunc\indfunc &= 1, \\
\kfofunc &= \forfunc\qipfunc\kiqfunc &
\forfunc\qipfunc\redfunc &= \forfunc, &
\limfunc\kiqfunc\comfunc &\cong \limfunc, &
\disfunc &\adjoint \kfofunc, \\
\comfunc &\adjoint \kiqfunc, &
\qipfunc &\adjoint \redfunc, &
\qipfunc\kiqfunc\disfunc &\adjoint \forfunc, &
\forfunc &\adjoint \indfunc.
\end{align*}
The functors \(\kiqfunc\), \(\qipfunc\), \(\disfunc\), and
\(\indfunc\) are fully faithful; \(\forfunc\) is faithful, whence also
\(\redfunc\) and \(\kfofunc\) are faithful.  The category \(\Kdcat\)
is a reflective subcategory of \(\Qdcat\).

We shall use some results from \cite[\S X]{hhgs} to transfer results
from \(\Qdcat\) to \(\Kdcat\).  These results refer to reflective
subcategories that are also full subcategories and closed under
isomorphisms.  These conditions are satisfied by \(\Kdcat\) inside
\(\Qdcat\).

We note in passing that in what follows we are essentially proving
that \(\Qdcat\) is \emph{topological} over \(\dcat\).  This is mildly
reassuring since our intention was to model a specific type of
topological space.

Let us start by defining push forward filtrations.

\begin{lemma}
\label{lem:pushlift}
Let \(\{\Qdobj_I\}\) be a class of \Qdobjs.  Let \((\dobj, f_I)\) be a
sink for the underlying class of \dobjs.  Then there is a reduced
projective filtration, say \(\Qdobj'\), on \(\dobj\) such that
\begin{enumerate}
\item each \dmor \(f_i \colon \sabs{\Qdobj_i} \to \dobj\) lifts to a
\Qdmor \(\mQo{f_i} \colon \Qdobj_i \to \Qdobj'\), and

\item if \(\Qdobj''\) is an \Qdobj and \(h \colon \dobj \to
\sabs{\Qdobj''}\) is a \dmor then \(h\) lifts to a \Qdmor \(\mQo{\,h\,}
\colon \Qdobj' \to \Qdobj''\) if and only if each \(h f_i \colon
\sabs{\Qdobj_i} \to \sabs{\Qdobj''}\) lifts to a \Qdmor.
\end{enumerate}
\end{lemma}

In the case that the class of \Qdobjs has only one element, say
\(\Qdobj\), and so the sink is just \(f \colon \sabs{\Qdobj} \to
\dobj\) then we shall refer to the resulting \Qdobj as the \emph{push
  forward} of \(\Qdobj\) along \(f\) and write it as
\(\cov{f}(\Qdobj)\).

When the class of \Qdobjs is empty, clearly the discrete filtration
has the required properties.

\begin{proof}
We define \(\Qdobj'\) as follows: it consists of all \dmors[g] with
source \(\dobj\) such that for each \(i \in I\), \(g f_i\) is in
\(\Qdobj_i\).

Firstly, \(\Qdobj'\) is not empty as \(\dcat\) has a terminal object
and so the terminal morphism from \(\dobj\) to this is in \(\Qdobj'\).

Secondly, \(\Qdobj'\) is directed.  To see this, suppose that \(g_1
\colon \dobj \to \dobj_1\) and \(g_2 \colon \dobj \to \dobj_2\) are in
\(\Qdobj'\).  Consider the \dmor \((g_1 \times g_2) \Delta \colon
\dobj \to \dobj_1 \times \dobj_2\) (which exists as \(\dcat\) is
complete).  If we can show that this lies in \(\Qdobj'\) then we are
done as it preceeds both \(g_1\) and \(g_2\).  Let \(i \in I\).  As
\(\Qdobj_i\) is directed and saturated, \((g_1 f_i \times g_2 f_i)
\Delta\) is in \(\Qdobj_i\).  By the functorality of products, this is
\((g_1 \times g_2) \Delta f_i\).  As this holds for all \(i\), \((g_1
\times g_2)\Delta\) is in \(\Qdobj'\).

Thirdly, \(\Qdobj'\) is saturated.  To see this, suppose that \(g_1\)
is in \(\Qdobj'\) and \(g_2 = k g_1\).  Then for \(i \in I\), \(g_1
f_i\) is in \(\Qdobj_i\) and so as this is saturated, \(k g_1 f_i\) is
in \(\Qdobj_i\).  Hence \(g_2\) is in \(\Qdobj'\).

We therefore have a projective filtration on \(\dobj\).  We shall now
show that it is reduced.  Let \(g \colon \dobj \to \dobj_1\) be in
\(\Qdobj'\).  As \(\dcat\) is an (extremal epi, mono) category it has
an (extremal epi, mono)\hyp{}factorisation \(g = m e\) with
intervening \dobj[\dobj_2]\added[AS]{,} say.  Let \(i \in I\).  The
\dmor \(g f_i \colon \sabs{\Qdobj_i} \to \dobj_1\) is in \(\Qdobj_i\)
and so since \(\Qdobj_i\) is reduced it factors through an extremal
epimorphism in \(\Qdobj_i\), say \(g f_i = h_i e_i\) with intervening
\dobj[\dobj_3] and \(e_i\) in \(\Qdobj_i\).  We therefore have a
commutative diagram in \(\dcat\).
\[
  \xymatrix{
    \sabs{\Qdobj_i} \ar[r]^{f_i} \ar[d]^{e_i} &
    \dobj \ar[d]^{g} \ar[r]^{e} &
    \dobj_2 \ar[ld]^{m} \\
    \dobj_3 \ar[r]^{h_i} &
    \dobj_1
  }
\]
As \(\dcat\) is an (extremal epi, mono) category, we can find a \dmor
\(\dobj_3 \to \dobj_2\) which fits into this diagram.  Hence \(e f_i\)
is in \(\Qdobj_i\).  Thus \(e\) is in \(\Qdobj'\) which is now shown
to be reduced.

By construction the \dmors \(f_i \colon \sabs{\Qdobj_i} \to \dobj\)
lift to \Qdmors \(\Qdobj_i \to \Qdobj'\).

Let \(h \colon \dobj \to \sabs{\Qdobj''}\) be a \dmor.  Let
\(\Qdelt''\) be an element of \(\Qdobj''\).  Then \(\Qdelt'' h\) is in
\(\Qdobj'\) if and only if \(\Qdelt'' h f_i\) is in \(\Qdobj_i\) for
all \(i \in I\).  Hence \(\con{h}(\Qdobj'') \subseteq \Qdobj'\) if and
only if \(\con{(h f_i)}(\Qdobj'') \subseteq \Qdobj_i\) for all \(i \in
I\).
\end{proof}

We already have the notion of pull back filtrations in the
not\hyp{}necessarily\hyp{}reduced case and it is easy to see that this
generalises.

\begin{lemma}
\label{lem:pulllift}
Let \(\{\Qdobj_I\}\) be a class of \Qdobjs.  Let \((\dobj, f_I)\) be a
source for the underlying class of \dobjs.  Then there is a reduced
projective filtration, say \(\Qdobj'\), on \(\dobj\) such that
\begin{enumerate}
\item each \dmor \(f_i \colon \dobj \to \sabs{\Qdobj_i}\) lifts to a
\Qdmor \(\mQo{f_i} \colon \Qdobj' \to  \Qdobj_i\), and

\item if \(\Qdobj''\) is an \Qdobj and \(h \colon \sabs{\Qdobj''} \to
\dobj\) is a \dmor then \(h\) lifts to a \Qdmor \(\mQo{\,h\,} \colon
\Qdobj'' \to \Qdobj'\) if and only if each \(f_i h \colon
\sabs{\Qdobj''} \to \sabs{\Qdobj_i}\) lifts to a \Qdmor.
\end{enumerate}
\end{lemma}

In the case that the class of \Qdobjs has one element we obtain the
reduced pull back filtration, \(\redfunc(\con{f}(\Qdobj))\).  If the
class of \Qdobjs is empty, we obtain the reduced indiscrete filtration
on \(\dobj\).

\begin{proof}
Each \dmor \(f_i \colon \dobj \to \sabs{\Qdobj_i}\) defines a pull
back filtration \(\con{f_i}(\Qdobj_i)\) on \(\dobj\).  The union of
these is a subclass of \(\mDo{\dcat}\).  As \(\dcat\) is
complete, we can find a smallest projective filtration containing this
subclass: first we include all finite products to ensure that it
is directed and non\hyp{}empty (via the empty product) and then we
saturate it.  Considered as an \Pdobj, this clearly has the required
properties.

We then apply \(\redfunc\) to this \Pdobj.  The properties then follow
from the fact that \(\redfunc\) is right adjoint to the inclusion
\(\qipfunc \colon \Qdcat \to \Pdcat\) and that both of these functors
cover the identity on \(\dcat\).
\end{proof}

As the forgetful functor \(\Qdcat \to \dcat\) is faithful and has both
a left and a right adjoint, it reflects and preserves monomorphisms
and epimorphisms.  Extremal epimorphisms are easy to characterise.

\begin{corollary}
\label{cor:qexepilift}
A \Qdmor \(f \colon \Qdobj_1 \to \Qdobj_2\) is an extremal epimorphism
if and only if \(\sabs{f}\) is an extremal epimorphism in \(\dcat\)
and \(\Qdobj_2 = \cov{\sabs{f}}(\Qdobj_1)\).
\end{corollary}

\begin{proof}
Let us show the ``only if'' part first, so that we suppose that \(f
\colon \Qdobj_1 \to \Qdobj_2\) is an extremal epimorphism.  We need to
show two things: that \(\sabs{f}\) is an extremal epimorphism and that
\(\Qdobj_2 = \cov{\sabs{f}}(\Qdobj_1)\).

As the forgetful functor \(\Qdcat \to \dcat\) has a right adjoint,
namely the reduced indiscrete functor, \(\sabs{f} \colon
\sabs{\Qdobj_1} \to \sabs{\Qdobj_2}\) is an epimorphism.  Suppose that
\(\sabs{f} = m g\) with \(m\) a monomorphism.  Let \(\dobj\) be the
intervening \dobj.  We put the reduced pull back filtration on
\(\dobj\) via \(m\).  Then by lemma~\ref{lem:pulllift}, \(m\) lifts to
a \Qdmor \(\mQo{m} \colon \redfunc(\con{m}(\Qdobj_2)) \to \Qdobj_2\)
and \(g\) lifts to a \Qdmor \(\mQo{\,g\,} \colon \Qdobj_1 \to
\redfunc(\con{m}(\Qdobj_2))\) with \(\mQo{m} \mQo{\,g\,} = f\).  As
the forgetful functor is faithful, \(\mQo{m}\) is a monomorphism.
Since \(f\) is an extremal epimorphism, \(\mQo{m}\) is thus an
isomorphism.  Hence \(m = \sabs{\mQo{m}}\) is an isomorphism.  Thus
\(\sabs{f}\) is an extremal epimorphism.

From lemma~\ref{lem:pushlift}, the identity on \(\sabs{\Qdobj_2}\)
underlies a \Qdmor \(\cov{\sabs{f}}(\Qdobj_1) \to \Qdobj_2\) and so
\(f\) factorises as \(\Qdobj_1 \to \cov{\sabs{f}}(\Qdobj_1) \to
\Qdobj_2\).  As the forgetful functor is faithful, the \Qdmor
\(\cov{\sabs{f}}(\Qdobj_1) \to \Qdobj_2\) is a monomorphism.  Hence as
\(f\) is an extremal epimorphism, it is an isomorphism.  As it covers
the identity on \(\sabs{\Qdobj_2}\), \(\cov{\sabs{f}}(\Qdobj_1)\) and
\(\Qdobj_2\) must in fact be the same projective filtrations.

Now let us show the ``if'' part.  Let \(f \colon \Qdobj_1 \to
\Qdobj_2\) be such that \(\sabs{f}\) is an extremal epimorphism and
\(\Qdobj_2 = \cov{f}(\Qdobj_1)\).  As the forgetful functor is
faithful, \(f\) is per force an epimorphism.  Suppose that we have a
factorisation of \(f\) as \(m g\) with \(m\) a monomorphism and
intervening \Qdobj[\Qdobj'].  We therefore have a factorisation of
\(\sabs{f}\) as \(\sabs{m} \sabs{g}\).  As the forgetful functor has a
left adjoint, \(\sabs{m}\) is a monomorphism.  Hence, as \(\sabs{f}\)
is an extremal epimorphism, \(\sabs{m}\) is an isomorphism.  Consider
the \dmor \(\sabs{m}^{-1} \sabs{f} \colon \sabs{\Qdobj_1} \to
\sabs{\Qdobj'}\).  This simplifies to \(\sabs{g}\) which lifts to a
\Qdmor.  Hence, by lemma~\ref{lem:pushlift} since \(\Qdobj_2 =
\cov{\sabs{f}}(\Qdobj_1)\), \(\sabs{m}^{-1}\) lifts to a \Qdmor and
thus, as the forgetful functor is faithful, \(m\) is an isomorphism.
Hence \(f\) is an extremal epimorphism.
\end{proof}

This characterisation helps us prove the required extremallity
properties of \(\Qdcat\). 

\begin{corollary}
\(\Qdcat\) is an (extremal epi, mono) category.
\end{corollary}

\begin{proof}
Let \(f \colon \Qdobj_1 \to \Qdobj_2\) be a \Qdmor.  The \dmor
\(\sabs{f}\) has an (extremal epi, mono)\hyp{}factorisation, say
\(\sabs{f} = m e\).  We can lift this to a factorisation of \(f\) as
\(\Qdobj_1 \to \cov{e}(\Qdobj_1) \to \Qdobj_2\).  This is an (extremal
epi, mono) factorisation by corollary~\ref{cor:qexepilift} and as
monomorphisms lift to monomorphisms.

To show uniqueness it is sufficient to show that we have the diagonal
property.  That is, suppose that we have a commutative square in
\(\Qdcat\),
\[
  \xymatrix{
  \Qdobj_1 \ar[r] \ar[d]^{e} &
    \Qdobj_2 \ar[d]^{m} \\
    \Qdobj_3 \ar[r] &
    \Qdobj_4
  }
\]
with \(e\) an extremal epimorphism and \(m\) a monomorphism.  The
underlying square in \(\dcat\) has the same properties and thus there
is a (unique) \dmor \(h \colon \sabs{\Qdobj_3} \to \sabs{\Qdobj_2}\)
which fits into the corresponding diagram in \(\dcat\).  Then \(h e\)
lifts to a \Qdmor so as \(\Qdobj_3 = \cov{e}(\Qdobj_1)\), \(h\) lifts
to a \Qdmor and thus \(\Qdcat\) has the (extremal epi,
mono)\hyp{}diagonalisation property.  Thus \(\Qdcat\) is an (extremal
epi, mono)\hyp{}category.
\end{proof}

\begin{corollary}
\(\Qdcat\) is extremally \co{well}\hyp{}powered.
\end{corollary}

\begin{proof}
Let \(\Qdobj\) be an \Qdobj.  From the characterisation of extremal
epimorphisms in \(\Qdcat\) we see that the forgetful functor \(\Qdcat
\to \dcat\) defines a bijection from the class of isomorphism classes
of extremal epimorphisms in \(\Qdcat\) with source \(\Qdobj\) to the
class of isomorphism classes of extremal epimorphisms in \(\dcat\)
with source \(\sabs{\Qdobj}\).  Hence the property of being extremally
\co{well}\hyp{}powered lifts from \(\dcat\) to \(\Qdcat\).
\end{proof}

From lemmas~\ref{lem:pushlift} and~\ref{lem:pulllift} we can deduce
that \(\Qdcat\) is both complete and \co{complete}. 

\begin{proposition}
\label{prop:qcomp}
\(\Qdcat\) is complete and \co{complete}.
\end{proposition}

\begin{proof}
This is a standard proof.  We form limits and \co{limits} in
\(\Qdcat\) by forming the limit or \co{limit} first in \(\dcat\) and
then putting the appropriate reduced filtration on the resulting
object: the pull back filtration for the limit and the push forward
for the \co{limit}.
\end{proof}

We therefore have all our required properties of \(\Qdcat\).  We now
turn to \(\Kdcat\).  Completeness and \co{completeness} follow
directly from proposition~\ref{prop:qcomp}.

\begin{corollary}
\label{cor:kcomp}
\(\Kdcat\) is complete and \co{complete}.
\end{corollary}

\begin{proof}
It is a reflective, full subcategory of \(\Qdcat\) which is closed
under isomorphism.  Hence by \cite[corollaries 36.14,18]{hhgs}, both
completeness and \co{completeness} descend from \(\Qdcat\) to
\(\Kdcat\).
\end{proof}

Note that \co{limits} in \(\Kdcat\) are not simply the \co{limits} of
the corresponding family in \(\Qdcat\) and therefore do not
necessarily project down to the corresponding \co{limit} in \(\dcat\).
Rather we form the \co{limit} in \(\Qdcat\) and then apply the functor
\(\comfunc\) to the resulting object.

Extremal epimorphisms in \(\Kdcat\) are more complicated than in
\(\Qdcat\) and so we need to work harder to prove that \(\Kdcat\) is
an (extremal epi, mono) category and is extremally
\co{well}\hyp{}powered.  

\begin{proposition}
\(\Kdcat\) is an (extremal epi, mono) category.
\end{proposition}

\begin{proof}
The proof of the factorisation property is an adaptation of the
standard proof that every morphism in a complete well\hyp{}powered
category is (extremal epi, mono)\hyp{}factorisable.

Let \(f \colon \Kdobj_1 \to \Kdobj_2\) be a \Kdmor.  We consider the
class of all factorisations \(f = m h\) with \(m\) a monomorphism.
This is not empty as it contains the factorisation \((1, f)\).
This is a quasi\hyp{}ordered class with \((m_1, h_1) \ge (m_2, h_2)\)
if there is a \Kdmor from the source of \(m_1\) to the source of
\(m_2\) making the obvious diagram commute.  If this morphism exists,
it is obviously unique and a monomorphism.  There is an obvious
functor from this quasi\hyp{}ordered class to \(\Kdcat\).

We wish to show that this functor has a limit.  We shall do this by
showing that the class has a small initial subclass.  Let
\(\sabs{f} = m_{\dobj} e_{\dobj}\) be the (extremal epi,
mono)\hyp{}factorisation of \(\sabs{f}\) in \(\dcat\) with intervening
\dobj[\dobj].  Let \(f = m h\) be a factorisation of \(f\) with \(m\)
a monomorphism and
intervening \Kdobj[\Kdobj].  Then \(\sabs{f} = \sabs{m}\sabs{h}\) is a
factorisation of \(\sabs{f}\).  As the forgetful functor \(\Kdcat \to
\dcat\) has a left adjoint (the discrete filtration functor) it takes
monomorphisms to monomorphisms and so \(\sabs{m}\) is a monomorphism.
Hence as \(\dcat\) is an (extremal epi, mono) category there is a
\dmor \(g \colon \dobj \to \sabs{\Kdobj}\) making the following
diagram commute.
\[
  \xymatrix{
    \sabs{\Kdobj_1} \ar[r]^{\sabs{h}} \ar[d]^{e_{\dobj}} &
    \sabs{\Kdobj} \ar[d]^{\sabs{m}} \\
    \dobj \ar[ur]^{g} \ar[r]^{m_{\dobj}} &
    \sabs{\Kdobj_2}
  }
\]

We pull back and reduce the projective filtration \(\Kdobj_2\) on
\(\sabs{\Kdobj_2}\) via \(m_{\dobj}\) to one on \(\dobj\); let us
write this as \(\Qdobj\).  The above diagram
then lifts to
\(\Qdcat\) with \(\Qdobj\) in the lower left corner.  Via the
adjunction \(\comfunc \adjoint \kiqfunc\), the \Qdmors with source
\(\Qdobj\) factor through the natural morphism \(\Qdobj \to
\kiqfunc\comfunc(\Qdobj)\).  We therefore have the diagram in
\(\Qdcat\)
\[
  \xymatrix{
    \kiqfunc(\Kdobj_1) \ar[rr]^{\kiqfunc(h)} \ar[d]^{\mQo{e_{\dobj}}}
    &&
    \kiqfunc(\Kdobj) \ar[d]^{\kiqfunc(m)} \\
    \Qdobj \ar[r]^{\iota} &
    \kiqfunc\comfunc(\Qdobj) \ar[ur]^{\kiqfunc(\mKo{g})}
    \ar[r]^{\kiqfunc(\mKo{m_{\dobj}})} &
    \kiqfunc(\Kdobj_2),
  }
\]
where \(\iota\) is the canonical morphism.

We claim that \(\mKo{g}\) is a monomorphism.  It is necessary and
sufficient to show that \(\sabs{\mKo{g}}\) is a monomorphism as the
forgetful functor is faithful and has a left adjoint.  Thus let
\(\delt_1, \delt_2 \colon \dobj' \to \sabs{\comfunc(\Qdobj)}\) be
\dmors such that \(\sabs{\mKo{g}}\delt_1 = \sabs{\mKo{g}}\delt_2\).  As
\(\sabs{\comfunc(\Qdobj)}\) is the underlying \dobj of an \Kdobj, it
is a limit and so \(\delt_1\) and \(\delt_2\) are completely
determined by their compositions with the morphisms into the
appropriate family.  This family is the projective filtration
\(\Qdobj\) and for \(\Qdelt\) in \(\Qdobj\) we have a \dmor
\(\tilde{\Qdelt} \colon \sabs{\comfunc(\Qdobj)} \to \Qdobj_{\Qdelt}\)
such that \(\tilde{\Qdelt} \sabs{\iota} = \Qdelt\).

The projective filtration \(\Qdobj\) was defined as the reduction
of the pull back of \(\Kdobj_2\) via \(m_{\dobj}\).  It
therefore has an initial family as follows: for each \(\Kdelt\) in
\(\Kdobj_2\) the fact that \(\dcat\) is an (extremal epi, mono) category
implies the existence of a commutative diagram, unique up to canonical
isomorphism,
\[
  \xymatrix{
    \dobj \ar[r]^{m_{\dobj}} \ar[d]^{\Kdelt_{\dobj}} &
    \sabs{\Kdobj_2} \ar[d]^{\Kdelt} \\
    \dobj_{\Kdelt} \ar[r]^{m_{\Kdelt}} &
    \Kdobj_{\Kdelt},
  }
\]
with \(\Kdelt_{\dobj}\) an extremal epimorphism and \(m_{\Kdelt}\) a
monomorphism.  The family \(\Kdelt_{\dobj}\) is initial for
\(\Qdobj\).  We therefore have the following commutative diagram. 
\[
  \xymatrix{
    \dobj \ar[r]^{\sabs{\iota}} \ar[rd]_{\Kdelt_{\dobj}} &
    \sabs{\comfunc(\Qdobj)} \ar[d]^{\tilde{\Kdelt}_{\dobj}}
    \ar[r]^{m\sabs{\mKo{g}}} &
    \sabs{\Kdobj_2} \ar[d]^{\Kdelt} \\
    & \dobj_{\Kdelt} \ar[r]^{m_{\Kdelt}} &
    \Kdobj_{\Kdelt}
  }
\]
Thus \(\delt_1\) and \(\delt_2\) are completely determined by the
compositions \(\tilde{\Kdelt}_{\dobj} \delt_i\) for \(\Kdelt\) in
\(\Kdobj_2\).  Now as they satisfy \(\sabs{\mKo{g}} \delt_1 =
\sabs{\mKo{g}} \delt_2\), for each \(\Kdelt\) in \(\Kdobj_2\) we have
\(\Kdelt m\sabs{\mKo{g}} \delt_1 = \Kdelt m\sabs{\mKo{g}} \delt_2\)
whence \(m_{\Kdelt} \tilde{\Kdelt}_{\dobj} \delt_1 = m_{\Kdelt}
\tilde{\Kdelt}_{\dobj} \delt_2\).  As \(m_{\Kdelt}\) is a
monomorphism, we therefore have \(\tilde{\Kdelt}_{\dobj} \delt_1 =
\tilde{\Kdelt}_{\dobj} \delt_2\) and thus \(\delt_1 = \delt_2\).
Hence \(\sabs{\mKo{g}}\) is a monomorphism and thus so is \(\mKo{g}\).

The \Kdmor \(\mKo{m_{\dobj}}\) is equal to \(m \mKo{g}\) and so is a
monomorphism.  Let \(\mKo{e_{\dobj}} \colon \Kdobj_1 \to
\comfunc(\Qdobj)\) be the \Kdmor which, under the inclusion \(\Kdcat
\subseteq \Qdcat\), maps to \(\iota \mQo{e_{\dobj}}\).  The
factorisation of \(f\) in \(\Kdcat\) as \(\mKo{m_{\dobj}}
\mKo{e_{\dobj}}\) is thus in our class of factorisations and it
preceeds the factorisation \(f = m h\).

The key property of this factorisation is that the morphism
\(\mKo{e_{\dobj}}\) is obtained by applying the functor \(\comfunc\)
to a \Qdmor \(\mQo{e_{\dobj}} \colon \kiqfunc(\Kdobj_1) \to \Qdobj\)
such that \(\sabs{\mQo{e_{\dobj}}} = e_{\dobj}\); note that \(e_{\dobj}\)
depends only on \(f\) and not on the factorisation that we were trying
to dominate.  This factorisation is therefore completely determined by
the reduced projective filtration \(\Qdobj\) on \(\dobj\).

As \(\dcat\) is extremally \co{well}\hyp{}powered, the class of all
reduced projective filtrations on a specified \dobj is actually a set.
To see this, observe that a reduced projective filtration is
completely determined by its subclass of extremal epimorphisms.
Moreover, this subclass is closed under isomorphism and so is a union
of equivalence classes of extremal epimorphisms, whence a reduced
projective filtration is completely determined by an element of the
power class of the class of equivalence classes of extremal
epimorphisms emanating from the original \dobj.  As \(\dcat\) is
extremally \co{well}\hyp{}powered, the class of equivalence classes of
extremal epimorphisms is actually a set and so its power set is also a
set.  Hence the class of reduced projective filtrations on a given
\dobj is a set.

Thus our class of factorisations of the \Kdmor[f] has an initial set
and so, as \(\Kdcat\) is complete, has a limit.  The proof that this
limit is an (extremal epi, mono)\hyp{}factorisation of \(f\) proceeds
exactly as in the analogous proof for a morphism in a complete
well\hyp{}powered category.  See, for example, \cite[17.8,17.16]{hhgs}.

The proof that \(\Kdcat\) is in fact an (extremal epi, mono)
category now follows since it is complete.  See, for example,
\cite[34.1]{hhgs}.
\end{proof}

It is worth pointing out that even if \(\dcat\) were, in fact, a
(regular epi, mono) category then it would not necessarily be true that
\(\Kdcat\) was a (regular epi, mono) category.

Buried within the above proof are all the necessary pieces to prove
the final property that we want.

\begin{proposition}
\(\Kdcat\) is extremally \co{well}\hyp{}powered.
\end{proposition}

\begin{proof}
It is easy to see from the proof that \(\Kdcat\) is an (extremal epi,
mono) category that every extremal epimorphism is obtained by applying
\(\comfunc\) to a \Qdmor of the form \(\kiqfunc(\Kdobj) \to \Qdobj\)
with underlying \dmor an extremal epimorphism (this is a necessary,
but not sufficient, condition).  To specify the isomorphism class of an
extremal epimorphism in \Kdcat it is therefore sufficient to specify
the isomorphism class of the corresponding extremal epimorphism in
\dcat and a reduced projective filtration on the target.

Hence for an \Kdobj[\Kdobj], the class of isomorphism classes of
extremal epimorphisms emanating from \(\Kdobj\) injects into
\[
  \coprod_{\text{Iso}(\text{ex epi } \delt \colon \sabs{\Kdobj} \to
    \dobj) } \text{Iso}(\Qdcat_{\dobj})
\]
where \(\text{Iso}(\text{ex epi } \delt \colon \sabs{\Kdobj} \to
\dobj)\) is the class of isomorphism classes of extremal epimorphisms
with source \(\sabs{\Kdobj}\) and \(\text{Iso}(\Qdcat_{\dobj})\) is
the class of isomorphism classes of the fibre category of \(\Qdcat \to
\dcat\) at \(\dobj\); that is, the class of reduced projective
filtrations on \(\dobj\).

As \(\dcat\) is extremally \co{well}\hyp{}powered, all of the classes
in this \co{product} are small.  Hence the \co{product} is small and
thus \(\Kdcat\) is extremally \co{well}\hyp{}powered.
\end{proof}

\subsection{The Canonical Filtration Functor}
\label{sec:filfunc}

In this section we shall construct a right adjoint to the forgetful
functor \(\kfofunc \colon \KKdcat \to \Kdcat\).  The indiscrete
filtration functor, and its reduction, provide right adjoints to the
forgetful functors \(\Pdcat \to \dcat\) and \(\Qdcat \to \dcat\) but
in general the forgetful functor \(\Kdcat \to \dcat\) does not have a
right adjoint.  In the specific case \(\KKdcat \to \Kdcat\), however,
we are able to construct one.  It is a straightforward adaptation of
the filtration functor for a filtered category as described in
example~\ref{ex:filtrations}(\ref{ex:filfunc}).  We assume that
\(\dcat\) has the properties of section~\ref{sec:catprop}.

It is simple to adapt the definition of
example~\ref{ex:filtrations}(\ref{ex:filfunc}) to define a functor
\(\Kdcat \to \PKdcat\).  For an \Pdobj[\Pdobj] the canonical
filtration was defined by taking the projective filtration on
\(\Pdobj\) with initial subclass the family of \Pdmors
\(\mPo{\;\Pdelt\,} \colon \Pdobj \to
\qipfunc\kiqfunc\disfunc(\Pdobj_{\Pdelt})\) (recall that we now regard
\(\disfunc\) as a functor into \(\Kdcat\)).  Similarly, we define a
functor \(\Kdcat \to \PKdcat\) by taking the projective filtration on
\(\Kdobj\) with initial subclass the family of \Kdmors \(\mKo{\Kdelt}
\colon \Kdobj \to \disfunc(\Kdobj_{\Kdelt})\).

\begin{proposition}
This functor factors through \(\KKdcat\).
\end{proposition}

\begin{proof}
We need to show first that the filtration defined above is reduced and
then that it is an iso\hyp{}filtration.

Let \(\Kdobj\) be an \Kdobj and let \(\PKdobj\) be the resulting
\PKdobj.  By construction, an initial subclass for the filtration
\(\PKdobj\) is given by taking the \Kdmors
\[
  \mKo{\Kdelt} \colon \Kdobj \to \disfunc(\Kdobj_{\Kdelt})
\]
for \(\Kdelt\) in \(\Kdobj\).  It is obvious that we may refine this
further and take \(\Kdelt\) in an initial subclass of \(\Kdobj\).

In particular, we can take those \(\Kdelt\) which are extremal
epimorphisms.  Let \(\Kdelt\) be one of these.  We wish to show that
\(\mKo{\Kdelt}\) is an extremal epimorphism. Since \(\Kdelt\) is an
epimorphism and the forgetful functor is faithful, \(\mKo{\Kdelt}\) is
an epimorphism. Now let \(\mKo{\Kdelt} = m f\) be a factorisation with
\(m\) a monomorphism.  Let \(\Kdobj'\) be the intervening \Kdobj.  By
applying the forgetful functor we obtain a factorisation of \(\Kdelt\)
as \(\sabs{m} \sabs{f}\).  The forgetful functor \(\Kdcat \to \dcat\)
has a left adjoint, namely the discrete filtration functor, so
preserves monomorphisms.  Hence \(\sabs{m}\) is a monomorphism and
thus, as \(\Kdelt\) is an extremal epimorphism, \(\sabs{m}\) is an
isomorphism.  Its inverse is a \dmor \(\Kdobj_{\Kdelt} \to
\sabs{\Kdobj'}\) and hence lifts to a \Kdmor
\(\disfunc(\Kdobj_{\Kdelt}) \to \Kdobj'\).  This lift is inverse to
\(m\) because the forgetful functor \(\Kdcat \to \dcat\) is faithful.
Hence \(\mKo{\Kdelt}\) is an extremal epimorphism and so \(\PKdobj\)
is a reduced projective filtration.

To show that it is an iso\hyp{}filtration we need to show that the
limit of \(\PKdobj\) is isomorphic to \(\Kdobj\) via the canonical
morphism.  This follows from the description of limits in \(\Kdcat\):
they are formed by taking the underlying limit in \(\dcat\) and
putting the reduced pull back filtration on the resulting \dobj.  In our case,
the resulting \dobj is (naturally isomorphic to) \(\sabs{\Kdobj}\) and
it is obvious that the reduced pull back filtration of the family
\(\sabs{\Kdobj} \to \sabs{\disfunc(\Kdobj_{\Kdelt})}\) is again
\(\Kdobj\).
\end{proof}

\begin{defn}
We shall refer to the functor defined above as the \emph{canonical
  filtration functor} and denote it by \(\kflfunc \colon \Kdcat \to
\KKdcat\).
\end{defn}

Let us now show that this functor is the required adjoint.

\begin{proposition}
The canonical filtration functor \(\kflfunc \colon \Kdcat \to
\KKdcat\) is right adjoint and right inverse to the forgetful functor
\(\kfofunc \colon \KKdcat \to \Kdcat\).
\end{proposition}

\begin{proof}
In this proof, \(\kfofunc\) will refer exclusively to the forgetful
functor \(\KKdcat \to \Kdcat\) and we will use the notation
\(\sabs{-}\) for the forgetful functor \(\Kdcat \to \dcat\).

By construction, \(\kfofunc\kflfunc\) is the identity functor on
\(\Kdcat\).  This provides the natural transformation which will be
the \co{unit} of the adjunction.

Let \(\KKdobj\) be an \KKdobj.  Both \(\KKdobj\) and
\(\kflfunc\kfofunc(\KKdobj)\) are iso\hyp{}filtrations on the same
underlying \Kdobj, namely \(\kfofunc(\KKdobj)\).  We shall show that
\(\kflfunc\kfofunc(\KKdobj) \subseteq \KKdobj\) as projective
filtrations on \(\kfofunc(\KKdobj)\).  By construction,
\(\kflfunc\kfofunc(\KKdobj)\) is the projective filtration on
\(\kfofunc(\KKdobj)\) with initial subclass
\[
  \mKo{\Kdelt} \colon \kfofunc(\KKdobj) \to \disfunc(\kfofunc(\KKdobj)_k)
\]
for \(\Kdelt\) in \(\kfofunc(\KKdobj)\).  Now \(\KKdobj\) is an
iso\hyp{}filtration on \(\kfofunc(\KKdobj)\) and so the canonical
morphism
\[
  \kfofunc(\KKdobj) \to \varprojlim_{\KKdelt} \KKdobj_{\KKdelt}
\]
is an isomorphism of \Kdobjs, where the limit is over \(\KKdelt\) in
\(\KKdobj\).  Limits in \(\Kdcat\) are formed by taking the
corresponding limit in \(\dcat\) and then putting the reduced pull
back filtration on the resulting object.  Thus an initial subclass of
\(\kfofunc(\KKdobj)\) consists of the extremal epimorphisms coming from
the (extremal epi, mono)\hyp{}factorisation of \dmors of the form
\[
  \sabs{\kfofunc(\KKdobj)} \xrightarrow{\sabs{\KKdelt}}
  \sabs{\KKdobj_{\KKdelt}} \xrightarrow{\KKdelt'} \KKdobj_{\KKdelt,
    \KKdelt'}
\]
for \(\KKdelt\) in \(\KKdobj\) and \(\KKdelt'\) in
\(\KKdobj_{\KKdelt}\).  We can assume that \(\KKdelt\) and
\(\KKdelt'\) are themselves in initial subclasses of their respective
filtrations and so we can assume that they are extremal epimorphisms
in their respective categories.

Let us show that the composition \(\KKdelt' \sabs{\KKdelt}\) is itself
an extremal epimorphism.  It is an epimorphism because \(\KKdelt'\)
and \(\sabs{\KKdelt}\) are both epimorphisms.  Let \(\KKdelt'
\sabs{\KKdelt} = m e\) be the (extremal epi, mono)\hyp{}factorisation
in \(\dcat\) of \(\KKdelt' \sabs{\KKdelt}\) with intervening
\dobj[\dobj].  By the above, \(e\) is in \(\kfofunc(\KKdobj)\) and thus
lifts to a \Kdmor \(\tilde{e} \colon \kfofunc(\KKdobj) \to
\disfunc(\dobj)\).  The \dmor \(m \colon \dobj \to \KKdobj_{\KKdelt,
\KKdelt'}\) lifts to a \Kdmor \(\disfunc(m) \colon \disfunc(\dobj) \to
\disfunc(\KKdobj_{\KKdelt, \KKdelt'})\).  This is again a monomorphism
as the forgetful functor is faithful.  Since \(\KKdelt'\) is in
\(\KKdobj_{\KKdelt}\) it also lifts to a \Kdmor \(\tilde{\KKdelt'}
\colon \KKdobj_{\KKdelt} \to \disfunc(\KKdelt_{\KKdelt, \KKdelt'})\).
As lifts of \dmors to \Kdmors are unique, we therefore have
the following
commutative diagram in \(\Kdcat\).
\[
  \xymatrix{
    \kfofunc(\KKdobj) \ar[r]^{\KKdelt} \ar[d]_{\tilde{e}} &
    \KKdobj_{\KKdelt} \ar[d]^{\tilde{\KKdelt}'} \\
    \disfunc(\dobj) \ar[r]^{\disfunc(m)} &
    \disfunc(\KKdobj_{\KKdelt,\KKdelt'})
  }
\]
Since \(\Kdcat\) is an (extremal epi, mono) category there is a \Kdmor
\(g \colon \KKdobj_{\KKdelt} \to \disfunc(\dobj)\) which fits into the
above diagram.  Applying the forgetful functor we see that \(\KKdelt'
= m \sabs{g}\).  As \(\KKdelt'\) is an extremal epimorphism we see
that \(m\) is an isomorphism.  Hence \(\KKdelt' \sabs{\KKdelt}\) is
isomorphic to \(e\) and thus is an extremal epimorphism.

Thus an initial subclass for \(\kfofunc(\KKdobj)\) is the family of
\dmors
\[
  \sabs{\kfofunc(\KKdobj)} \xrightarrow{\sabs{\KKdelt}}
  \sabs{\KKdobj_{\KKdelt}} \xrightarrow{\KKdelt'} \KKdobj_{\KKdelt,
    \KKdelt'}
\]
with \(\KKdelt\) in \(\KKdobj\) and \(\KKdelt'\) in
\(\KKdobj_{\KKdelt}\).  Thus an initial subclass of
\(\kflfunc\kfofunc(\KKdobj)\) consists of the family of \Kdmors
\[
   \sabs{\KKdobj} \xrightarrow{\KKdelt} \KKdobj_{\KKdelt}
  \xrightarrow{\tilde{\KKdelt'}} \disfunc(\KKdobj_{\KKdelt,
    \KKdelt'}).
\]
As each of these factors through \(\KKdelt\) it is in \(\KKdobj\).
Hence \(\kflfunc\kfofunc(\KKdobj) \subseteq \KKdobj\) and so the
identity \Kdmor on \(\kfofunc(\KKdobj)\) lifts to a \KKdmor \(\KKdobj \to
\kflfunc\kfofunc(\KKdobj)\).

These lifts fit together to define a natural transformation of
functors from the identity on \(\KKdcat\) to \(\kflfunc\kfofunc\):
all the necessary diagrams commute because they do in \(\Kdcat\).
This will be the unit of our adjunction.

Our functors, \(\kfofunc\) and \(\kflfunc\), are both lifts of the
identity functor on \(\Kdcat\) along the forgetful functor \(\KKdcat
\to \Kdcat\), one lifting the source and the other the target.  The
natural transformations, \(\kfofunc\kflfunc \to I\) and \(I \to
\kflfunc\kfofunc\), are both lifts of the identity natural
transformation \(I \to I\) on \(\Kdcat\).  Therefore for an
\KKdobj[\KKdobj] and \Kdobj[\Kdobj] the forgetful functor \(\kfofunc
\colon \KKdcat \to \Kdcat\) induces a commutative diagram of morphisms
of hom\hyp{}sets.
\[
  \xymatrix{
    \Hom{\Kdcat}{\kfofunc(\KKdobj)}{\Kdobj} \ar[r] \ar[d]^{=} &
    {\rule[-6pt]{0pt}{6pt}} \Hom{\KKdcat}{\KKdobj}{\kflfunc(\Kdobj)}
    \ar[r] \ar[d]^{\kfofunc} &
    \Hom{\Kdcat}{\kfofunc(\KKdobj)}{\Kdobj} \ar[d]^{=} \\
    \Hom{\Kdcat}{\kfofunc(\KKdobj)}{\Kdobj} \ar[r]^{=} &
    \Hom{\Kdcat}{\kfofunc(\KKdobj)}{\Kdobj} \ar[r]^{=} &
    \Hom{\Kdcat}{\kfofunc(\KKdobj)}{\Kdobj}
  }
\]
As the forgetful functor is faithful, the morphisms in the upper line
are isomorphisms and hence \(\kflfunc\) is right adjoint to
\(\kfofunc\).
\end{proof}

\subsection{Lifts of Functors}

In this section we shall examine certain lifts of functors involving
filtered categories.  The two lifts that we shall consider are 
described
in the next definition.  Recall that for a functor \(\func{G} \colon
\dcat \to \ecat\) we defined, in section~\ref{sec:projind}, a
corresponding functor \(\pfunc{G} \colon \Pdcat \to \Pecat\).

The conditions that we impose on our categories in the following
theorems are not minimal.  Recall that we write \(\kfofunc \colon
\KKdcat \to \Kdcat\) for the forgetful functor.

\begin{defn}
\label{def:isolift}
Let \(\dcat\) and \(\ecat\) be complete, \co{complete},
extremally \co{well}\hyp{}powered, (extremal epi,
mono) categories.  Let \(\func{G} \colon \Kdcat \to \ecat\) be a
covariant functor.  We define \(\kfunc{G} \colon \Kdcat \to \Kecat\)
by
\[
\kfunc{G} \coloneqq \comfunc \redfunc \pfunc{G}
\qipfunc \kiqfunc \kflfunc
\colon \Kdcat \xrightarrow{\kflfunc} \KKdcat \xrightarrow{\kiqfunc} \QKdcat
  \xrightarrow{\qipfunc} \PKdcat \xrightarrow{\pfunc{G}}
  \Pecat \xrightarrow{\redfunc} \Qecat \xrightarrow{\comfunc}
  \Kecat.
\]

Let \(\func{H} \colon \ecat \to \Kdcat\) be a covariant functor.  We
define \(\kfunc{H} \colon \Kecat \to \Kdcat\) by
\[
\kfunc{H} \coloneqq \kfofunc \comfunc \redfunc
\pfunc{H} \qipfunc \kiqfunc 
\colon \Kecat \xrightarrow{\kiqfunc} \Qecat \xrightarrow{\qipfunc} \Pecat
\xrightarrow{\pfunc{H}} \PKdcat \xrightarrow{\redfunc}
\QKdcat \xrightarrow{\comfunc} \KKdcat \xrightarrow{\kfofunc} \Kdcat.
\]
\end{defn}

We trust that there will be no confusion with using the same notation
for two different constructions.  Note that in these definitions we
have two different instances of various functors.  In the definition
of \(\kfunc{G}\) the inclusion functor \(\kiqfunc\) is
from \(\KKdcat\) to \(\QKdcat\) whereas in the definition of
\(\kfunc{H}\) it is from \(\Kecat\) to \(\Qecat\).  We
trust that this also will not cause confusion.

We wish to prove two results about these constructions.  The first
gives a condition whereby the first construction is associative.  The
second relates to adjunctions.

\begin{theorem}
\label{th:complift}
Let \(\dcat\), \(\ecat\), and \(\fcat\) be complete, \co{complete},
extremally \co{well}\hyp{}powered, (extremal epi, mono) categories.
Let \(\func{G} \colon \Kdcat \to \ecat\) and \(\func{H} \colon \Kecat
\to \fcat\) be covariant functors.  Then if \(\func{H}\) preserves
monomorphisms, there is a natural isomorphism of functors \(\Kdcat \to
\Kfcat\)
\[
  \mKo{\func{H}\kfunc{G}} \cong \kfunc{H} \kfunc{G}
\]
satisfying the obvious coherence for triples.
\end{theorem}

\begin{proof}
Let us expand out the two sides to make clear what we have to prove.
\begin{align*}
\mKo{\func{H}\kfunc{G}} &= \comfunc\redfunc \mPo{\func{H}\kfunc{G}}
\qipfunc\kiqfunc \kflfunc =  \comfunc\redfunc \pfunc{H} \pkfunc{G}
\qipfunc\kiqfunc \kflfunc, \\
\kfunc{H}\kfunc{G} &= \comfunc\redfunc \pfunc{H} \qipfunc \kiqfunc
\kflfunc \comfunc \redfunc \pfunc{G} \qipfunc \kiqfunc \kflfunc =
\comfunc \redfunc \pfunc{H} \qipfunc \kiqfunc \kflfunc \kfunc{G}.
\end{align*}
From this it is clear that the first step is to compare \(\pkfunc{G}
\qipfunc \kiqfunc \kflfunc\) with \(\qipfunc \kiqfunc \kflfunc
\kfunc{G}\).  These are functors \(\Kdcat \to \PKecat\).  Using the
fact that \(\sabs{\pkfunc{G}(-)} = \kfunc{G}(\sabs{-})\) we see that
both \(\pkfunc{G} \qipfunc\kiqfunc \kflfunc\) and \(\qipfunc \kiqfunc
\kflfunc \kfunc{G}\) are lifts of \(\kfunc{G}\) along the forgetful
functor \(\PKecat \to \Kecat\).  To compare these lifts we need to
examine the resulting filtrations on \(\kfunc{G}(\Kdobj)\) for an
\Kdobj[\Kdobj].

Firstly, let us establish some notation.  Applying \(\redfunc\) to an
\Pdobj[\Pdobj] does not change the underlying \dobj, it merely alters
the filtration.  An initial class for \(\redfunc(\Pdobj)\) is given by
taking the extremal epimorphisms which come from the (extremal epi,
mono)\hyp{}factorisations of the \dmors in \(\Pdobj\).  As before, for
\(\Pdelt \colon \sabs{\Pdobj} \to \Pdobj_{\Pdelt}\) in \(\Pdobj\) let
us write
\[
  \mQo{\;\Pdelt\;} \colon \sabs{\Pdobj} \to \mQo{\Pdobj_{\Pdelt}}
\]
for the corresponding extremal epimorphism.

Let \(\Kdobj\) be an \Kdobj.  Let us examine \(\kfunc{G}(\Kdobj)\).
An initial subclass of the filtration \(\kflfunc(\Kdobj)\)
is given by the family
\[
  \Kdobj \to \disfunc(\Kdobj_{\Kdelt})
\]
for \(\Kdelt\) in \(\Kdobj\).  Here and henceforth we will suppress
the label for the morphism as the notation rapidly becomes unwieldy;
in each case it will be the obvious morphism derived from \(\Kdelt\).

An initial subclass of the filtration \(\pfunc{G}\qipfunc \kiqfunc
\kflfunc(\Kdobj)\) is thus given by the family
\begin{equation}
\label{eq:pfunc}
  \func{G}(\Kdobj) \to
  \func{G}\disfunc(\Kdobj_{\Kdelt})
\end{equation}
and of \(\redfunc \pfunc{G} \qipfunc \kiqfunc \kflfunc(\Kdobj)\) by
\[
  \func{G}(\Kdobj) \to
  \mQo{\func{G}\disfunc(\Kdobj_{\Kdelt})}.
\]
To get \(\kfunc{G}(\Kdobj)\) we apply \(\comfunc\) which replaces the
source of these morphisms by the appropriate limit.  This produces an
initial family for \(\kfunc{G}(\Kdobj)\) consisting of the \PKemors
\[
  \sabs{\kfunc{G}(\Kdobj)} \to
  \mQo{\func{G}\disfunc(\Kdobj_{\Kdelt})}.
\]

We can read off from this an initial family for \(\qipfunc \kiqfunc
\kflfunc \kfunc{G}(\Kdobj)\).  It consists of the \Kemors
\[
  \kfunc{G}(\Kdobj) \to
  \disfunc(\mQo{\func{G}\disfunc(\Kdobj_{\Kdelt})}).
\]
By applying the above as far as \eqref{eq:pfunc} with \(\kfunc{G}\) in
place of \(\func{G}\) we can also read off an initial family for
\(\pkfunc{G} \qipfunc \kiqfunc \kflfunc(\Kdobj)\).  It consists of the
\Kemors
\[
  \kfunc{G}(\Kdobj) \to \kfunc{G} \disfunc(\Kdobj_{\Kdelt}).
\]
It is obvious from the construction that \(\kfunc{G} \disfunc =
\disfunc \func{G} \disfunc\).  Thus we can rewrite the above \Kemors as
\[
  \kfunc{G}(\Kdobj) \to \disfunc \func{G} \disfunc(\Kdobj_{\Kdelt}).
\]

The \eobjs[\mQo{\func{G}\disfunc(\Kdobj_{\Kdelt})}] are defined (up to
canonical isomorphism) by the (extremal epi, mono)\hyp{}factorisations
\[
  \func{G}(\Kdobj) \to \mQo{\func{G}\disfunc(\Kdobj_{\Kdelt})} \to
  \func{G}\disfunc(\Kdobj_{\Kdelt}).
\]
From the construction of \(\kfunc{G}\) as a limit (via the functor
\(\comfunc\)) we see that we can replace \(\func{G}(\Kdobj)\) by
\(\sabs{\kfunc{G}(\Kdobj)}\) in this.  Thus, as projective filtrations
on \(\kfunc{G}(\Kdobj)\),
\(\pkfunc{G}\qipfunc\kiqfunc\kflfunc(\Kdobj)\) is contained in
\(\qipfunc\kiqfunc\kflfunc\kfunc{G}(\Kdobj)\).  From this we deduce
that the identity on \(\kfunc{G}(\Kdobj)\) lifts to a \PKemor
\[
\qipfunc\kiqfunc\kflfunc\kfunc{G}(\Kdobj) \to
\pkfunc{G}\qipfunc\kiqfunc\kflfunc(\Kdobj).
\]
As this is a lift of the identity on \(\kfunc{G}\) it defines a
natural transformation \(\qipfunc\kiqfunc\kflfunc\kfunc{G} \to
\pkfunc{G}\qipfunc\kiqfunc\kflfunc\).

We claim that this becomes a natural isomorphism after reduction.
To prove this claim, consider the following diagrams in \(\Kecat\).
\[
  \xymatrix{
    \kfunc{G}(\Kdobj) \ar[r] \ar[d] &
    \disfunc \func{G} \disfunc (\Kdobj_{\Kdelt}) \\
    \disfunc \mQo{\func{G} \disfunc (\Kdobj_{\Kdelt})} \ar[ru]
  }
\]

\[
  \xymatrix{
    \kfunc{G}(\Kdobj) \ar[r] \ar[d] &
    \mQo{\disfunc \func{G} \disfunc (\Kdobj_{\Kdelt})} \ar[r] &
    \disfunc \func{G} \disfunc (\Kdobj_{\Kdelt}) \\
    \mQo{\disfunc \mQo{\func{G} \disfunc (\Kdobj_{\Kdelt})}} \ar[d] \\
    \disfunc \mQo{\func{G} \disfunc (\Kdobj_{\Kdelt})} \ar[rruu] 
  }
\]

The second is derived from the first by taking the (extremal epi,
mono)\hyp{}factorisations of the vertical and horizontal morphisms.
The diagonal arrow in each diagram is a monomorphism since
it is \(\disfunc\) applied to a monomorphism and \(\disfunc\) preserves
monomorphisms. Thus
by
the diagonalisation property of an (extremal epi, mono) category,
there is an isomorphism
\[
    \mQo{\disfunc \mQo{\func{G} \disfunc (\Kdobj_{\Kdelt})}} \cong
    \mQo{\disfunc \func{G} \disfunc (\Kdobj_{\Kdelt})} 
\]
fitting in to the second diagram.

Thus the natural transformation
\[
  \qipfunc\kiqfunc\kflfunc\kfunc{G}
\to \pkfunc{G}\qipfunc\kiqfunc\kflfunc
\]
induces a natural
\emph{isomorphism}
\[
  \redfunc \qipfunc\kiqfunc\kflfunc\kfunc{G} \to
\redfunc \pkfunc{G}\qipfunc\kiqfunc\kflfunc.
\]

This is still not quite what is needed as there is an occurence of the
functor \(\redfunc\) which is not present in the expansions of the two
functors that we wish to compare.  What we shall now see is that we
could easily insert \(\redfunc\) at the appropriate juncture without
changing the result; actually we insert \(\qipfunc\redfunc\).

There is a natural transformation \(\qipfunc\redfunc \to 1\) coming
from the adjunction \(\qipfunc \adjoint \redfunc\).  Applying
\(\redfunc\pfunc{H}\) yields a natural transformation \(\redfunc
\pfunc{H} \qipfunc\redfunc \to \redfunc \pfunc{H}\).  We wish to show
that this is a natural isomorphism.  Let \(\PKeobj\) be an \PKeobj.
Using notation as above, we have the following initial classes of the
projective filtrations
\begin{align*}
 \big( 
\func{H}(\sabs{\PKeobj}) \to \mQo{\func{H}(\PKeobj_{\PKeelt})}
\big) & \text{\hspace{1em}for\hspace{1em}} \redfunc \pfunc{H} (\PKeobj), \\
\big( \func{H}(\sabs{\PKeobj}) \to \mQo{\func{H}(\mQo{\PKeobj_{\PKeelt}})}
\big) & \text{\hspace{1em}for\hspace{1em}} \redfunc \pfunc{H} \qipfunc\redfunc (\PKeobj)
\end{align*}
and the natural transformation is the identity on
\(\func{H}(\sabs{\PKeobj})\).  This natural transformation comes from
the fact that the diagonalisation property of an (extremal epi, mono)
category allows us to add in the required morphism on the inner
diagonal of this diagram:
\[
  \xymatrix{
    \func{H}(\sabs{\PKeobj}) \ar[r] \ar[d] &
    \mQo{\func{H}(\PKeobj_{\PKeelt})} \ar[r] &
    \func{H}(\PKeobj_{\PKeelt}) \\
    \mQo{\func{H}(\mQo{\PKeobj_{\PKeelt}})} \ar[d] \ar[ur] \\
    \func{H}(\mQo{\PKeobj_{\PKeelt}}) \ar[uurr].
  }
\]
By assumption, \(\func{H}\) preserves monomorphisms and thus the outer
diagonal is a monomorphism.  Using the diagonalisation property again,
the morphism on the inner diagonal is therefore an isomorphism.  Thus
the identity on \(\func{H}(\sabs{\PKeobj})\) underlies an
\(\Qfcat\)\hyp{}isomorphism \(\redfunc \pfunc{H} \qipfunc \redfunc
(\PKeobj) \to \redfunc \pfunc{H} (\PKeobj)\) and we therefore have the
required natural isomorphism.

Hence we have natural isomorphisms
\begin{align*}
  \mKo{\func{H}\kfunc{G}} &= \comfunc \redfunc \pfunc{H}
  \mPo{\kfunc{G}} \qipfunc \kiqfunc \kflfunc \\
& \cong \comfunc \redfunc
  \pfunc{H} \qipfunc \redfunc \mPo{\kfunc{G}} \qipfunc \kiqfunc
  \kflfunc \\
& \cong \comfunc \redfunc \pfunc{H} \qipfunc \redfunc
  \qipfunc\kiqfunc\kflfunc \kfunc{G} \\
& = \comfunc \redfunc \pfunc{H}
  \qipfunc \kiqfunc \kflfunc \kfunc{G} \\
& = \kfunc{H} \kfunc{G}
\end{align*}
as required.

That this satisfies the required coherence axiom is a simple
deduction from the fact that the natural transformations involved in
the above all derive from lifting identity morphisms along faithful
functors.  Thus one can faithfully map the required diagram
down to one which obviously commutes.
\end{proof}

The constructions of definition~\ref{def:isolift} extend in the
obvious way to morphisms (that is, natural transformations). Thus the
first construction gives a functor from the category of covariant functors
from \(\Kdcat\) to \(\ecat\) to the category of covariant functors
from \(\Kdcat\) to \(\Kecat\), which on objects is given by \(\func{G}
\mapsto \kfunc{G}\). Similarly, the second construction gives a functor
from the category of covariant functors
from \(\ecat\) to \(\Kdcat\) to the category of covariant functors
from \(\Kecat\) to \(\Kdcat\), which on objects is given by
\(\func{H}\mapsto \kfunc{H}\).

\medskip

The rest of this section in concerned with proving the following
result on adjunctions.

\begin{theorem}
\label{th:isoadj}
Let \(\dcat\) and \(\ecat\) be complete, \co{complete}, extremally
\co{well}\hyp{}powered, (extremal epi, mono) categories.  Let
\(\func{G} \colon \Kdcat \to \ecat\) and \(\func{H} \colon \ecat \to
\Kdcat\) be covariant functors such that \(\func{H}\) is left adjoint
to \(\func{G}\).  If \(\func{H}\) preserves extremal epimorphisms then
\(\kfunc{H}\) is left adjoint to \(\kfunc{G}\).
\end{theorem}

We shall prove this theorem in stages.  We shall work
``functor\hyp{}by\hyp{}functor'' using the definitions.  If we read
\(\kfunc{G}\) from right to left and
\(\kfunc{H}\) from left to right then we see that the
functors in the definitions pair\hyp{}up as adjoint pairs.  For four
of the six functors these adjunctions are the correct way around.
Two, however, are the wrong way round.  Using \(\kfunc{G}\)
as our primary reference these are the \(\comfunc\) and the
\(\qipfunc\).  Therefore we shall need to prove that these functors
are ``ineffective''.

\begin{lemma}
\label{lem:isoirr}
There is a natural isomorphism of covariant functors \(\Kdcat \to
\Qecat\)
\[
  \kiqfunc \comfunc \redfunc \pfunc{G} \qipfunc \kiqfunc
  \kflfunc \cong \redfunc \pfunc{G} \qipfunc
  \kiqfunc \kflfunc.
\]
\end{lemma}

\begin{proof}
Let \(\Kdobj\) be an \Kdobj and consider \(\pfunc{G} \qipfunc \kiqfunc
\kflfunc (\Kdobj)\).  The functors \(\qipfunc\) and \(\kiqfunc\) are
simply the inclusions of full subcategories and so the resulting
\Peobj has initial class
\[
  \func{G}(\Kdobj) \to
  \func{G}\disfunc(\Kdobj_{\Kdelt}).
\]
Using our earlier notation, the reduction of this has initial class
\[
  \func{G}(\Kdobj) \to
  \mQo{\func{G}\disfunc(\Kdobj_{\Kdelt})}.
\]
We claim that this \Qeobj is actually an \Keobj.

As \(\kflfunc(\Kdobj)\) is an \KKdobj, the natural morphism
\[
  \Kdobj \to \varprojlim_{\Kdelt} \disfunc(\Kdobj_{\Kdelt})
\]
is an isomorphism.  Since \(\func{G}\) has a left adjoint, it
preserves limits and so the limit in \(\ecat\) exists and the natural
morphism
\[
  \func{G}(\Kdobj) \to \varprojlim_{\Kdelt}
  \func{G}\disfunc(\Kdobj_{\Kdelt})
\]
is an isomorphism.  To complete the proof we just need to show that
this limit is unchanged upon passing to the reduction.

Let \(\{\eobj_\mu\}\) be a directed family of \eobjs with limit
\(\eobj\) and canonical morphisms \(\eelt_\mu \colon \eobj \to
\eobj_\mu\).  For each \(\mu\) in the indexing set, let \(\eelt_\mu =
m_\mu \eelt'_\mu\) be the (extremal epi,
mono)\hyp{}factorisation of \(\eelt_\mu\) with intervening object
\(\tilde{\eobj}_\mu\).  As \(\ecat\) is an (extremal epi, mono)
category the family \(\{\tilde{\eobj}_\mu\}\) is directed under the
same ordering as on the original family.  Let
\(\tilde{\eobj}\) be the limit of this family with \(\tilde{\eelt}_\mu
\colon \tilde{\eobj} \to \tilde{\eobj}_\mu\) the canonical morphisms.

The morphisms \(\eelt'_\mu \colon \eobj \to \tilde{\eobj}_\mu\) define
a morphism \(\eelt' \colon \eobj \to \tilde{\eobj}\) such that
\(\tilde{\eelt}_\mu \eelt' = \eelt'_\mu\).  The morphisms
\(m_\mu \tilde{\eelt}_\mu \colon \tilde{\eobj} \to \eobj_\mu\) define
a morphism \(m \colon \tilde{\eobj} \to \eobj\) such that \(\eelt_\mu
m = m_\mu \tilde{\eelt}_\mu\).

Consider \(m \eelt' \colon \eobj \to \eobj\).  As \(\eobj\) is the
limit of the family \(\{\eobj_\mu\}\) this morphism is completely
determined by the composisions \(\eelt_\mu m \eelt'\).  This
simplifies as follows.
\[
  \eelt_\mu m \eelt' = m_\mu \tilde{\eelt}_\mu \eelt' = m_\mu
  \eelt'_\mu = \eelt_\mu.
\]
Hence \(m \eelt'\) is the identity on \(\eobj\).  Conversely, consider
\(\eelt' m \colon \tilde{\eobj} \to \tilde{\eobj}\).  Similarly, this
is completely determined by \(\tilde{\eelt}_\mu \eelt' m\) which in turn is
determined by \(m_\mu \tilde{\eelt}_\mu \eelt' m\)  as \(m_\mu\) is
a monomorphism.  This simplifies as follows.
\[
  m_\mu \tilde{\eelt}_\mu \eelt' m = m_\mu \eelt'_\mu m =
  \eelt_\mu m = m_\mu \tilde{\eelt}_\mu.
\]
Hence \(\eelt' m\) is the identity on \(\tilde{\eobj}\).  Thus
\(\eobj\) and \(\tilde{\eobj}\) are naturally isomorphic.

Applying this to the case in hand, we see that the natural morphism
\[
  \func{G}(\Kdobj) \to \varprojlim_{\Kdelt}
  \mQo{\func{G}\disfunc(\Kdobj_{\Kdelt})}
\]
is an isomorphism and hence \(\redfunc \pfunc{G} \qipfunc \kiqfunc
\kflfunc(\Kdobj)\) is actually an \Keobj, albeit viewed as an \Qeobj.
Hence applying \(\kiqfunc\comfunc\) to this object results in the same
object again.  A similar result holds for morphisms and thus there is
a natural isomorphism
\[
  \kiqfunc \comfunc \redfunc \pfunc{G} \qipfunc \kiqfunc
  \kflfunc \cong \redfunc \pfunc{G} \qipfunc
  \kiqfunc \kflfunc
\]
as required.
\end{proof}

\begin{corollary}
\label{cor:ladjpart}
There is a natural isomorphism of bifunctors \(\Kecat \times \Kdcat
\to \scat\)
\[
  \Hom{\Kecat}{-}{\kfunc{G}(-)} \cong
  \Hom{\PKdcat}{\pfunc{H} \qipfunc \kiqfunc(-)}{\qipfunc
    \kiqfunc \kflfunc(-)}.
\]
\end{corollary}

\begin{proof}
Lemma~\ref{lem:isoirr} provides us with the crucial step in the
following chain of natural isomorphisms.
\begin{align*}
\Hom{\Kecat}{-}{\kfunc{G}(-)} &=
\Hom{\Kecat}{-}{\comfunc \redfunc \pfunc{G} \qipfunc
  \kiqfunc \kflfunc(-)} && \text{definition of
  \(\kfunc{G}\),} \\
&=\Hom{\Qecat}{\kiqfunc(-)}{\kiqfunc \comfunc \redfunc
  \pfunc{G} \qipfunc \kiqfunc \kflfunc(-)} &&
\text{\(\Kecat\) full subcategory of \(\Qecat\),}\\
&\cong \Hom{\Qecat}{\kiqfunc(-)}{ \redfunc \pfunc{G}
\qipfunc \kiqfunc \kflfunc(-)} && \text{by lemma~\nmref{lem:isoirr},}\\
&\cong \Hom{\PKdcat}{\pfunc{H} \qipfunc \kiqfunc (-)}{\qipfunc
  \kiqfunc \kflfunc(-)} && \text{by the adjuntion pairings} \\
&&& \qipfunc \adjoint \redfunc \text{ and } \pfunc{H} \adjoint
  \pfunc{G} \\
&&& \text{(the latter by proposition~\nmref{prop:padjoint})}. \qedhere
\end{align*}
\end{proof}

We now reach the second point where the adjunction is the wrong way
around: we have \(\qipfunc \adjoint \redfunc\) but \(\qipfunc\)
currently appears on the \emph{right}.  When dealing with the
problematic \(\comfunc\) our strategy was to show that it was
essentially not there.  Our strategy here is instead to show that
there is an ``invisible'' \(\redfunc\).

\begin{lemma}
\label{lem:alredy}
There is a natural isomorphism of covariant functors \(\Kecat \to
\PKdcat\)
\[
  \pfunc{H} \qipfunc \kiqfunc \cong \qipfunc \redfunc
  \pfunc{H} \qipfunc \kiqfunc.
\]
\end{lemma}

\begin{proof}
Let \(\Keobj\) be an \Keobj.  As \(\kiqfunc\) and \(\qipfunc\) are
simply the inclusions of full subcategories,
\(\pfunc{H}\qipfunc\kiqfunc(\Keobj)\) is the \PKdobj with initial class
\[
  \func{H}(\sabs{\Keobj}) \to
  \func{H}(\Keobj_{\Kdelt}).
\]
Since \(\Keobj\) is an \Keobj, it has an initial class of \emors
\(\sabs{\Keobj} \to \Keobj_{\Kdelt}\) which are extremal epimorphisms.
Hence by assumption on \(\func{H}\), the corresponding morphisms
\(\func{H}(\sabs{\Keobj}) \to \func{H}(\Keobj_{\Kdelt})\) are also
extremal epimorphisms.  Thus the filtered object
\(\pfunc{H}\qipfunc\kiqfunc(\Keobj)\) is already reduced.  That is to
say, there is a natural isomorphism \(\qipfunc \redfunc \pfunc{H}
\qipfunc \kiqfunc \cong \pfunc{H} \qipfunc \kiqfunc\).
\end{proof}

\begin{corollary}
There is a natural isomorphism of bifunctors \(\Kecat \times \Kdcat
\to \scat\)
\[
  \Hom{\Kecat}{-}{\kfunc{G}(-)} \cong
  \Hom{\Kdcat}{\kfunc{H}(-)}{-}.
  \]
That is to say, \(\kfunc{H} \adjoint
\kfunc{G}\).
\end{corollary}

\begin{proof}
Lemma~\ref{lem:alredy} provides us with the crucial step in the
following chain of natural isomorphisms.
\begin{align*}
  \Hom{\Kecat}{-}{\kfunc{G}(-)} &\cong
  \Hom{\PKdcat}{\pfunc{H} \qipfunc \kiqfunc(-)}{\qipfunc
    \kiqfunc \kflfunc(-)} && \text{by corollary~\nmref{cor:ladjpart}} \\
&\cong
\Hom{\PKdcat}{\qipfunc \redfunc \pfunc{H} \qipfunc
  \kiqfunc(-)}{\qipfunc \kiqfunc \kflfunc(-)} && \text{by
  lemma~\nmref{lem:alredy}}\\
&\cong \Hom{\QKdcat}{\redfunc \pfunc{H} \qipfunc
  \kiqfunc(-)}{\kiqfunc \kflfunc(-)} &&\text{\(\QKdcat\) full
    subcategory of \(\PKdcat\)} \\
&\cong  \Hom{\Kdcat}{\kfofunc
    \comfunc \redfunc \pfunc{H} \qipfunc
  \kiqfunc(-)}{-} &&\text{by the adjunction pairings} \\
&&& \comfunc \adjoint \kiqfunc \text{ and } \kfofunc \adjoint \kflfunc
  \\
&= \Hom{\Kdcat}{\kfunc{H}(-)}{-} &&\text{by definition of
    \(\kfunc{H}\)}. \qedhere
\end{align*}
\end{proof}

This completes the proof of theorem~\ref{th:isoadj}.

\subsection{Inductive Filtrations}

We have focused in this section on projective filtrations; that is,
filtrations defined by morphisms out of a particular object.  Using
the standard categorical idea of duality we can 'op' everything in
this section to obtain results on filtrations defined by morphisms
\emph{into} a particular object.  We call these filtrations
\emph{inductive} filtrations.  Our main interest in inductive
filtrations arises from the obvious fact that a contravariant functor
takes an inductive filtration to a projective one; therefore if we
have a contravariant functor and wish to end up with projectively
filtered objects then the correct starting point is inductively
filtered ones.  All of our results dualise with no difficulty and we
shall use them in the following without further comment.  We shall
denote the inductive versions of the various categories of
projectively filtered objects by the corresponding left\hyp{}pointing
arrow.  Thus \(\Idcat\) is \Idcat, \(\Jdcat\) is \Jdcat, and
\(\Ldcat\) is \Ldcat.

%% file: hopf.tall.tex
\section{Filtered Tall\hyp{}Wraith Monoids}
\label{sec:filtall}

We now bring together the work of the previous two sections; that is,
we explain the interaction between Tall\hyp{}Wraith monoids and
filtrations.  Let \(\Gvcat\) be a variety of graded algebras (with
grading set \(\iset\)) and \(\dcat\) a complete, \co{complete},
extremally \co{well}\hyp{}powered, (extremal epi, mono) category.  Let
\(\KdCGvcat\) be \KdCGvcat.  \KdCGvobjus represent covariant functors
\(\Kdcat \to \Gvcat\).  Recall from theorem~\ref{th:gvcatprop} that
\(\Gvcat\) is also a complete, \co{complete}, extremally
\co{well}\hyp{}powered, (extremal epi, mono) category and so we can
lift these to functors \(\Kdcat \to \KGvcat\) using the first
construction in definition~\ref{def:isolift}.  In particular we can
take \(\dcat = \Gvcat\) whence we obtain functors \(\KGvcat \to
\KGvcat\).

The category of all covariant functors \(\KGvcat \to \KGvcat\) has an
obvious monoidal structure coming from composition which we would like
to transfer to \(\KGvCGvcat\).  In the more general situation, we can
consider the composition of a functor \(\Kdcat \to \KGvcat\) with a
functor \(\KGvcat \to \KGvcat\) and we would like to show that this
transfers to a pairing \(\KdCGvcat \times \KGvCGvcat \to \KdCGvcat\).
The groundwork of the previous section allows us to do this.

\begin{defn}
\label{def:comprep}
We say that a covariant functor \(\Kdcat \to \KGvcat\) is
\emph{representable} if it is isomorphic to a functor of the form
\(\Kcov{\KdCGvobj}\) for a \KdCGvobj[\KdCGvobj].  We write
\(\covrep(\Kdcat,\KGvcat)\) for the category of such functors.
\end{defn}

It follows from the definition that the lift functor given by the
first construction in definition~\ref{def:isolift}, \(\covrep(\Kdcat,
\Gvcat) \to \covrep(\Kdcat, \KGvcat)\), is dense.

\begin{proposition}
The functor \(\mOo{(\KdCGvcat)} \to \covrep(\Kdcat,\KGvcat)\), 
given on objects by
\(\KdCGvobj \mapsto \Kcov{\KdCGvobj}\), is an equivalence of
categories.
\end{proposition}

\begin{proof}
The standard Yoneda argument gives us an equivalence between the
categories \(\mOo{(\KdCGvcat)}\) and \(\covrep(\Kdcat, \Gvcat)\) so we
focus on the lift functor
\[
  \covrep(\Kdcat, \Gvcat) \to \covrep(\Kdcat, \KGvcat). 
\]

We define a functor \(\covrep(\Kdcat,\KGvcat) \to \covrep(\Kdcat,
\Gvcat)\) by post composition with the forgetful functor \(\kfofunc
\KGvcat \to \Gvcat\).  Recall that \(\kfofunc = \forfunc \qipfunc
\kiqfunc\).

The category \(\Kdcat\) is \co{complete} by corollary~\ref{cor:kcomp}
and so by corollary~\ref{cor:adjrepgrade} a functor, say \(\func{G}\),
in \(\covrep(\Kdcat,\Gvcat)\) has a left adjoint.  Thus
lemma~\ref{lem:isoirr} applies and after
composing with \(\forfunc
\qipfunc\) we have a natural isomorphism
\[
\kfofunc \kfunc{G} = \forfunc\qipfunc\kiqfunc\kfunc{G} \cong
  \forfunc\qipfunc\redfunc \pfunc{G} \qipfunc\kiqfunc\kflfunc.
\]
Since \(\redfunc\) does not change the underlying \Gvobj,
\(\forfunc\qipfunc\redfunc = \forfunc\).  Hence there is a natural
isomorphism
\[
  \kfofunc\kfunc{G} \cong
  \forfunc\pfunc{G}\qipfunc\kiqfunc\kflfunc.
\]
From the definition of \(\pfunc{G}\) we see that \(\forfunc\pfunc{G}
= \func{G}\forfunc\) and from the definition of \(\kflfunc\) we see
that \(\forfunc\qipfunc\kiqfunc\kflfunc\) is the identity on
\(\KGvcat\).  Thus we have a natural isomorphism
\[
  \kfofunc\kfunc{G} \cong \func{G}.
\]
A similar story occurs on morphisms (i.e.~natural transformations) and
so it is clear that the composition
\[
  \covrep(\KGvcat,\Gvcat) \to \covrep(\KGvcat,\KGvcat) \to
  \covrep(\KGvcat,\Gvcat)
\]
is naturally isomorphic to the identity on \(\covrep(\KGvcat,
\Gvcat)\).  The functor
\[
  \covrep(\KGvcat,\KGvcat)
\to \covrep(\KGvcat, \Gvcat)
\]
is therefore dense.  Using the fact that \(\covrep(\KGvcat, \Gvcat)
\to \covrep(\KGvcat, \KGvcat)\) is dense, we also see that the above
functor is full.  Furthermore, it is faithful because the forgetful
functor, \(\forfunc \colon \KGvcat \to \Gvcat\), is faithful.  Hence
it is an equivalence of categories.
\end{proof}

Exactly as previously, we wish to use this equivalence to define a
monoidal structure on \(\KGvCGvcat\) which corresponds to composition of
functors. 

\begin{theorem}
\label{th:kcomprep}
Under the stated assumptions on \(\dcat\), the composition of
representable functors, in the sense of definition~\ref{def:comprep},
is representable.
\end{theorem}

The proof of this depends on the results of the previous section.  We
therefore need to show that they apply.  Let \(\KdCGvobj\) be a
\KdCGvobj.  This represents the covariant functor
\(\cov{\KdCGvobj}\) and this functor has a left adjoint
\(\ladj{\KdCGvobj}\) by corollary~\ref{cor:adjrepgrade}.

\begin{lemma}
\label{lem:exepialg}
The functor \(\ladj{\KdCGvobj} \colon \Gvcat \to
\Kdcat\) preserves extremal epimorphisms.
\end{lemma}

\begin{proof}
Let \(f \colon \Gvobj_1 \to \Gvobj_2\) be an extremal epimorphism in
\(\Gvcat\).  As \(\ladj{\KdCGvobj}\) is a  left adjoint it takes
epimorphisms to epimorphisms and thus \(\ladj{\KdCGvobj}(f)\) is an
epimorphism.  The crucial part is therefore to show that it is
extremal.

To prove this we need to examine the definition of \(\ladj{\KdCGvobj}\)
from theorem~\ref{th:adjgrade}.  It is defined using the left adjoint of
the underlying functor 
\({\sabs{\cov{\KdCGvobj}}} \colon \Kdcat \to\Zscat\); 
let us write this adjoint as \(\ladj{\sabs{\KdCGvobj}}\).
The relationship is that there is a \co{equaliser} sequence
in \(\Kdcat\), natural
in \(\Gvobj\) in \(\Gvcat\),
\[
  \xymatrix{
    \ladj{\sabs{\KdCGvobj}}(\sabs{\free{Gv}(\sabs{\Gvobj})})
    \ar@<2pt>[r]^(.6){r_{\Gvobj}} \ar@<-2pt>[r]_(.6){s_{\Gvobj}} &
    \ladj{\sabs{\KdCGvobj}}(\sabs{\Gvobj}) \ar[r]^{p_{\Gvobj}} &
    \ladj{\KdCGvobj}(\Gvobj).
  }
\]

Therefore for a \(\Gvcat\)\hyp{}morphism \(f \colon \Gvobj_1 \to
\Gvobj_2\) we have a diagram in \(\Kdcat\) with the obvious commuting
properties:
\[
  \xymatrix{
    \ladj{\sabs{\KdCGvobj}}(\sabs{\free{Gv}(\sabs{\Gvobj_1})})
    \ar@<2pt>[r]^(.6){r_1} \ar@<-2pt>[r]_(.6){s_1}
    \ar[d]^{\ladj{\sabs{\KdCGvobj}}(\free{Gv}(\sabs{f}))} &
    \ladj{\sabs{\KdCGvobj}}(\sabs{\Gvobj_1}) \ar[r]^{p_1}
    \ar[d]^{\ladj{\sabs{\KdCGvobj}}(\sabs{f})} &
    \ladj{\KdCGvobj}(\Gvobj_1) \ar[d]^{\ladj{\KdCGvobj}(f)} \\
    \ladj{\sabs{\KdCGvobj}}(\sabs{\free{Gv}(\sabs{\Gvobj_2})})
    \ar@<2pt>[r]^(.6){r_2} \ar@<-2pt>[r]_(.6){s_2} &
    \ladj{\sabs{\KdCGvobj}}(\sabs{\Gvobj_2}) \ar[r]^{p_2} &
    \ladj{\KdCGvobj}(\Gvobj_2).
  }
\]    

We know that \(\ladj{\KdCGvobj}(f)\) is an epimorphism.  Suppose that
we have a factorisation \(\ladj{\KdCGvobj}(f) = m e\) with \(m
\colon \Kdobj \to \ladj{\KdCGvobj}(\Gvobj_2)\) a monomorphism.  As
\(f\) is an extremal epimorphism, theorem~\ref{th:gvcatprop} says that
\(\sabs{f} \colon \sabs{\Gvobj_1} \to \sabs{\Gvobj_2}\) is an
epimorphism on the underlying \Zsobjs and thus admits a section, say
\(g \colon \sabs{\Gvobj_2} \to \sabs{\Gvobj_1}\).  Hence there is a
\(\Kdcat\)\hyp{}morphism \(h \colon
\ladj{\sabs{\KdCGvobj}}(\sabs{\Gvobj_2}) \to \Kdobj\) given by
\[
  \ladj{\sabs{\KdCGvobj}}(\sabs{\Gvobj_2})
  \xrightarrow{\ladj{\sabs{\KdCGvobj}}(g)}
  \ladj{\sabs{\KdCGvobj}}(\sabs{\Gvobj_1}) \xrightarrow{p_1}
  \ladj{\KdCGvobj}(\Gvobj_1) \xrightarrow{e} \Kdobj.
\]
Then
\[
  m h = m e p_1 \ladj{\sabs{\KdCGvobj}}(g) = \ladj{\KdCGvobj}(f) p_1
  \ladj{\sabs{\KdCGvobj}}(g) = p_2 \ladj{\sabs{\KdCGvobj}}(\sabs{f})
  \ladj{\sabs{\KdCGvobj}}(g) = p_2
\]
and so we have a factorisation of \(p_2\) as \(m h\).  Thus
\(m h r_2 = m h s_2\) and so, as \(m\) is a monomorphism, \(h r_2 = h
s_2\).  Thus \(h\) \co{equalises} \(r_2\) and \(s_2\) so factors
uniquely through the \co{equaliser} which is
\(\ladj{\KdCGvobj}(\Gvobj_2)\).  That is, there is a unique
\(\Kdcat\)\hyp{}morphism \(k \colon \ladj{\KdCGvobj}(\Gvobj_2) \to
\Kdobj\) such that \(k p_2 = h\).  Then \(m k p_2 = m h = p_2\)
whence, as \(p_2\) is an epimorphism, \(m k\) is the identity on
\(\ladj{\KdCGvobj}(\Gvobj_2)\).  Thus \(m\) is a monomorphism which
admits a section and hence is an isomorphism.

In conclusion, \(\ladj{\KdCGvobj}(f)\) is an extremal epimorphism.
\end{proof}

Thus theorem~\ref{th:isoadj} applies to the pair
\((\ladj{\KdCGvobj}, \cov{\KdCGvobj})\).

\begin{corollary}
\label{cor:liftadj}
The functor \(\mKo{\ladj{\KdCGvobj}}\) is left adjoint to
\(\Kcov{\KdCGvobj}\). \noproof
\end{corollary}

We extend the functor \(\mKo{\ladj{\KdCGvobj}} \colon \KGvcat \to
\Kdcat\) to \(\mZo{\mKo{\ladj{\KdCGvobj}}} \colon \ZKGvcat \to
\ZKdcat\) in the obvious way.  As \(\mKo{\ladj{\KdCGvobj}}\) is a left
adjoint it preserves \co{products} and thus so does
\(\mZo{\mKo{\ladj{\KdCGvobj}}}\).  This therefore lifts to a functor
\(\mCo{\mKo{\ladj{\KdCGvobj}}} \colon \KGvCGvcat \to \KdCGvcat\).

\begin{corollary}
\label{cor:repkv}
Let \(\KdCGvobj\) be a \KdCGvobj and \(\KGvCGvobj\) a \KGvCGvobj.  The
composition
\[
  \cov{\KGvCGvobj} \Kcov{\KdCGvobj} \colon \Kdcat \to \Gvcat
\]
is representable with representing object
\(\mCo{\mKo{\ladj{\KdCGvobj}}}(\KGvCGvobj)\).
\end{corollary}

As we are considering a functor into \(\Gvcat\) we are using the usual
meaning of ``representable'' here.

\begin{proof}
Both \(\cov{\KGvCGvobj}\) and \(\Kcov{\KdCGvobj}\) have left
adjoints and so their composition has a left adjoint.  It is therefore
representable by corollary~\ref{cor:adjrepgrade}.  That the
representing object is as given is straightforward to see.
\end{proof}

As \(\cov{\KGvCGvobj}\) is a right adjoint, it preserves
monomorphisms.  Hence theorem~\ref{th:complift} applies to the pair
\((\cov{\KGvCGvobj}, \cov{\KdCGvobj})\) to prove
theorem~\ref{th:kcomprep}.

\begin{corollary}
Let \(\KdCGvobj\) be a \KdCGvobj and \(\KGvCGvobj\) be a \KGvCGvobj.
The composition
\[
  \Kcov{\KGvCGvobj} \Kcov{\KdCGvobj} \colon \Kdcat \to
  \KGvcat
\]
is representable with representing object
\(\mCo{\mKo{\ladj{\KdCGvobj}}}(\KGvCGvobj)\).
\end{corollary}

\begin{proof}
By theorem~\ref{th:complift}, the composition \(\Kcov{\KGvCGvobj}
\Kcov{\KdCGvobj}\) is isomorphic to \(\mKo{\cov{\KGvCGvobj}
  \Kcov{\KdCGvobj}}\).  Corollary~\ref{cor:repkv} says that
\(\cov{\KGvCGvobj}\Kcov{\KdCGvobj}\) is representable (in the
traditional sense) by \(\mCo{\mKo{\ladj{\KdCGvobj}}}(\KGvCGvobj)\) and
hence, by definition, \(\Kcov{\KGvCGvobj}
\Kcov{\KdCGvobj}\) is representable with representing object
\(\mCo{\mKo{\ladj{\KdCGvobj}}}(\KGvCGvobj)\).
\end{proof}

By putting \(\dcat = \Gvcat\) we have a pairing \(\KGvCGvcat \times
\KGvCGvcat \to \KGvCGvcat\) corresponding to composition of functors.
The next step in showing that \(\KGvCGvcat\) is monoidal is to prove
that the identity functor on \(\KGvcat\) is representable.

\begin{lemma}\label{disclift}
The discrete filtration functor, \(\disfunc \colon \Gvcat \to \KGvcat\),
lifts to a functor
\[
  \mCo{\disfunc} \colon \GvCGvcat \to \KGvCGvcat.
\]
For a \GvCGvobj[\GvCGvobj] and an \KGvobj[\KGvobj] there is an
isomorphism of \Gvalgs, natural in both \(\GvCGvobj\) and \(\KGvobj\),
\[
  \Hom{\KGvcat}{\mCo{\disfunc}(\GvCGvobj)}{\KGvobj} \cong
  \Hom{\Gvcat}{\GvCGvobj}{\kfofunc(\KGvobj)}.
\]
\end{lemma}

\begin{proof}
There is an obvious extension
\[
  \mZo{\disfunc} \colon \ZGvcat \to \ZKGvcat
\]
which is left adjoint to the corresponding extension of the forgetful
functor
\[
  \mZo{\kfofunc} \colon \ZKGvcat \to \ZGvcat.
\]
As \(\mZo{\disfunc}\) is a left adjoint, it preserves \co{product}s and
hence takes \co{algebra} objects to \co{algebra} objects.  It
therefore extends to \(\mCo{\disfunc} \colon \GvCGvcat \to \KGvCGvcat\).

Let \(\GvCGvobj\) be a \GvCGvobj and \(\KGvobj\) an \KGvobj.  The underlying
\Zsobj of \(\Hom{\KGvcat}{\mCo{\disfunc}(\GvCGvobj)}{\KGvobj}\) is
\begin{align*}
\isetelt \mapsto &
\Hom{\KGvcat}{\sabs{\mCo{\disfunc}(\GvCGvobj)}(\isetelt)}{\KGvobj} \\
&= \Hom{\KGvcat}{\mZo{\disfunc}(\sabs{\GvCGvobj})(\isetelt)}{\KGvobj} \\
&= \Hom{\KGvcat}{\disfunc(\sabs{\GvCGvobj}(\isetelt))}{\KGvobj} \\
&\cong \Hom{\Gvcat}{\sabs{\GvCGvobj}(\isetelt)}{\kfofunc(\KGvobj)}.
\end{align*}
Thus there is a bijection from the underlying \Zsobj of
\(\Hom{\KGvcat}{\mCo{\disfunc}(\GvCGvobj)}{\KGvobj}\) to that of
\(\Hom{\Gvcat}{\GvCGvobj}{\kfofunc(\KGvobj)}\).  It is clear that this is
an isomorphism of \Gvalgs.
\end{proof}

\begin{proposition}
Let \(I\) be the unit of the Tall\hyp{}Wraith monoidal structure on
\(\GvCGvcat\).  The \KGvCGvobj, \(\mCo{\disfunc}(I)\), represents the
identity functor on \(\KGvcat\).
\end{proposition}

Recall that as the functor \(\GvCGvcat \to \covrep(\Gvcat, \Gvcat)\)
is strong monoidal, \(I\) represents the identity functor on
\(\Gvcat\).

\begin{proof}
Let \(\KGvobj\) be an \KGvobj.  Using lemma~\ref{disclift}, we have
the following isomorphisms of \Gvalgs
\[
\Hom{\KGvcat}{\mCo{\disfunc}(I)}{\KGvobj} \cong
\Hom{\Gvcat}{I}{\kfofunc(\KGvobj)} \cong
  \kfofunc(\KGvobj).
\]
Thus \(\mCo{\disfunc}(I)\) represents the forgetful
functor, \(\kfofunc \colon \KGvcat \to \Gvcat\).  We need to show
that \(\mKo{\kfofunc} \colon \KGvcat \to \KGvcat\) is isomorphic to the
identity.

Let us write out the stages in defining \(\mKo{\kfofunc}(\KGvobj)\).  It is
\begin{align*}
\KGvobj = (\sabs{\KGvobj} \to \KGvobj_\lambda) &\xmapsto{\kflfunc}
(\KGvobj \to \disfunc(\KGvobj_\lambda)) \\
&\xmapsto{\kiqfunc} (\KGvobj \to \disfunc(\KGvobj_\lambda)) \\
&\xmapsto{\qipfunc} (\KGvobj \to \disfunc(\KGvobj_\lambda)) \\
&\xmapsto{\mPo{\forfunc}} (\sabs{\KGvobj} \to \KGvobj_\lambda) \\
&\xmapsto{\redfunc} (\sabs{\KGvobj} \to \KGvobj_\lambda) \\
&\xmapsto{\comfunc} (\sabs{\KGvobj} \to \KGvobj_\lambda) = \KGvobj.
\end{align*}
Hence \(\mKo{\kfofunc}\) is the identity functor on \(\KGvcat\) and thus the
identity functor is representable.
\end{proof}

The final step in showing that we have a monoidal structure on
\(\KGvcat\) is to show that the pairing and the unit satisfy the
various coherences.  These are automatic as these coherences are
satisfied in the category of functors from \(\KGvcat\) to itself.

Let us gather together all the main results into one theorem.  By
working with the appropriate opposite categories we gain one
additional result -- note that this involves a lot of switching, even
of the direction of the filtration; recall that we use
left\hyp{}pointing arrows for the categories of inductively filtred
objects.  By another obvious alteration, we can also have two graded
varieties of algebras involved, even with different grading sets.

\begin{theorem}
Let \(\Gvcat\) be a variety of graded algebras.
\begin{enumerate}
\item
\KGvCGvcatu, \(\KGvCGvcat\), has a monoidal structure with pairing
\[
  (\KGvCGvobj_1, \KGvCGvobj_2) \mapsto
  \mCo{\mKo{\ladj{\KGvCGvobj_2}}}(\KGvCGvobj_1)
\]
and unit \(\mCo{\disfunc}(I)\).

The functor \(\KGvCGvcat \to \covfun(\KGvcat, \KGvcat)\), \(\KGvCGvobj
\mapsto \Kcov{\KGvCGvobj}\), is strong monoidal.

\item Let \(\dcat\) be a complete, \co{complete}
extremally \co{well}\hyp{}powered, (extremal epi, mono) category.
There is a pairing
\[
  \KGvCGvcat \times \KdCGvcat \to \KdCGvcat
\]
compatible with the monoidal structure on \(\KGvCGvcat\) and with
composition of functors.

\item Let \(\dcat\) be a complete, \co{complete}, extremally
well\hyp{}powered, (epi, extremal mono)
category.  There is a pairing
\[
  \mOo{(\KGvCGvcat)} \times \LdGvcat \to \LdGvcat
\]
compatible with the monoidal structure on \(\KGvCGvcat\) and with
composition of functors.

\item Let \(\Gwcat\) be a graded variety of algebras, possibly with a
different grading set to that of \(\Gvcat\).  There is a pairing
\[
  \KGvCGwcat \times \KGvCGvcat \to \KGvCGwcat
\]
compatible with the monoidal structure on \(\KGvCGvcat\) and with
composition of functors. \noproof
\end{enumerate}
\end{theorem}

Again, we shall use the notation \(- \twprod -\) to denote the pairing
in each case.

We can deal with the variance shift on the middle pairing in the usual
manner.

\begin{proposition}
There is a pairing
\[
  \KGvCGvcat \times \LdGvcat \to \LdGvcat, \qquad (\KGvCGvobj,
  \LdGvobj) \mapsto \KGvCGvobj \circledast \LdGvobj
\]
compatible with the monoidal structure on \(\KGvCGvcat\)
and a natural isomorphism
\[
  \Hom{\LdGvcat}{\KGvCGvobj \circledast \LdGvobj}{\LdGvobj'} \cong
  \Hom{\LdGvcat}{\LdGvobj}{\KGvCGvobj \twprod \LdGvobj'}
\]
for \(\KGvCGvobj\) in \(\KGvCGvcat\) and \(\LdGvobj, \LdGvobj'\) in
\(\LdGvcat\). \noproof
\end{proposition}

We have the obvious definition.

\begin{defn}
Let \(\Gvcat\) be a variety of graded algebras.  A
\emph{Tall\hyp{}Wraith \(\KGvcat\)\hyp{}monoid} is a monoid in
\(\KGvCGvcat\).
\end{defn}

The remarks following the definition of a Tall\hyp{}Wraith
\(\vcat\)\hyp{}monoid regarding module and \co{module} objects clearly
apply here as well.

%% file: hopf.algtop.tex
\section{Tall\hyp{}Wraith Monoids in Algebraic Topology}
\label{sec:algtop}

In this section we come to the heart of the matter: operations on
cohomology theories.  For a suitable cohomology theory we shall show
that the set of unstable operations on that theory has the structure
of a Tall\hyp{}Wraith monoid.  This encodes all of the information
available except for the suspension isomorphism.  We shall also show
how to interpret the associated \emph{enriched Hopf ring} of the
\co{operations}.

\subsection{Operations}

We shall start by considering operations.  We shall do this in a
general setting.  Let us start by listing the ingredients and the
properties that they have to satisfy.

We have the following.
\begin{enumerate}
\item A category \(\ecat\) which is closed under finite products.

\item A graded variety of algebras, \(\Gvcat\), with indexing set
\(\iset\).

\item A filtration functor \(\eflfunc \colon \ecat \to \Iecat\).

\item A \eGvobj[\eGvobj].
\end{enumerate}

Before stating the conditions let us introduce some notation.  Using
the above data we construct a contravariant functor \(\Bcon{\eGvobj}
\colon \ecat \to \KGvcat\) as the composition
\[
  \ecat \xrightarrow{\eflfunc} \Iecat
  \xrightarrow{\Pcon{\eGvobj}} \PGvcat \xrightarrow{\redfunc}
  \QGvcat \xrightarrow{\comfunc} \KGvcat.
\]

The \eGvobj[\eGvobj] has underlying \Zeobj[\sabs{\eGvobj}] and so for
each \(\isetelt \in \iset\) we have an
\eobj[\sabs{\eGvobj}(\isetelt)].  
As \(\ecat\) is closed under finite
products, for each finite subset \(\iset_0 \subseteq \iset\)
we can form the product \(\prod_{\isetelt \in \iset_0}
\sabs{\eGvobj}(\isetelt)\). 
As \(\Bcon{\eGvobj}\) is a contravariant
functor, for \(\isetelt_0 \in \iset_0\) it takes the projection \emor
\(\prod_{\isetelt \in \iset_0}
\sabs{\eGvobj}(\isetelt)\to \sabs{\eGvobj}(\isetelt_0)\) to
 a \KGvmor
\[
\Bcon{\eGvobj}(\sabs{\eGvobj}(\isetelt_0)) \to
\Bcon{\eGvobj}\big(\prod_{\isetelt \in \iset_0}
\sabs{\eGvobj}(\isetelt)\big) 
\]
and thus, as \(\KGvcat\) is \co{complete}, we have a \KGvmor
\[
 \coprod_{\isetelt \in \iset_0}
\Bcon{\eGvobj}(\sabs{\eGvobj}(\isetelt)) \to
\Bcon{\eGvobj}\big(\prod_{\isetelt \in \iset_0}
\sabs{\eGvobj}(\isetelt)\big).
\]

We can now state our conditions on the ingredients for our
construction.

\begin{enumerate}
\item For each finite subset \(\iset_0 \subseteq \iset\), the above
  \KGvmor
\[
 \coprod_{\isetelt \in \iset_0}
\Bcon{\eGvobj}(\sabs{\eGvobj}(\isetelt)) \to
\Bcon{\eGvobj}\big(\prod_{\isetelt \in \iset_0}
\sabs{\eGvobj}(\isetelt)\big)
\]
is an isomorphism.

\item The filtration functor \(\eflfunc \colon \ecat \to \Iecat\) has
  the following property.  The inductive filtration
  \(\eflfunc(\eobj)\) for an \eobj[\eobj] has a final family, say
  \(\{\eobj_\lambda\}\), such that each \(\eflfunc(\eobj_\lambda)\) is
  an iso\hyp{}filtration.
\end{enumerate}

The first condition obviously implies the following
result.

\begin{lemma}
The \ZKGvobj[\ZBcon{\eGvobj}\!(\abs{\eGvobj})] underlies a
\KGvCGvobj. \noproof
\end{lemma}

The full notation for this \KGvCGvobj ought to be
\[
\ZBcon{\eGvobj}(\eGvobj)
\]
but that rapidly becomes unwieldy.  As we have fixed our initial data,
including \(\eGvobj\), we shall write \(\ZBcon{\eGvobj}(\eGvobj)\) as
just \(\KGvCGvTtobj\).  We shall also write \(\sabs{\KGvCGvTtobj}\)
for \(\ZBcon{\eGvobj}(\sabs{\eGvobj})\).

\begin{proposition}
\label{prop:elift}
Let \(\eobj\) be an \eobj.  The \Zsmoralt
\[
  \Hom{\ecat}{\eobj}{\sabs{\eGvobj}} \to
  \Hom{\KGvcat}{\ZBcon{\eGvobj}(\sabs{\eGvobj})}{%
    \Bcon{\eGvobj}(\eobj)} =
  \Hom{\KGvcat}{\sabs{\KGvCGvTtobj}}{%
    \Bcon{\eGvobj}(\eobj)}
\]
induced by \(\Bcon{\eGvobj}\) defines a morphism of \KGvobjs
\[
  \Bcon{\eGvobj}(\eobj) \to
  \Kcov{\KGvCGvTtobj}
  \big(\Bcon{\eGvobj}(\eobj)\big)
\]
which is natural in \(\eobj\).
\end{proposition}

The last expression could probably do with some explanation.  The
\KGvCGvobj[\KGvCGvTtobj] represents a covariant functor
\(\cov{\KGvCGvTtobj} \colon \KGvcat \to \Gvcat\).  By the first
construction of definition~\ref{def:isolift} we can lift this to a
covariant functor \(\KGvcat \to \KGvcat\) which, using the notation of
definition~\ref{def:isolift}, we write as \(\Kcov{\KGvCGvTtobj}\).

\begin{proof}
For \eobjs[\eobj_1 \text{ and } \eobj_2], the functor
\(\Bcon{\eGvobj}\) defines a morphism of sets, natural in \(\eobj_1\)
and \(\eobj_2\),
\[
  \Hom{\ecat}{\eobj_1}{\eobj_2} \to
  \Hom{\KGvcat}{\Bcon{\eGvobj}(\eobj_2)}{\Bcon{\eGvobj}(\eobj_1)}.
\]
For \(\isetelt \in \iset\) we therefore have a morphism of sets
\[
  \Hom{\ecat}{\eobj}{\sabs{\eGvobj}(\isetelt)} \to
  \Hom{\KGvcat}{\Bcon{\eGvobj}(\sabs{\eGvobj}(\isetelt))}{%
    \Bcon{\eGvobj}(\eobj)}
\]
and thus a \Zsmoralt
\[
  \Hom{\ecat}{\eobj}{\sabs{\eGvobj}} \to
  \Hom{\KGvcat}{\ZBcon{\eGvobj}(\sabs{\eGvobj})}{%
    \Bcon{\eGvobj}(\eobj)} =
  \Hom{\KGvcat}{\sabs{\KGvCGvTtobj}}{\Bcon{\eGvobj}(\eobj)}
\]
which can be rewritten as
\[
  \con{\sabs{\eGvobj}}(\eobj) \to
  \cov{\sabs{\KGvCGvTtobj}}
  \big(\Bcon{\eGvobj}(\eobj)\big).
\]

Both sides underlie \Gvalgs.  The structure on the left comes from the
occurence of \(\eGvobj\) as this is a \eGvobj.  That on the right
comes from the occurence of \(\sabs{\KGvCGvTtobj}\) as this underlies
a \KGvCGvobj.  This structure in turn comes from the fact that
\(\sabs{\KGvCGvTtobj} = \ZBcon{\eGvobj}(\sabs{\eGvobj})\) and
\(\sabs{\eGvobj}\) underlies a \eGvobj.  Both structures therefore
come from the same source and so the above morphism of \Zsobjs lifts
to a morphism of \Gvalgs
\begin{equation}
\label{eq:egvmor}
  \con{\eGvobj}(\eobj) \to
  \cov{\KGvCGvTtobj}
  \big(\Bcon{\eGvobj}(\eobj)\big)
\end{equation}
which is natural in \(\eobj\).

Both sides have filtrations and we wish to compare these.  The
filtration on the left comes from applying the functor \(\eflfunc\) to
\(\eobj\).  That on the right comes from applying the canonical
filtration functor, \(\qipfunc\kiqfunc\kflfunc \colon \KGvcat \to
\PKGvcat\), to \(\Bcon{\eGvobj}(\eobj)\).  We shall show that the
above morphism of \Gvalgs lifts to a morphism of \QGvalgs
\begin{equation}
\label{eq:redcomp}
  \redfunc\Pcon{\eGvobj}\eflfunc(\eobj) \to
  \redfunc\Pcov{(\KGvCGvTtobj)}
  \qipfunc\kiqfunc\kflfunc
  \big(\Bcon{\eGvobj}(\eobj)\big).
\end{equation}
The reduced projective filtrations on both sides are constructed from
the inductive filtration \(\eflfunc(\eobj)\) on \(\eobj\).
  Initial
subclasses can be explicitly constructed as follows.  We shall start
with the reduced projective filtration on the left hand side of
\eqref{eq:redcomp}.  Let \(\eelt \colon \eobj_{\eelt} \to \eobj\) be
in \(\eflfunc(\eobj)\).  Applying the contravariant functor
\(\con{\eGvobj}\) to \(\eelt\) results in a \Gvmor
\(\con{\eGvobj}(\eobj) \to \con{\eGvobj}(\eobj_{\eelt})\).  We choose
an (extremal epi, mono)\hyp{}factorisation of this with intervening
\Gvobj[\mQo{\con{\eGvobj}(\eobj_{\eelt})}].  We therefore have an
extremal epimorphism \(\con{\eGvobj}(\eobj) \to
\mQo{\con{\eGvobj}(\eobj_{\eelt})}\).  Doing this for each \(\eelt\)
in \(\eflfunc(\eobj)\) defines the required initial subclass for
\(\redfunc\Pcon{\eGvobj}\eflfunc(\eobj)\).

Now let us consider the reduced projective filtration on the right
hand side.  The \KGvalg, \(\Bcon{\eGvobj}(\eobj)\), as a projectively
filtered \Gvalg, has initial subclass given by the \Gvmors
\[
  \sabs{\Bcon{\eGvobj}(\eobj)} \to \mQo{\con{\eGvobj}(\eobj_{\eelt})}
\]
for \(\eelt\) in \(\eflfunc(\eobj)\).  An initial subclass for the
canonical filtration on \(\Bcon{\eGvobj}(\eobj)\) is therefore given
by the \KGvmors
\[
  \Bcon{\eGvobj}(\eobj) \to
  \disfunc\left(\mQo{\con{\eGvobj}(\eobj_{\eelt})}\right).
\]
Therefore an initial subclass for the reduced projective filtration on
the right hand side of \eqref{eq:redcomp} is given by the \Gvmors
\[
  \cov{\KGvCGvTtobj}
  \big(\Bcon{\eGvobj}(\eobj)\big) \to
  \mQo{\cov{\KGvCGvTtobj} \disfunc
    \big(\mQo{\con{\eGvobj}(\eobj_{\eelt})}\big)}.
\]
  
Let \(\eelt\) be in \(\eflfunc(\eobj)\).  We claim that
the following is a
commutative diagram in \(\Gvcat\).
\[
  \xymatrix{
    \con{\eGvobj}(\eobj) \ar[d] \ar[r] &
      \cov{\KGvCGvTtobj}
  \big(\Bcon{\eGvobj}(\eobj)\big) \ar[d] \\
    \mQo{\con{\eGvobj}(\eobj_{\eelt})} \ar[dd] &
  \mQo{\cov{\KGvCGvTtobj} \disfunc
    \big(\mQo{\con{\eGvobj}(\eobj_{\eelt})}\big)} \ar[d] \\
    &   \cov{\KGvCGvTtobj} \disfunc
    \big(\mQo{\con{\eGvobj}(\eobj_{\eelt})}\big) \ar[d] \\
    \con{\eGvobj}(\eobj_{\eelt}) \ar[r] &
     \cov{\KGvCGvTtobj}
    \big(\Bcon{\eGvobj}(\eobj_{\eelt})\big)
  }
\]
The lower horizontal morphism comes from applying \eqref{eq:egvmor}
with \(\eobj_{\eelt}\) in place of \(\eobj\).  The bottom vertical
morphism on the right arises as follows.  From the definition of
\(\mQo{\con{\eGvobj}(\eobj_{\eelt})}\) there is a morphism
\[
  \mQo{\con{\eGvobj}(\eobj_{\eelt})} \to \con{\eGvobj}(\eobj_{\eelt}).
\]
From the definition of \(\Bcon{\eGvobj}(\eobj_{\eelt})\) there is a morphism
\[
  \con{\eGvobj}(\eobj_{\eelt}) \to \sabs{\Bcon{\eGvobj}(\eobj_{\eelt})}
\]
whence we have a morphism
\[
  \disfunc \big(  \mQo{\con{\eGvobj}(\eobj_{\eelt})}\big) \to
  \disfunc\big( \con{\eGvobj}(\eobj_{\eelt})\big) \to
  \Bcon{\eGvobj}(\eobj_{\eelt})
\]
to which we apply \(\cov{\KGvCGvTtobj}\).

It is routine to check that the composition on the right hand side is
\(\cov{\KGvCGvTtobj}\big(\Bcon{\eGvobj}(\eelt)\big)\)
and thus the above diagram commutes by
naturality of \eqref{eq:egvmor}.

We claim that our assumption on \(\eflfunc\) ensures that the bottom
vertical morphism on the right hand side is a monomorphism.  As
\(\con{\eGvobj}\) is representable, it is one of a mutually right
adjoint pair and thus takes epi\hyp{}sinks to mono\hyp{}sources.
Hence \(\con{\eGvobj}(\eobj_{\eelt})\) is a mono\hyp{}source for
\(\Pcon{\eGvobj}\eflfunc(\eobj_{\eelt})\).  This obviously remains
true after applying the reduction functor and hence the natural
morphism \(\con{\eGvobj}(\eobj_{\eelt}) \to
\sabs{\Bcon{\eGvobj}(\eobj_{\eelt})}\) is a monomorphism.  Since
\(\mQo{\con{\eGvobj}(\eobj_{\eelt})} \to
\con{\eGvobj}(\eobj_{\eelt})\) is a monomorphism by construction we
therefore have a monomorphism \(\mQo{\con{\eGvobj}(\eobj_{\eelt})}
\to\sabs{\Bcon{\eGvobj}(\eobj_{\eelt})}\).  This morphism underlies
the \KGvmor \(\disfunc \big(\mQo{\con{\eGvobj}(\eobj_{\eelt})}\big)
\to \Bcon{\eGvobj}(\eobj_{\eelt})\) which is therefore a monomorphism
as the forgetful functor is faithful.  Since
\(\cov{\KGvCGvTtobj}\) is a right adjoint it
takes monomorphisms to monomorphisms.

Thus the composition of the lower two vertical morphisms on the right
hand side is a monomorphism and so as \Gvcat is an (extremal epi,
mono) category we have a horizontal morphism 
\[
    \mQo{\con{\eGvobj}(\eobj_{\eelt})} \to
  \mQo{\cov{\KGvCGvTtobj} \disfunc
    \big(\mQo{\con{\eGvobj}(\eobj_{\eelt})}\big)}
\]
fitting into the diagram.

This shows that the \Gvmor in \eqref{eq:egvmor} lifts to a \QGvmor
\[
  \redfunc\Pcon{\eGvobj}\eflfunc(\eobj) \to
  \redfunc\Pcov{\KGvCGvTtobj}
  \qipfunc\kiqfunc\kflfunc
  \big(\Bcon{\eGvobj}(\eobj)\big)
\]
as required.  Applying \(\comfunc\) to both sides results in a \KGvmor
\[
  \comfunc\redfunc\Pcon{\eGvobj}\eflfunc(\eobj) \to
  \comfunc\redfunc\Pcov{\KGvCGvTtobj}
  \qipfunc\kiqfunc\kflfunc
  \big(\Bcon{\eGvobj}(\eobj)\big).
\]

By definition, this is a morphism of \KGvalgs
\[
  \Bcon{\eGvobj}(\eobj) \to
  \Kcov{\KGvCGvTtobj}
  \big(\Bcon{\eGvobj}(\eobj)\big)
\]
as required.  It is obvious from its construction that it is natural
in \(\eobj\).
\end{proof}

\begin{corollary}
\label{cor:tw}
The \KGvCGvobj, \(\KGvCGvTtobj\), is the underlying
\KGvCGvobjalt of a Tall\hyp{}Wraith \(\KGvcat\)\hyp{}monoid.

For an \eobj[\eobj], \(\Bcon{\eGvobj}(\eobj)\) is a module for this
monoid.
\end{corollary}

\begin{proof}
Using the adjunction \(\mKo{\ladj{\KdCGvobj}} \adjoint
\Kcov{\KdCGvobj}\) from corollary~\ref{cor:liftadj}, the morphism from
proposition~\ref{prop:elift} defines a morphism of \KGvalgs
\begin{equation}
\label{eq:elift}
  \Kladj{\KGvCGvTtobj}
  \big(\Bcon{\eGvobj}(\eobj)\big) \to
  \Bcon{\eGvobj}(\eobj)
\end{equation}
By applying this to the underlying \eobjs of \(\eGvobj\) we obtain a
morphism of \ZKGvalgs
\[
  \Kladj{\KGvCGvTtobj}
  (\sabs{\KGvCGvTtobj}) \to
  \sabs{\KGvCGvTtobj}.
\]
As left adjoints preserve \co{products} this lifts to a morphism of
\KGvCGvobjs
\[
  \Kladj{\KGvCGvTtobj}
  (\KGvCGvTtobj) \to
  \KGvCGvTtobj.
\]
From the definition of \(\twprod\) we can rewrite this as
\[
\KGvCGvTtobj \twprod
\KGvCGvTtobj \to
  \KGvCGvTtobj.
\]  

That this is associative comes from its origins.  This morphism
started life as the morphism of \Zsobjs
\[
  \Hom{\ecat}{\sabs{\eGvobj}}{\sabs{\eGvobj}} \to
  \Hom{\Gvcat}{\con{\eGvobj}(\sabs{\eGvobj})}{%
    \con{\eGvobj}(\sabs{\eGvobj})}
  \subseteq \Hom{\Zscat}{
    \Hom{\ecat}{\sabs{\eGvobj}}{\sabs{\eGvobj}}}{
    \Hom{\ecat}{\sabs{\eGvobj}}{\sabs{\eGvobj}}}
\]
which is adjoint to the morphism
\[
  \Hom{\ecat}{\sabs{\eGvobj}}{\sabs{\eGvobj}} \times
  \Hom{\ecat}{\sabs{\eGvobj}}{\sabs{\eGvobj}} \to
  \Hom{\ecat}{\sabs{\eGvobj}}{\sabs{\eGvobj}}
\]
and this is clearly associative.

For the unit, we observe that the identity morphism on
\(\sabs{\eGvobj}\) defines a morphism of \ZZsobjs
\begin{equation}
\label{eq:idtw}
  \diafunc(\{*\}) \to \Hom{\ecat}{\sabs{\eGvobj}}{\sabs{\eGvobj}}
\end{equation}
where \(\diafunc \colon \scat \to \ZZscat\) is the diagonal functor
defined on objects by
\[
  \diafunc(\sobj) = \Bigg(\isetelt \mapsto \big( \isetelt'
  \mapsto \begin{cases} \sobj &\text{if } \isetelt = \isetelt' \\
  \emptyset &\text{otherwise} 
  \end{cases}
  \big)\Bigg)
\]
and in the obvious way on morphisms.

The right hand side of \eqref{eq:idtw} is the underlying
\ZZsobj of an \ZGvobj and so we have
an adjoint \ZGvmor
\begin{equation}
\label{eq:idtwgv}
  \free[mZo]{Gv} \diafunc(\{*\}) \to \con{\eGvobj}(\sabs{\eGvobj}).
\end{equation}

The \free{Gv} on \(\diafunc(\{*\})\) is (up to canonical isomorphism)
the underlying \ZGvalg of the unit, \(I\), of the Tall\hyp{}Wraith
monoidal structure on \(\GvCGvcat\).  It is straightforward to show
that the morphism in \eqref{eq:idtwgv} lifts to a morphism of
\GvCGvobjs
\[
  I \to \con{\eGvobj}(\eGvobj).
\]

The \GvCGvobjalt[\con{\eGvobj}(\eGvobj)] is the underlying \GvCGvobjalt
of the \PGvCGvobjalt[\Pcon{\eGvobj}\eflfunc(\eGvobj)], whence we obtain a
morphism of \PGvCGvobjalts
\[
  \mCo{\disfunc}(I) \to \Pcon{\eGvobj}\eflfunc(\eGvobj).
\]
Applying \(\redfunc\) and \(\comfunc\) to this does not change the
source as that is already an \KGvobj.  We therefore have a morphism of
\KGvCGvobjalts
\[
  \mCo{\disfunc}(I) \to \ZBcon{\eGvobj}(\eGvobj) = \KGvCGvTtobj.
\]

As this originated from the inclusion of the identity morphism on
\(\sabs{\eGvobj}\), it is clear that this is a unit for the product
defined above.  We therefore have the structure of a Tall\hyp{}Wraith
\(\KGvcat\)\hyp{}monoid.

Returning to \eqref{eq:elift} we see that it defines a
\KGvmoralt
\[
\KGvCGvTtobj \twprod
  \big(\Bcon{\eGvobj}(\eobj)\big) \to
  \Bcon{\eGvobj}(\eobj)
\]
and it is obvious that this is compatible with the monoidal 
structure, thus making \(\Bcon{\eGvobj}(\eobj)\)
a module for the Tall\hyp{}Wraith \(\KGvcat\)\hyp{}monoid.
\end{proof}

Applying this to a cohomology functor we obtain the following result,
giving theorems~\ref{th:twcoh} and~\ref{th:twcohmod} of the
introduction.

\begin{corollary}
Let \(\Efunc\) be a graded, multiplicative, commutative cohomology
theory with representing objects \(\Erep\).  Let \(\Evar\) be the
variety of graded, commutative, \Ealgs.  Suppose that
\(\Ehom(\Erep[k])\) is free as an \(\Ecoef\)\hyp{}module for all
\(k\).  Then \(\Efunc(\Erep)\) is a Tall\hyp{}Wraith
\(\KEvar\)\hyp{}monoid and for a space \(X\), \(\BEfunc(X)\) is
a module for this monoid.
\end{corollary}

\begin{proof}
In this case we have the following ingredients.
\begin{enumerate}
\item \(\ecat\) is the category of spaces of the homotopy type of a
  \(\m{C W}\)\hyp{}complex, which is closed under finite products.

\item \(\Gvcat\) is \(\Evar\), the variety of graded, commutative
\Ealgs, with indexing set \Z.

\item The filtration functor, \(\eflfunc \colon \ecat \to \Iecat\), is
 the \emph{ind\hyp{}finite} filtration functor defined by the
 subcategory of finite \(\m{C W}\)\hyp{}complexes; see
 example~\ref{ex:profin} of~\ref{ex:filtrations} for the projective
 version.

\item The \eGvobj[\Erep] is provided by Brown's representability
 theorem.
\end{enumerate}

The freeness condition on the cohomology theory in the statement of this
corollary ensures that the cohomology theory satisfies the (completed)
K{\"u}nneth formula for its own representing spaces and thus satisfies
the required ``product to \co{product}'' condition.  Incidentally,
this condition also ensures that the cohomology operations are already
\KGvobjs and so the \(\comfunc\) functor simply regards them in
\(\KGvcat\) rather than \(\QGvcat\).

The filtration functor satisfies its required condition since every
filtration obtained by applying the filtration functor has a final
class of finite objects and it is clear that such objects have the
discrete filtration.

We can therefore apply the work of this section to obtain the desired
result.
\end{proof}

\subsection{\Co{operations}}
\label{subsec:co-ops}

Now we turn to \co{operations}; that is, to the Hopf ring associated
to a suitable cohomology theory.  In fact, we do not need to assume
that the homology theory used to produce the Hopf ring is that
associated to the original cohomology theory.  Moreover, by using the
standard categorical trick of swapping a category for its dual we can
also consider the set of operations from one suitable cohomology
theory to another.

In the general setting, we have the following ingredients.
\begin{enumerate}
\item A category \(\ecat\) which is closed under finite products.

\item A complete, \co{complete}, extremally well\hyp{}powered,
(epi, extremal mono) category, \(\dcat\), with the property that
epimorphisms and filtered \co{limits} commute with finite products.

\item A graded variety of algebras, \(\Gvcat\).

\item A filtration functor \(\eflfunc \colon \ecat \to \Iecat\).

\item A covariant functor \(\homfunc \colon \ecat \to \dcat\).

\item A \eGvobj[\eGvobj].
\end{enumerate}

The conditions that we require are as follows.
\begin{enumerate}
\item For each finite subset \(\iset_0 \subseteq \iset\), the natural
\KGvmor
\[
 \coprod_{\isetelt \in \iset_0}
\Bcon{\eGvobj}(\sabs{\eGvobj}(\isetelt)) \to
\Bcon{\eGvobj}\big(\prod_{\isetelt \in \iset_0}
\sabs{\eGvobj}(\isetelt)\big)
\]
is an isomorphism.
\item For each finite subset \(\iset_0 \subseteq \iset\), the natural
  \Ldmor
\[
\mBo{\homfunc}\big(\prod_{\isetelt \in \iset_0}
\sabs{\eGvobj}(\isetelt)\big) \to
 \prod_{\isetelt \in \iset_0}
\mBo{\homfunc}(\sabs{\eGvobj}(\isetelt)) 
\]
is an isomorphism; where \(\mBo{\homfunc} \colon \ecat \to \Ldcat\) is
the composition
\[
  \ecat \xrightarrow{\eflfunc} \Iecat
  \xrightarrow{\mIo{\,\homfunc\,}} \Idcat \xrightarrow{\redfunc}
  \Jdcat \xrightarrow{\comfunc} \Ldcat.
\]

\item The filtration functor \(\eflfunc \colon \ecat \to \Iecat\) has
  the following property.  The inductive filtration
  \(\eflfunc(\eobj)\) for an \eobj[\eobj] has a final family, say
  \(\{\eobj_\lambda\}\), such that each \(\eflfunc(\eobj_\lambda)\) is
  an iso\hyp{}filtration.
\end{enumerate}

The second condition obviously implies the following result.

\begin{lemma}
The \ZLdobj[\mZo{\mBo{\homfunc}}(\abs{\eGvobj})] underlies a \LdGvobj.
\noproof
\end{lemma}

The full notation for this \LdGvobj would be
\(\mZo{\mBo{\homfunc}}(\eGvobj)\).  We shall shorten this to
\(\LdGvobj\).

We therefore have a contravariant functor
\(\con{\LdGvobj} \colon \Ldcat \to
\Gvcat\).  Using the standard categorical trick, we apply the first
construction of definition~\ref{def:isolift} to this contravariant
functor to obtain a lift
\[
  \Kcon{\LdGvobj} \colon \Ldcat \to
  \KGvcat.
\]

Recall that \( \KGvCGvTtobj = \ZBcon{\eGvobj}(\eGvobj)\) is,
by corollary~\ref{cor:tw}, the underlying
\KGvCGvobjalt of a Tall\hyp{}Wraith \(\KGvcat\)\hyp{}monoid.
We have the following result by a similar argument to that in the
previous section.

\begin{proposition}
\label{prop: covtwmod}
The \LdGvobj[\LdGvobj]
is a module for the
Tall\hyp{}Wraith \(\KGvcat\)\hyp{}monoid,
\(\KGvCGvTtobj\).  For an \eobj[\eobj], the \KGvalg
\[
  \Kcon{\LdGvobj}(\mBo{\homfunc}(\eobj))
\]
is a module for \(\KGvCGvTtobj\) and the morphism of
\KGvalgs
\[
  \Bcon{\eGvobj}(\eobj) \to
  \Kcon{\LdGvobj}(\mBo{\homfunc}(\eobj))
\]
induced by the morphism of \Zsobjs
\[
  \Hom{\ecat}{\eobj}{\sabs{\eGvobj}} \to
  \Hom{\Ldcat}{\mBo{\homfunc}(\eobj)}{\mZo{\mBo{\homfunc}}(\sabs{\eGvobj})}
  = \Hom{\Ldcat}{\mBo{\homfunc}(\eobj)}{\sabs{\LdGvobj}}
\]
is a morphism of \(\KGvCGvTtobj\)\hyp{}modules. \noproof
\end{proposition}

Notice that we use the variance switch to turn the more obvious
\co{action} morphism
\[
  \LdGvobj \to \LdGvobj
  \twprod \KGvCGvTtobj
\]
into an action morphism
\[
  \KGvCGvTtobj \circledast
  \LdGvobj \to \LdGvobj.
\]

By replacing \(\dcat\) by \(\Odcat\) we can consider contravariant
functors \(\homfunc \colon \ecat \to \dcat\).  The above conditions on
\(\dcat\) are now conditions on \(\Odcat\) and so need to be replaced
by their duals to get the required conditions on \(\dcat\).  With the
obvious notation we have the following result.

\begin{proposition}
\label{prop:contwmod}
The \KdCGvobj[\KdCGvobj] is a module for the
Tall\hyp{}Wraith \(\KGvcat\)\hyp{}monoid,
\(\KGvCGvTtobj\).  For an \eobj[\eobj], the \KGvalg
\[
  \Kcov{\KdCGvobj}(\mBo{\homfunc}(\eobj))
\]
is a module for \(\KGvCGvTtobj\) and the morphism of
\KGvalgs
\[
  \Bcon{\eGvobj}(\eobj)\to
  \Kcov{\KdCGvobj}(\mBo{\homfunc}(\eobj))
\]
induced by the morphism of \Zsobjs
\[
  \Hom{\ecat}{\eobj}{\sabs{\eGvobj}} \to
  \Hom{\Kdcat}{\mZo{\mBo{\homfunc}}(\sabs{\eGvobj})}{\mBo{\homfunc}(\eobj)}
  = \Hom{\Kdcat}{\sabs{\KdCGvobj}}{\mBo{\homfunc}(\eobj)}
\]
is a morphism of \(\KGvCGvTtobj\)\hyp{}modules. \noproof
\end{proposition}

We now obtain theorem~\ref{th:twhopf} of the introduction by applying
the above results to suitable homology and cohomology functors.

Firstly, applying proposition~\ref{prop: covtwmod} to an appropriate
homology functor we obtain the following description of an enriched
Hopf ring.

\begin{proposition}
\label{prop:enhopf}
Let \(\Efunc\) and \(\Ffunc\) be graded, multiplicative,
commutative cohomology theories with representing objects
\(\Erep\) and \(\Frep\) respectively.  Let \(\Evar\) be the
variety of graded, commutative, \Ealgs.

Suppose that \(\Fhom(\Erep[k])\) is free as an
\(\Fcoef\)\hyp{}module for all \(k\) and that
\(\Ehom(\Erep[k])\) is free as an
\(\Ecoef\)\hyp{}module for all \(k\).  Then
\(\Fhom(\Erep)\) is a module for the Tall\hyp{}Wraith
\(\KEvar\)\hyp{}monoid \(\Efunc(\Erep)\).

For a space \(X\) with the property that \(\Fhom(X)\) is free as
an \Fmod, the obvious morphism
\[
  \Efunc(X) \to \Kcon{\big(
    \mZo{\mBo{\Fhom}}{(\Erep)}\big)}(\Fhom(X))
\]
is a morphism of \(\Efunc(\Erep)\)\hyp{}modules. \noproof
\end{proposition}

In the last part of this proposition our intermediate category is that
of \Fcoalgs.  The conditions on \(\Efunc\) and \(\Ffunc\) ensure that
each \(\Fhom(\Erep[k])\) is an \Fcoalg and thus that \(\Fhom(\Erep)\)
is an \Ealg object in \Fcoalgs.  The condition on the space \(X\) also
ensures that \(\Fhom(X)\) is an \Fcoalg.  As we
show in~\cite{assw2}, \Fcoalgcat
does have all the categorical properties necessary for our results to
apply.

Secondly, applying proposition~\ref{prop:contwmod} to an appropriate
cohomology functor we obtain the following description of the
structure of the operations from one cohomology theory to another.

\begin{proposition}
\label{prop:optoop}
Let \(\Efunc\) and \(\Ffunc\) be graded, multiplicative,
commutative cohomology theories with representing objects
\(\Erep\) and \(\Frep\) respectively.  Let \(\Evar\) be the
variety of graded, commutative, \Ealgs.

Suppose that \(\Fhom(\Erep[k])\) is free as an
\Fmod for all \(k\) and that
\(\Ehom(\Erep[k])\) is free as an
\Emod for all \(k\).  Then
\(\Ffunc(\Erep)\) is a module for the Tall\hyp{}Wraith
\(\KEvar\)\hyp{}monoid \(\Efunc(\Erep)\).

For a space \(X\), the obvious morphism
\[
  \Efunc(X) \to \Kcov{\big(
    \mZo{\mKo{\Ffunc}}{(\Erep)}\big)}(\Ffunc(X))
\]
is a morphism of \(\Efunc(\Erep)\)\hyp{}modules. \noproof
\end{proposition}

As is well\hyp{}known, the assumptions in
proposition~\ref{prop:enhopf} are sufficient to ensure that the object
\(\Fhom(\Erep)\) is a \emph{Hopf ring}; that is, an
\Ealg object in \Fcoalgcat;
see \cite{dcrwsw} and \cite{wsw00}.  In \cite{jbdjww} the authors
expanded this to the notion of an \emph{enriched Hopf ring}.  The
\emph{enriched} part corresponds to the action of the Tall\hyp{}Wraith
monoid \(\Efunc(\Erep)\).

It is not usual to view homology as topologised.  Indeed, the
inductive filtration does not really give these structures a topology.
Even though there is a notion of an inductive topology, this is not
how one should regard the inductive filtration on the homology groups
of a space.  Rather, it is a \emph{dual topology} in that it is the
structure required to induce a topology on the dual, or more generally
on the result of applying some other contravariant functor.  Thus we
have not altered the usual view of homology as being discrete.

This concept of an inductive filtration was not used explicitly in
\cite{jb4} and \cite{jbdjww}.  On the other hand, it is there under
the surface.  In our language, what the authors of \cite{jb4} and
\cite{jbdjww} do is to define an inductive filtration functor on the
category of modules over some ring.  They then show that this is
compatible with taking the usual filtration functor on the category of
topological spaces and so rather than having to lift the homology
functor to the filtered category they first apply the (unfiltered)
homology functor and follow it by the filtration functor on the
category of modules.

The Milnor short exact sequence for homology shows that the homology
of a space is always iso\hyp{}filtered and so the issue of passing to
the ``completion'' does not arise in homology.

\medskip

As noted above, in~\cite{assw2} we show that the category of
\co{commutative} \Fcoalgs has all the properties
necessary for our results to apply.  This means that one can ``do''
general algebra in this category and thus we can apply the results of
section~\ref{sec:genalg}.  From proposition~\ref{prop:freegrade} and
theorem~\ref{th:impidgrade} we deduce the existence of free \Gvalgobjs
for any variety of graded algebras.  In particular, if \(\Gvcat\) is
the category of commutative \Ealgs we recover the
result of Ravenel and Wilson, \cite[\S 1]{dcrwsw}, that free Hopf
rings exist.  If \(\Gvcat\) is the category of
\Emods we recover the result of Hunton and Turner,
\cite[\S 2]{jrhprt}, that free \co{algebraic} modules exist.

\medskip

In future work we expect to give descriptions of the Tall\hyp{}Wraith
monoids of unstable operations for various cohomology theories,
including the Morava \(K\)-theories.

%% file: hopf.bbl
\begin{thebibliography}{KPMS82}

\bibitem[Ber98]{gb}
George~M. Bergman.
\newblock {\em An invitation to general algebra and universal constructions}.
\newblock Henry Helson, Berkeley, CA, 1998.

\bibitem[BH96]{gbah}
George~M. Bergman and Adam~O. Hausknecht.
\newblock {\em Co-groups and co-rings in categories of associative rings},
  volume~45 of {\em Mathematical Surveys and Monographs}.
\newblock American Mathematical Society, Providence, RI, 1996.

\bibitem[BJW95]{jbdjww}
J.~Michael Boardman, David~Copeland Johnson, and W.~Stephen Wilson.
\newblock Unstable operations in generalized cohomology.
\newblock In {\em Handbook of algebraic topology}, pages 687--828.
  North-Holland, Amsterdam, 1995.

\bibitem[Boa95]{jb4}
J.~Michael Boardman.
\newblock Stable operations in generalized cohomology.
\newblock In {\em Handbook of algebraic topology}, pages 585--686.
  North-Holland, Amsterdam, 1995.

\bibitem[BW05]{jbbw}
James Borger and Ben Wieland.
\newblock Plethystic algebra.
\newblock {\em Adv. Math.}, 194(2):246--283, 2005.

\bibitem[CMS02]{dscmns}
Dena~S. Cowen~Morton and Neil Strickland.
\newblock The {H}opf rings for {$K{\rm O}$} and {$K{\rm U}$}.
\newblock {\em J. Pure Appl. Algebra}, 166(3):247--265, 2002.

\bibitem[Fre66]{pf2}
P.~Freyd.
\newblock Algebra valued functors in general and tensor products in particular.
\newblock {\em Colloq. Math.}, 14:89--106, 1966.

\bibitem[HS73]{hhgs}
Horst Herrlich and George~E. Strecker.
\newblock {\em Category theory: an introduction}.
\newblock Allyn and Bacon Series in Advanced Mathematics. Allyn and Bacon Inc.,
  Boston, Mass., 1973.

\bibitem[HT98]{jrhprt}
John~R. Hunton and Paul~R. Turner.
\newblock Coalgebraic algebra.
\newblock {\em J. Pure Appl. Algebra}, 129(3):297--313, 1998.

\bibitem[KPMS82]{wkjmmpis}
W.~K{\"u}hnel, M.~Pfender, J.~Meseguer, and I.~Sols.
\newblock Algebras with actions and automata.
\newblock {\em Internat. J. Math. Math. Sci.}, 5(1):61--85, 1982.

\bibitem[RW96]{dcrwsw96}
Douglas~C. Ravenel and W.~Stephen Wilson.
\newblock The {H}opf ring for {$P(n)$}.
\newblock {\em Canad. J. Math.}, 48(5):1044--1063, 1996.

\bibitem[RW77]{dcrwsw}
Douglas~C. Ravenel and W.~Stephen Wilson.
\newblock The {H}opf ring for complex cobordism.
\newblock {\em J. Pure Appl. Algebra}, 9(3):241--280, 1976/77.

\bibitem[SW06]{math.AT/0605471}
Andrew Stacey and Sarah Whitehouse.
\newblock {Stable and Unstable Operations in mod \(p\) Cohomology Theories},
  2006.

\bibitem[SW07a]{assw2}
Andrew Stacey and Sarah Whitehouse.
\newblock General algebra in a category of coalgebras, 2007.

\bibitem[SW07b]{assw3}
Andrew Stacey and Sarah Whitehouse.
\newblock Tall-{W}raith monoids, 2007.

\bibitem[TW70]{dtgw}
D.~O. Tall and G.~C. Wraith.
\newblock Representable functors and operations on rings.
\newblock {\em Proc. London Math. Soc. (3)}, 20:619--643, 1970.

\bibitem[Wil84]{wsw84}
W.~Stephen Wilson.
\newblock The {H}opf ring for {M}orava {$K$}-theory.
\newblock {\em Publ. Res. Inst. Math. Sci.}, 20(5):1025--1036, 1984.

\bibitem[Wil00]{wsw00}
W.~Stephen Wilson.
\newblock Hopf rings in algebraic topology.
\newblock {\em Expo. Math.}, 18(5):369--388, 2000.

\end{thebibliography}
